\documentclass[11pt,leqno]{amsart}
\usepackage{amssymb}
\usepackage{mathrsfs}
\usepackage{color}
\usepackage{soul}

\newcommand\N{{\mathbb N}}
\newcommand\R{{\mathbb R}}

\newcommand\PPP{{\bf P}}

\newcommand\Ee{{\mathbb E}}

\def\AA{{\mathcal A}}
\def\BB{{\mathcal B}}
\def\CC{{\mathcal C}}

\def\HH{{\mathcal H}}
\def\II{{\mathcal I}}
\def\JJ{{\mathcal J}}
\def\KK{{\mathcal K}}

\def\MM{{\mathcal M}}
\def\NN{{\mathcal N}}
\def\OO{{\mathcal O}}
\def\PP{{\mathcal P}}

\def\SS{{\mathcal S}}

\def\VV{{\mathcal V}}
\def\WW{{\mathcal W}}
\def\XX{{\mathcal X}}
\def\YY{{\mathcal Y}}
\def\ZZ{{\mathcal Z}}

\def\PPP{{\mathbf P}}

\def\BBB{{\mathscr{B}}}
\def\LLL{{\mathscr{L}}}
\def\MMM{{\mathscr{M}}}

\def\Ps{{\bf P}_{\! sym}}

\def\PEj{{\bf P}(E^j)}

\def\PPE{{\bf P} ({\bf P} (E))}

\def\WWs{{\WW_{H^{-s}}}}

\DeclareMathOperator{\diver}{div}
\DeclareMathOperator{\Supp}{Supp}

\def\eps{{\varepsilon}}

\newcommand{\wto}{\rightharpoonup}

\def\SN{\mathfrak{S}_N}

\newtheorem{theo}{Theorem}
\newtheorem{prop}[theo]{Proposition}
\newtheorem{lem}[theo]{Lemma}
\newtheorem{cor}[theo]{Corollary}
\newtheorem{rem}[theo]{Remark}

\newtheorem{defin}[theo]{Definition}

\newcommand{\beqn}{\begin{equation}}
\newcommand{\eeqn}{\end{equation}}
\newcommand{\bear}{\begin{eqnarray}}
\newcommand{\eear}{\end{eqnarray}}
\newcommand{\bean}{\begin{eqnarray*}}
\newcommand{\eean}{\end{eqnarray*}}

\newcommand{\ds}{\displaystyle}

\newcommand{\indiq}{{\bf 1}}

\newcommand{\Black}{\color{black}}

\def\signmh{\bigskip \begin{center} {\sc Maxime Hauray\par\vspace{3mm}
Universit\'e d'Aix-Marseille \par
LATP, UMR CNRS 7353\par
Centre de Math\'ematiques et Informatique (CMI) \par
39, rue F. Joliot Curie 13453 Marseille Cedex 13 \par
FRANCE\par\vspace{3mm}
e-mail:} \tt{maxime.hauray@univ-amu.fr} \end{center}}

\def\signsm{\bigskip \begin{center} {\sc St\'ephane Mischler\par\vspace{3mm}
Universit\'e Paris-Dauphine \& IUF \par
CEREMADE, UMR CNRS 7534\par
Place du Mar\'echal de Lattre de Tassigny
75775 Paris Cedex 16\par
FRANCE\par\vspace{3mm}
e-mail:} \tt{mischler@ceremade.dauphine.fr} \end{center}}

\begin{document}

\title[Chaos]{On Kac's chaos and related problems}

\author{M. Hauray}
\author{S. Mischler}

\begin{abstract}
This paper is devoted to establish quantitative and qualitative estimates
related to the notion of  chaos as firstly formulated by M. Kac \cite{Kac1956}
in his study of  mean-field limit for systems of $N$ undistinguishable particles
as $N\to\infty$.

First, we quantitatively liken three usual measures of {\it Kac's chaos}, some
involving all the $N$ variables, others 
 involving a finite fixed number of variables. 

Next, we define the notion of {\it entropy chaos } and  {\it Fisher information
chaos } in a similar
way as defined by Carlen et al \cite{CCLLV}. We show that  {\it Fisher
information chaos } is 
stronger than  {\it entropy chaos}, which in turn is stronger than  {\it Kac's
chaos}. 
We also establish 
that  {\it Kac's chaos} plus  {\it Fisher
information bound } implies  {\it entropy chaos}.  
 
 We then extend our analysis 
 to
the framework of 
probability measures with support on the {\it
Kac's spheres}, revisiting  \cite{CCLLV} and giving 	a
 possible answer to \cite[Open problem 11]{CCLLV}. 
 Last, we consider the 
 context of probability measures mixtures introduced by De Finetti, Hewitt and
Savage.    We define the (level 3)  Fisher information for mixtures and
prove  that it is l.s.c. and affine, as that was done in \cite{RR} for the level 3 Boltzmann's entropy.

\end{abstract}

\maketitle

\begin{center} {\bf  \today}
\end{center}

\smallskip


94A17 Measures of information, entropy, 

60F05 Central limit and other weak theorems,

82C40 Kinetic theory of gases. 

43A15 $L^p$-spaces and other function spaces on groups, semigroups

52A40 Inequalities and extremum problems
	
\vspace{0.3cm}

\bigskip

\tableofcontents


\section{Introduction and main results}
\label{sec:intro}
\setcounter{equation}{0}
\setcounter{theo}{0}


The Kac's notion of chaos rigorously formalizes the intuitive idea for a family of stochastic valued vectors with $N$ coordinates  to have
asymptotically independent coordinates as $N$ goes to infinity. We refer to \cite{S6}  for an introduction to that topics from a probabilistic point of view, 
as well as to \cite{Mxedp} for a recent and short survey.

\begin{defin}\cite[section 3]{Kac1956} \label{def:chaosKac1956} Consider $E \subset \R^d$, $f \in \PPP(E)$ a probability measure on $E$ and $G^N \in \Ps(E^N)$ a sequence of probability measures on $E^N$, $N \ge 1$, which are invariant under coordinates permutations. We say that $(G^N)$ is $f$-Kac's chaotic (or has the {\it ``Boltzmann property"}) if  
\beqn\label{defKac1}
\forall \, j \ge 1, \quad G^N_j \wto f^{\otimes j} \quad \hbox{weakly in } \PPP(E^j) \quad \hbox{as} \quad N\to\infty,
\eeqn
where $G^N_j$ stands for the $j$-th marginal of $G^N$ defined by 
$$
  G^N_j := \int_{E^{N-j}} G^N \, dx_{j+1} \, ... \, dx_N.
$$
\end{defin}

\smallskip
Interacting $N$-indistinguishable particle systems are naturally described by exchangeable random variables (which corresponds to
the fact that their associated probability laws are symmetric, i.e. invariant under coordinates permutations) but they are not described by 
random variables with independent coordinates (which corresponds to the fact that their associated probability laws are tensor products)
except for situations  with no interaction!  
Kac's chaos is therefore a well adapted concept to formulate and investigate the infinite number of particles limit $N\to\infty$
for these  systems  as it has been illustrated by many works since the seminal
article by Kac \cite{Kac1956}. Using the above definition of chaos, it  is shown
in  \cite{Kac1956,McKean1967,McKean-TCL,GM,MM-KacProgram} that if  $f(t)$ evolves
according to the nonlinear space homogeneous Boltzmann equation, $G^N(t)$
evolves according to the linear Master/Kolmogorov equation associated to the
stochastic Kac-Boltzmann jumps (collisions) process and $G^N(0)$ is
$f(0)$-chaotic, 
then for any later time $t>0$ the sequence  $G^N(t)$
is also $f(t)$-chaotic: in other words \emph{propagation of chaos} holds for
that model. 
As it is explained in the latest reference and using the uniqueness of
statistical solutions proved in \cite{ArkerydCI99}, 
some of these propagation of chaos results can be seen as an illustration of the
``BBGKY hierarchy method" whose most famous success is the Lanford's proof of 
the ``Boltzmann-Grad limit" \cite{Lanford}. 

\smallskip
In order to investigate quantitative version of Kac's chaos, the above weak convergence in \eqref{defKac1} can be formulated in terms of the 
Monge-Kantorovich-Wasserstein (MKW) transportation distance between $G^N_j$ and $f^{\otimes j}$. 
More precisely, given $d_E$ a bounded distance on $E$, we define the normalized
distance $d_{E^j}$ on $E^j$, $j \in \N^*$,  by setting
\beqn\label{def:distEj}
\forall \, X=(x_1, ...,x_j),Y=(y_1, ...,y_j) \in E^j \quad    d_{E^j} (X,Y)  := {1 \over j} \sum_{i=1}^j d_E(x_i,y_i),
\eeqn
and then we define  $W_1$ (without specifying the dependence on $j$)  the associated  MKW distance in $ \PPP (E^j)$ (see the definition \eqref{def:W1} below).  
With the notations of Definition~\ref{def:chaosKac1956}, $G^N$ is $f$-Kac's chaotic if, and only if, 
$$
\forall \, j \ge 1, \quad \Omega_j(G^N;f) := W_1(G^N_j,f^{\otimes j}) \, \to \, 0  \quad \hbox{as} \quad N\to\infty.
$$

\medskip
Let us introduce now another formulation of Kac's chaos which we firstly formulate in a probabilistic language. For any $X = (x_1, ...,x_N) \in E^N$, we define the associated  empirical measure 
\beqn\label{def:EmpirMeasure}
\mu^N_X (dy) := {1 \over N} \sum_{i=1}^N \delta_{x_i}(dy) \in \PPP(E). 
\eeqn
We say that an exchangeable $E^N$-valued  random vector  $\XX^N$  is $f$-chaotic  if the associated  $\PPP(E)$-valued random variable  $\mu^N_{\XX^N} $ converges  to the deterministic random variable $f$  in law in $\PPP(E)$:
\beqn\label{defKacProba} 
\mu^N_{\XX^N}  \, \Rightarrow \, f \quad \hbox{in law as} \quad N\to\infty. 
\eeqn

In the framework of Definition~\ref{def:chaosKac1956},  the convergence \eqref{defKacProba} can be equivalently formulated in the following way. 
Introducing $G^N := \LLL(\XX^N)$ the law of $\XX^N$, the exchangeability hypothesis means that $G^N \in \PPP_{\! sym}(E^N)$. Next the law $\hat G^N:= \LLL(\mu^N_{\XX^N} )$ of $\mu^N_{\XX^N} $ is nothing but the (unique) measure $ \hat G^N \in \PPP(\PPP(E))$ such that 
$$
\langle \hat G^N, \Phi \rangle = \int_{E^N} \Phi(\mu^N_X) \, G^N(dX) \qquad
\forall \, \Phi \in C_b(\PPP(E)),
$$
or equivalently the push-forward of $G^N$ by the ``empirical
distribution"  application.

Then the convergence \eqref{defKacProba} just means that 
\beqn\label{defKac2}
\hat G^N \to \delta_f \quad \hbox{weakly in } \PPP(\PPP(E))  \hbox{ as } N \to \infty,
\eeqn
where this definition does not refer anymore to the random variables $\XX^N$ or
 $\mu^N_{\XX^N} $.
It is  well known (see for instance  \cite[section 4]{Grunbaum}, \cite{Kac1956,T2,S4} and \cite[Proposition 2.2]{S6}) 
that for a sequence $(G^N)$ of $\PPP_{\! sym}(E^N)$  and a probability measure $f \in \PPP(E)$ the three following assertions are equivalent:

\smallskip
(i) convergence \eqref{defKac1} holds for any $j \ge 1$; 

\smallskip
(ii) convergence \eqref{defKac1} holds for some $j \ge 2$;

\smallskip
(iii) convergence \eqref{defKac2} holds;

\smallskip\noindent
so that in particular \eqref{defKac1} and \eqref{defKac2} are indeed equivalent formulations of Kac's chaos.  
The chaos formulation (ii) has been used since \cite{Kac1956}, while the chaos formulation (iii) is widely used in the works by 
Sznitman \cite{S1}, see also  \cite{S4,Meleard1996,MR839794},   where the chaos property is established by proving 
that the ``empirical process" $\mu^N_{\XX^N}$ converges to a limit process with values
in $\PPP(E)$   which is a solution to a nonlinear martingale problem associated to the mean-field limit equation. 
Formulation  \eqref{defKac2}  is also well adapted for proving quantitative propagation of chaos for deterministic dynamics associated to 
the Vlasov equation with regular interaction force  \cite{Dobrusin79} as well as singular interaction force
 \cite{MR2278413,HaurayJabinVlasov2}. Let us briefly explain this point now, see also \cite[section 1.1]{Mxedp}. 
 On the one hand, introducing the MKW transport distance $\WW_1 := \WW_{W_1}$ on $\PPP(\PPP(E))$ 
 based on the MKW distance $W_1$ on $\PPP(E)$, (see definition \eqref{def:WW1} below), the weak convergence  \eqref{defKac2} is nothing but the fact that 
$$
 \Omega_\infty (G^N; f) := \WW_1(\hat G^N,\delta_f) \, \to \, 0  \quad \hbox{as} \quad N\to\infty.
$$
On the other hand, for the Vlasov equation with smooth and bounded force term, it is proved in  \cite{Dobrusin79} that 
\beqn\label{eq:VlasovIneg}
\forall \, T > 0, \,\, \forall \, t \in [0,T] \qquad W_1(\mu^N_{\XX^N_t},f_t) \le C_T \, W_1(\mu^N_{\XX^N_0},f_0),
\eeqn
where $f_t \in \PPP(E)$ is the solution to the Vlasov equation with initial datum $f_0$ and $\XX^N_t \in E^N$ is the solution to the associated system of
ODEs with initial datum $\XX^N_0$.  Inequality \eqref{eq:VlasovIneg}  is a
consequence of the fact that $t \mapsto \mu^N_{\XX^N_t}$ solves the  Vlasov
equation and that a local $W_1$ stability result holds for such an equation.
When $\XX_0$ is distributed according to an initial density $G^N_0 \in \PPP_{\!
sym}(E^N)$ we may 
show that $\XX_t$ is distributed according to $G^N_t \in \PPP_{\! sym}(E^N)$ obtained as the transported measure along the  flow associated to the above
mentioned system of ODEs or equivalently $G^N_t$ is the solution to the
associated Liouville equation with initial condition $G^N_0$. Taking the
expectation in both sides of \eqref{eq:VlasovIneg}, we get 
\bean
\int_{E^N} W_1(\mu^N_Y,f_t) \, G^N_t(dY) &=& \Ee[ W_1(\mu^N_{\XX^N_t},f_t) ]  \\
&\le& C_T \, \Ee[ W_1(\mu^N_{\XX^N_0},f_0)]
= C_T \int_{E^N} W_1(\mu^N_Y,f_0) \, G^N_0(dY),
\eean
for any $t \in [0,T]$. We conclude with the following quantitative chaos propagation estimate 
$$
\forall \, t \in [0,T] \qquad \Omega_\infty (G^N_t; f_t) \le C_T \, \Omega_\infty (G^N_0;f_0) . 
$$
It is worth mentioning that  partially inspired from \cite{Grunbaum}, it is shown in \cite{MMWchaos,MM-KacProgram} a similar inequality as above for 
more general models including drift, diffusion and collisional interactions where however the estimate may mix several chaos quantification
quantities  as  $\Omega_\infty$ and $\Omega_2$ for instance.

\medskip
There exists at least one more way to guaranty chaoticity which is very popular because that chaos formulation naturally appears in the 
probabilistic coupling technique, see \cite{S6}, as well as 
\cite{MR1847094,BGV,MR2731396}  and the references therein.

Thanks to  the coupling techniques  we typically may show that an exchangeable
$E^N$-valued  random vector  $\XX^N$ satisfies 
 $$
\Ee \Bigl( {1 \over N} \sum_{i=1}^N | \XX^N_i - \YY^N_i | \Bigr) \to 0 \quad\hbox{as}\quad N \, \to \, \infty,
$$
for some   $E^N$-valued  random vector  $\YY^N$ with independent coordinates. Denoting by  $G^N\in\Ps(E^N)$ the law of $\XX^N$, $f$ the law of one 
coordinate $\YY^N_i$, and $W_1$ the MKW transport distance on $\PPP(E^N)$ based
on the normalized distance $d_{E^N}$ in $E^N$ defined by \eqref{def:distEj}, the
above convergence readily implies 
\beqn\label{defKac3}
\Omega_N(G^N;f) := W_1(G^N,f^{\otimes N}) \, \to \, 0 
\quad\hbox{as}\quad N \, \to \, \infty,
\eeqn
which in turn guaranties that $(G^N)$ is $f$-chaotic.  It is generally agreed
that the convergence \eqref{defKac3} is a strong version of chaos, maybe because it involves the all $N$ variables, while the Kac's original definition only involves
a finite fixed number of variables.

\medskip
\noindent
{\bf Summary of Section 2. } \Black
The first natural question we consider is about the
equivalence between these definitions of chaos, 
and more precisely the possibility to liken them in a quantitative way. The 
following result gives a positive answer, we also refer to
Theorem~\ref{Th:EquivChaosAccurate}  in section~\ref{sec:Kacchaos}
for a more accurate statement. 

\begin{theo}[Equivalence of measure for Kac's chaos] \label{Th:EquivChaos} 
For any moment order $k > 0$ and any positive exponent $\gamma <
(d+1+d/k)^{-1}$, 
there exists   a
constant  $C= C(d,k,\gamma) \in (0,\infty)$ such that for any  $f \in \PPP(E)$,  any
$G^N \in \PPP_{sym} (E^N)$, $N \ge 1$,  and any $j,\ell \in \{1, ... , N \} \cup \{ \infty \}$, 
$\ell \not= 1$, there holds
\[
\Omega_j (G^N;f) \le C \, \MMM_k^{1 / k} \Bigl( \Omega_\ell(G^N;f)+ {1 \over N} \Bigr) ^{\gamma} , 
\]     
where $\MMM_k = M_k(f) + M_k(G^N_1)$ is the sum of the moments of order 
$k$ of $f$ and $G^N_1$.
\end{theo}

It is worth emphasizing that the above inequality is definitively false 
in general for $\ell =1$. 
The first outcome of our theorem is that it shows that,  regardless of  the rate, 
the propagation of chaos results obtained 
by the coupling method is of the same nature as the propagation of chaos result 
obtained by  the ``BBGKY hierarchy method" and  the ``empirical
measures method". 

\smallskip
The proof of Theorem~\ref{Th:EquivChaosAccurate} (from which Theorem~\ref{Th:EquivChaos} follows) will be presented in section~\ref{sec:Kacchaos}. Let us briefly explain the strategy. 
First, the fact that we may control
$\Omega_j$ by $\Omega_\ell$ for $1 \le j \le \ell \le N$ is classical and quite easy. Next, we will establish an estimate of 
$\Omega_\infty$ by $\Omega_2$ following an idea introduced in \cite{MM-KacProgram}: we begin to prove a similar estimate where
we replace $\Omega_\infty$ by the MKW distance in $\PPP(\PPP(E))$ associated to
the $H^{-s}(\R^d)$  norm, $s > (d+1)/2$, on $\PPP(E)$
in order to take advantage  of the good algebraic structure of that Hilbert norm and then we come back to $\Omega_\infty$ thanks to the 
{\it ``uniform topological equivalence" }   of metrics in $\PPP(E)$ and the H\"older 
inequality.  Finally, and that is the other key new result, we compare $\Omega_\infty$ and $\Omega_N$: that is  direct consequence
of the following identity 
$$
\forall \, F^N, G^N \in \Ps(E^N) \qquad W_1(G^N,F^N)  = \WW_1(\hat G^N,\hat F^N) 
$$
applied to $F^N := f^{\otimes N}$ and a functional version of the law of large numbers. 

\bigskip\noindent
{\bf Summary of section 3.}
A somewhat stronger notion of chaos can be formulated in terms of entropy functionals. 
Such a notion has been explicitly introduced by Carlen, Carvahlo, Le Roux, Loss,  Villani in \cite{CCLLV} 
(in the context of probability measures  with support on the ``Kac's spheres")
but it is reminiscent of the works \cite{MR1872737,BenArousZeitouni99}. We also refer to 
 \cite{RR,MesserSpohn82,CLMP1,CLMP2} where the $N$ particles entropy functional below 
is widely used in order to identify the possible limits for a system of $N$ particles as $N\to\infty$. 
Consider $E \subset \R^d$ an open set or the adherence of a open space, in order that the gradient
of a function may be  well defined. For a (smooth  and/or decaying enough) probability measure $G^N \in \Ps(E^N)$ we define 
(see section~\ref{sec:Entropy} for the suitable definitions)
the  Boltzmann's entropy and the Fisher information by 
$$
H(G^N) := {1 \over N} \int_{E^N} G^N \log G^N \, dX, \qquad
I (G^N) := {1 \over N} \int_{E^N} {|\nabla G^N|^2 \over  G^N} \, dX.
$$
It is worth emphasizing that contrarily to the  most usual convention, adopted for instance in \cite[Definition
8]{CCLLV}, 
we have put the normalized factor $1/N$ in the definitions of the entropy and the Fisher information. Moreover we use the 
same notation for these functionals whatever is the dimension. As a consequence, we have $H(f^{\otimes N} ) = H(f)$ and
$I(f^{\otimes N} ) = I(f)$  for any probability measures $f \in \PPP(E)$. 
 
\begin{defin} \label{def:Entrop&FisherChaos}
Consider   $(G^N)$ a sequence of $\PPP_{\! sym}(E^N)$ such that for $k > 0$ the $k$-th moment $M_k(G^N_1)$ is  uniformly bounded in $N$, and $f \in \PPP(E)$. 
We say that  

(a) $(G^N)$ is $f$-entropy chaotic  (or $f$-chaotic in the sense of the Boltzmann's entropy)   if  
$$
G^N_1 \wto f \,\,\,\hbox{weakly in}\,\,\, \PPP(E) \quad\hbox{and}
\quad H(G^N) \to H(f), \,\,\, H(f) < \infty; 
$$

(b) $(G^N)$ is $f$-Fisher information chaotic (or $f$-chaotic in the sense of the Fisher information) if 
$$
G^N_1 \wto f \,\,\,\hbox{weakly in}\,\,\, \PPP(E) \quad\hbox{and}
\quad I(G^N) \to I(f), \,\,\, I(f) < \infty. 
$$
\end{defin}

Our second main result is the following qualitative comparison of the three above notions of chaos convergence. 

\begin{theo}\label{theo:LesDiffChaos} Assume $E = \R^d$, $d \ge 1$, or $E$ is a bi-Lipschitz volume preserving   deformation of a convex set of
$\R^d$, $d \ge 1$. Consider $(G^N)$ a sequence of $\PPP_{\! sym}(E^N)$ such that  the $k$-th moment $M_k(G^N_1)$ is  bounded, $k > 2$, and $f \in \PPP(E)$.  

\noindent
In the  list of assertions below, each one implies the assertion which follows:

(i) $(G^N)$ is $f$-Fisher information chaotic; 

(ii) $(G^N)$ is  $f$-Kac's chaotic and $I(G^N)$ is bounded;  

(iii) $(G^N)$ is $f$-entropy chaotic;

(iv) $(G^N)$ is  $f$-Kac's chaotic. 

\noindent
More precisely, the following quantitative estimate of the implication 
$(ii) \Rightarrow (iii)$ holds:  
\beqn
\label{eq:HGN-Hf<WN}
|H(G^N) - H(f)|\le  C_E \,  K \,\Omega_N(G^N;f)^\gamma , 
\eeqn
with $\gamma := 1/2-1/k$,   $K :=  \sup_N  I(G^N)^{1/2} \,   \sup_N M_k(G^N_1)^{1/k}$  and $C_E$ is a constant depending on the set $E$ 
(one can choose $C_E = 8$ when $E = \R^d$). 
\end{theo}

The implication $(ii) \Rightarrow (iii)$  is the most interesting part and hardest step in the proof of Theorem~\ref{theo:LesDiffChaos}. 
It  is based on estimate  \eqref{eq:HGN-Hf<WN} which is a mere consequence of the 
HWI inequality of Otto and Villani proved in \cite{OttoVillani} when $E = \R^d$

together with our  equivalence of chaos convergences previously established. It is also the most restrictive one in term of moment bound: the implication  $(ii) \Rightarrow (iii)$ 
requires a $k$-th moment bound of order $k >2$ while the other implications only require $k$-th moment bound of order $k >0$ or no moment bound condition  (we refer to the proof of Theorem~\ref{theo:LesDiffChaos}  in section~3 for details). 
  The proofs of the implications $(i)  \Rightarrow (ii)$ and $(iii)  \Rightarrow (iv)$ use the fact that the subadditivity inequalities of the Fisher information and of the entropy are saturated if and only if the probability measure is a tensor product. For  functionals  involving the entropy, similar ideas are classical and they have been used in \cite{MesserSpohn82,GPV,Yau,CLMP1,OVY} 
for instance.  

  We
believe that this result gives a better understanding of the different notions
of chaos.  Other but related notions of entropy chaos are introduced and discussed in 
\cite{CCLLV,MMgranular}. The entropy chaos definition in \cite{CCLLV}, which consists in asking
for point $(iii)$ and $(iv)$ above, is in fact equivalent to ours thanks to Theorem~\ref{theo:LesDiffChaos}.
 
 It is worth emphasizing that  Theorem~\ref{theo:LesDiffChaos} may be very useful in order to obtain entropic propagation of chaos (possibly with rate estimate) in contexts  where some bound on the Fisher information is available and propagation of Kac's chaos is already proved.  Unfortunately, a bound on the Fisher information is not   easy to propagate for $N$ particle systems. However, for the so-called {\it ``Maxwell molecules cross-section"},  following the proof of the fact that the Fisher information decreases along time for solutions to the homogeneous nonlinear Boltzmann equation  \cite{McK1,T92,VillaniFishBolt} and for solutions to the homogeneous nonlinear Landau equation \cite{VillaniFishLand}, it has been established that the $N$ particle Fisher information also decreases along time for  the law of solutions to  the stochastic Kac-Boltzmann jumps process in  \cite[Lemma 7.4]{MM-KacProgram} 
and for the law of solutions to the stochastic Kac-Landau diffusion process in  \cite{Carrapatoso2}.
 In these particular cases,   Theorem~\ref{theo:LesDiffChaos} provides  a quantitative version of the entropic propagation of chaos  proved  in \cite{MM-KacProgram}, and we refer to  \cite{Carrapatoso1,Carrapatoso2} for details.

\medskip\noindent
{\bf Summary of Section~\ref{sec:Sphere}.}

Here we consider the framework of probability measures with
support on the ``Kac's spheres" $\KK\SS_N$
defined by
$$
\KK\SS_N :=  \{ V = (v_1, ..., v_N) \in \R^N,  \,\,v_1^2 + ... + v_N^2 = N \},
$$
as firstly introduced by Kac in \cite{Kac1956}. Our aim is mainly to revisit the recent work \cite{CCLLV}
and to develop ``quantitative" versions of the chaos analysis. 

\smallskip
We start proving  a quantified  {\it ``Poincar\'e Lemma" } establishing that  the sequence of uniform probability measures $\sigma^N$ 
on $\KK\SS_N$ is $\gamma$-Kac's chaotic, with $\gamma$ the standard gaussian on $\R$, i.e. 
 $\gamma(v) = (2\pi)^{-1/2} \, \exp(-|v|^2/2)$, 
in the sense that we prove a rate of converge 
to $0$ for the  quantification of chaos $\Omega_N(\sigma^N;\gamma)$. We also prove  that for a large class of probability densities $f \in \PPP(E)$ the corresponding
sequence  $(F^N)$ of  {\it  ``conditioned to the Kac's spheres product measures" } (see section~\ref{subsec:condtensprod} for the precise definition)  is $f$-Kac's chaotic
in the sense  that we prove a rate of converge  to $0$ for the  quantification of chaos $\Omega_2(F^N;f)$.
That last result generalizes  the  {\it ``Poincar\'e Lemma" } since $f= \gamma$ implies $F^N = \sigma^N$. The main argument  in the last result is a (maybe new) $L^\infty$ optimal rate version of the Berry-Esseen theorem, also called local central limit theorem, which is nothing but an accurate (but less general) version of \cite[Theorem 27]{CCLLV}. Together with Theorem~\ref{Th:EquivChaos}, or the more accurate version of  it  stated in section~\ref{sec:Kacchaos}, we obtain the following
estimates.

\begin{theo}\label{theo:EmpiricalMeasure} 
The sequence $(\sigma^N)$ of uniform probability measures  on the ``Kac's spheres" is $\gamma$-Kac's chaotic, and more precisely 
\beqn\label{eq:ChaosEstimSigmaN}
\forall \, N \ge 1 \qquad 
\Omega_2(\sigma^N;\gamma) \le {C_1 \over N }, \quad \Omega_N(\sigma^N;\gamma) \le {C_2 \over N^{1 \over 2}} , \quad \Omega_\infty(\sigma^N;\gamma) \le
C_3   {(\ln N)^{1 \over 2} \over N^{1 \over 2}}, 
\eeqn
for some numerical constants $C_i$, $i=1,2,3$. 

\smallskip
More generally, consider $f \in \PPP(\R)$ with bounded moment $M_k(f)$ of order $k \ge 6$ and  bounded Lebesgue norm
$\| f \|_{L^p}$ of exponent $p > 1$. Then, the sequence $(F^N)$ of  associated ``conditioned (to the Kac's spheres) product measures"  is $f$-Kac's  chaotic, and more precisely
\beqn\label{eq:ChaosEstimFN}
\forall \, N \ge 1 \qquad 
\Omega_2(F^N;f) \le {C_4 \over N^{1 \over 2}}, \quad \Omega_N(F^N;f) \le {C_5 \over N^{\gamma \over 2}} , \quad \Omega_\infty(F^N;f) \le {C_6 \over N^{\gamma \over 2}},
\eeqn
for any $\gamma \in (0,(2+2/k)^{-1})$ and for some constants $C_i =
C_i(f,\gamma,k)$, $i=4,5,6$. 
 \end{theo}

Let us briefly discuss that last result. The question of establishing the convergence
 for the empirical law of large numbers associated to i.i.d. samples is an important question in theoretical statistics
 known as Glivenko-Cantelli theorem,  and the historical references seems to be  \cite{Glivenko,Cantelli,MR0094839}.
 Next the question of establishing rates of convergence in MKW distance in the above convergence has been addressed 
 for  instance in  \cite{MR0236977,MR779885,DubricYukich,BookRachev,MM-KacProgram,BoissardLeGouic}, 
 while the optimality of that rates have been considered for instance in \cite{MR779885,MR1189420,DubricYukich,BartheBordenave}. 
 We refer to \cite{BartheBordenave,BoissardLeGouic} and the references therein for a recent discussion on that topics. 
 With our notations, the question consists in establishing the
 estimate 
\beqn\label{eq:EmpiricLLN}
 \Ee (W_1(\mu^N_{\XX^N},f)) = \Omega_\infty(f^{\otimes N}; f) \le {C \over N^\zeta},
\eeqn
 for some constants $C = C(f)$ and $\zeta = \zeta(f)$. In the above left hand side term,   $\XX^N$ is a $E^N$-valued  random vector
 with independent coordinates with identical law $f$ or equivalently $\XX^N=X$ is the identity vector in $E^N$ and $\Ee$ is the expectation
 associated to the tensor product probability measure $f^{\otimes N}$.  When $E=\R^d$, estimate \eqref{eq:EmpiricLLN}   has been proved to hold with $\zeta = 1/d$,  if  
 $d \ge 3$ and supp$\, f$ is compact in \cite{DubricYukich}, 
 with $\zeta < \zeta_c := (d'+d'/k)^{-1}$, $d' = \max(d,2)$, if $d\ge1$ and  $M_k(f) < \infty$ in \cite{MM-KacProgram}
 and with $\zeta = \zeta_c $ if furthermore $d \ge 3$ in \cite{BoissardLeGouic}. 
 
 To our knowledge, \eqref{eq:ChaosEstimSigmaN} and
 \eqref{eq:ChaosEstimFN} are the first rates of convergence in MKW distance for the empirical law of large numbers associated to triangular
 array $\XX^N$ which coordinates are not i.i.d. random variables but only Kac's chaotic exchangeable random variables. The question of the 
 optimality of the rates in \eqref{eq:ChaosEstimSigmaN} and  \eqref{eq:ChaosEstimFN} is an open (and we believe interesting) problem. 
  
 \smallskip
 
Now, following \cite{CCLLV},  we introduce the notion of {\it entropy chaos} and {\it Fisher information chaos }
  in the context of the ``Kac's spheres" as follows. 
 For any $j \in \N$, and $f,g \in \PPP(E^j)$,
we define the usual relative entropy and usual relative Fisher information 
$$
H(f|g) := {1 \over j} \int_{E^j}  u \, \log u \, g(dv), \quad
I (f|g)  := {1 \over j} \int_{E^j} {|\nabla u|^2 \over  u}  \, g(dv), \quad u :=
\frac{df}{dg}, 
$$
where $ u = \frac{df}{dg}$ stands for the Radon-Nikodym derivative
 of $f$ with respect to $g$.

\noindent
For $f \in \PPP(E)$ and  $G^N \in \Ps(\KK\SS_N)$ such that $G^N_1 \wto f$ weakly in $\PPP(E)$, 
we say that  $(G^N)$ is 

\smallskip
(a$'$)  $f$-entropy chaotic  if $H(G^N|\sigma^N) \to H(f|\gamma)$, $H(f|\gamma) < \infty$;

\smallskip
(b$'$) $f$-Fisher information chaotic if $I(G^N|\sigma^N) \to I(f|\gamma)$, $I(f|\gamma) < \infty$. 

\smallskip 
In a next step, we prove that for a large class of probability measures  $f \in \PPP(\R)$ the sequence $(F^N)$ of  associated ``conditioned (to the Kac's spheres) product measures" is  $f$-entropy chaotic as well as
 $f$-Fisher information chaotic,  and we exhibit again rates for these convergences. 
 The proof is mainly  a careful rewriting and simplification of the proofs of the similar results 
(given without rate) in Theorems~9, 10, 19, 20 \& 21 in \cite{CCLLV}.

\smallskip
We next generalize Theorem~\ref{theo:LesDiffChaos} to the Kac's spheres context. 
Additionally to the yet mentioned arguments, we use a general version of the HWI inequality proved by Lott and Villani in \cite{MR2480619}, 
see also \cite[Theorem 30.21]{VillaniOTO&N}, and some entropy and Fisher inequalities on the Kac's spheres established by Carlen et al. \cite{CarlenLL2004}
and improved by Barthe et al. \cite{Barthe&co2004}. 

\smallskip
All these results are motivated by the question of giving quantified strong version of propagation of chaos 
for Boltzmann-Kac jump model studied in \cite{MM-KacProgram} by Mouhot and the
second author, where only quantitative uniform in time Kac's chaos is established. As a matter of fact, K. Carrapatoso in \cite{Carrapatoso1} extends
the present analysis to the probability measures with support to the {\it Boltzmann's spheres} and proves a quantitative propagation result of entropy chaos.  

\smallskip
Another outcome of our results is that we are able to give the following possible answer to \cite[Open problem 11]{CCLLV}:   

\begin{theo}\label{theo:OpenPb} Consider $(G^N)$ a sequence of $\PPP_{\! sym}(\R^N)$ with support on the Kac's spheres $\KK\SS_N$ such that 
\beqn\label{eq:MkGN1&IGN}
 M_k(G^N_1) \le C, \quad 
I(G^N|\sigma^N) \le C,
\eeqn
for  some $k  \ge 2$ and 
$C > 0$. Also consider  $f \in \PPP( \R)$,  satisfying $\int v^2 f(v)\,dv=1$ and 
\beqn\label{eq:f>exp}
f \ge \exp ( - \alpha \, |v|^{k'} + \beta) \quad\hbox{on}\quad \R,
\eeqn
with $0 < k' < k$, $\alpha > 0$, $\beta \in \R$. If $(G^N)$ is $f$-Kac's chaotic,   then for any fixed $j \ge 1$, there holds 
$$
H(G^N_j | f^{\otimes j}) \, \to \, 0 \quad \hbox{as}\quad N \to\infty,
$$
where $H(\cdot|\cdot)$ stands for the usual relative entropy functional defined in the flat space $E^j$. 
Remark that the boundedness of the $k$-th moment of $G^N$ is useless when $k\le2$ (because the support condition implies $M_2(G^N_1) = 1$) while
the condition on the second moment of $f$ is useless if $k >2$ (because it is inherited from the properties  of $(G^N_1)$).  
\end{theo}

Contrarily to the conditioned tensor product assumption made in \cite[Theorem 9]{CCLLV} which can be assumed at initial time for the 
stochastic Kac-Boltzmann process but which is not propagated along time, our assumptions \eqref{eq:MkGN1&IGN} and \eqref{eq:f>exp} in Theorem~\ref{theo:OpenPb}, which may seem to be stronger in some sense,  are in fact more natural since they are propagated along time. 
We refer to \cite{MM-KacProgram,Carrapatoso1} where such problems are studied.

\medskip\noindent
{\bf Summary of Section~\ref{sec:mixtures}.}
Here we investigate how our techniques can be
used in the context of probability measures mixtures as introduced
by De Finetti, Hewitt and Savage \cite{deFinetti1937,HewittSavage} and general sequences of probability 
densities $G^N$  of $N$ undistinguishable particles  as $N\to\infty$, without assuming chaos, as  it is the case in  \cite{MesserSpohn82,CLMP1,CLMP2} for instance. The results developed in that section are also used in a fundamental way in the recent work \cite{FHM}. 

In a first step, we give a new proof of  De Finetti, Hewitt and Savage theorem which is based on the use of the law of the empirical measure associated to
the $j$ first coordinates like in Diaconis and Freedman's proof \cite{DiaconisFreedman1987} or  Lions' proof \cite{PLLcoursCF}, but where the compactness
arguments are replaced by an argument of completeness. As a back product, we give a quantified equivalence of several notions of convergences 
of sequences of $\Ps(E^N)$ to its possible mixture limit. 

 In a second step, we revisit the level 3 entropy and level 3 Fisher information theory for a probability measures mixture as developed since the work by Robinson and Ruelle \cite{RR} at least. We give a comprehensive and elementary proof of the fundamental result 
\beqn\label{eq:introRR}
 \KK(\pi) := \int_{\PPP(E)} K(\rho) \, \pi(d\rho) = \lim_{j\to\infty} {1 \over j} \int_{E^j} K(\pi_j) 
\eeqn
 for any probability measures mixture $\pi \in \PPP_k(\PPP(E))$, $k > 0$ (see paragarph~\ref{subsec:DeFHS} where the space $ \PPP_k(\PPP(E))$ is defined), where $\pi_j$ stands for the  De Finetti, Hewitt and Savage projection of $\pi$ on the $j$ first coordinates
 and $K$  stands for the Boltzmann's entropy or the Fisher information functional. It is worth noticing that while the representation formula \eqref{eq:introRR} is well known when $K$ stands for the Boltzmann's entropy, we believe that it is new when $K$ stands for the Fisher information. The representation formula for the Fisher information is interesting for its own sake and it has also found an application as a key argument in the proof of propagation of chaos for system of vortices established in \cite{FHM}. 

 In our last result we establish a rate of convergence for the above limit \eqref{eq:introRR} when $K$ is the entropy functional mainly under a boundedness of the Fisher information hypothesis and we generalize such a quantitative result establishing links between several weak notions  of convergence as well as strong (entropy) notion of  convergence for sequences of probability densities $G^N \in \Ps(E^N)$  as $N\to\infty$, without assuming chaos.

\medskip
\noindent
{\bf Acknowledgement. }
 The authors would like to thank F. Bolley and C. Mouhot for many stimulating discussions about mean field limit and chaos,
as well as  N. Fournier for his suggestions that make possible to improve the statement and simplify the proof
of the result on the level-3 Fisher information in section~\ref{sec:mixtures} and A.~Einav for having pointing out a mistake in a previous version of the work.  
The second author also would like to acknowledge I. Gentil and C. Villani for discussions 
about the HWI inequality and \break P.-L.~Lions for discussions about entropy and mollifying tricks in infinite dimension. 
Finally, we would like to thank the anonymous referees  for their comments and suggestions making possible significant improvements in 
the presentation of  the article.  
The second author acknowledges support from the project ANR-MADCOF.


\section{Kac's chaos}
\label{sec:Kacchaos}
\setcounter{equation}{0}
\setcounter{theo}{0}


In this section we show the equivalence between several ways to measure  Kac's chaos as stated in Theorem~\ref{Th:EquivChaos}.
We start  presenting the framework we will deal with in the sequel, and thus making precise the definitions and notations 
used in the introductory section.

\subsection{Definitions and notations}\label{subsect:def}

In all the sequel, we denote by $E$ a closed subset of $\R^d$, $d \ge1$, endowed
with the usual topology, so that it is a locally compact Polish space. We denote
by $\PPP(E)$ the space of probability measures on the Borel $\sigma$-algebra
$\BBB_E$ of $E$. 

\medskip \noindent
{\bf Monge-Kantorovich-Wasserstein (MKW) distances.}

As they will be a cornerstone in that article, used in different setting, we
briefly recall their definition and main properties, and refer to
\cite{VillaniTOT} for a very nice presentation. 

On a general Polish space $Z$, for any  distance $D :Z \times Z \to \R^+$
and  $p \in [1,\infty)$,we define $W_{D,p}$ on $\PPP(Z)\times \PPP(Z) $ by setting for any 
$\rho_1, \rho_2 \in \PPP(Z)$
$$
[ W_{D,p}(\rho_1,\rho_2) ]^p := \inf_{\pi \in \Pi(\rho_1,\rho_2)}   
\int_{Z \times Z} D(x,y)^p \, \pi(dx, dy) 
$$
where $\Pi(\rho_1,\rho_2)$ is the set of proability measures 
$\pi \in \PPP(Z \times Z)$ with first marginal $\rho_1$ and second
marginal  $\rho_2$, that is $\pi (A \times Z) = \rho_1(A)$ and $\pi (Z \times A) = \rho_2(A)$  for any Borel set $A \subset Z$. 
It defines a distance on $\PPP(Z)$. 

\medskip \noindent
{ \bf The phase spaces $E^N$ (its marginal's space $E^j$) and $\PPP(E)$.}

When we study system of $N$ particles, the natural phase space is $E^N$. 
The space of marginals $E^j$ for $1 \leq j \leq N$ are also important.
We present here the different distances we shall use on these spaces.

\noindent $\bullet$ On $E$  we will use mainly two distances :
\begin{itemize}
\renewcommand{\labelitemi}{$-$}
\item the usual  Euclidean distance denoted by $|x-y|$; 
\item a bounded version of the square distance : $d_E(x,y) = |x-y| \wedge1$ for
any $x,y \in E$.
\end{itemize}

\noindent $\bullet$ On the space $E^j$ for $1 \leq j$, we will also use the
two distances
\begin{itemize}
\renewcommand{\labelitemi}{$-$}
\item the \emph{normalized} square distance $|X-Y|_2$ defined for any $X=(x_1,\ldots,x_j) \in E^j$ and $Y=(y_1,\ldots,y_j)\in E^j$ by
$$
|X-Y|_2^2 := \frac1j \sum_{i=1}^j
|x_i-y_j|^2;
$$
\item the \emph{normalized} bounded distance $d_j = d_{E^j}$ defined  by
\beqn\label{def:dEj}
d_{E^j} (X,Y) := {1 \over j} \sum_{i=1}^j
d_E(x_i,y_i) .
\eeqn
\end{itemize}
It is worth emphasizing that the normalizing factor $1/j$ is important in
the sequel in order to obtain formulas independant of the number $j$ of
variables.

\smallskip

\noindent $\bullet$ The introduction of the empirical measures allows to 
``identify"  our phase space $E^N$ to a subspace of $\PPP(E)$. To be more
precise,  we denote by $\PP_N (E)$ the set of empirical measures
$$
\PP_N(E) := \left\{ \mu^N_X  , \,\, X =
(x_1, ... , x_N) \in E^N \right\} \subset \PPP(E),
$$
where $ \mu^N_X$ stands for the empirical measure defined by \eqref{def:EmpirMeasure} 
and associated to the configuration $X
=(x_1,\ldots,x_n) \in E^N$.
We denote by $p_N : E^N \rightarrow \PP_N(E)$ the application that maps a
configuration to its empirical measure : $p_N(X) := \mu^N_X $. 

\smallskip
\noindent $\bullet$ On our phase space $\PPP(E)$, we will use three different
distances

\smallskip \noindent
- The usual MKW distance of order two $W_2$
defined as above with the choice $D(x,y) =|x-y|^2$
$$
W_2(\rho_1,\rho_2)^2 =  W_{|\cdot|_2,2} (\rho_1,\rho_2)^2 := \inf_{\pi \in \Pi(\rho_1,\rho_2)}   
\int_{E \times E} |x-y|^2 \, \pi(dx, dy) .
$$
\smallskip \noindent 
- The MKW distance $W_1$ associated to $d_E$ defined by
\beqn\label{def:W1} 
W_1(\rho_1,\rho_2) = W_{d_E,1} (\rho_1,\rho_2) := \inf_{\pi \in \Pi(\rho_1,\rho_2)}   
\int_{E \times E} d_E(x,y) \, \pi(dx, dy) .
\eeqn
From the Kantorovich-Rubinstein duality theorem (see for instance 
\cite[Theorem 1.14]{VillaniTOT})  we have the following alternative
characterization 
\beqn\label{def:W1KR}
\forall \, \rho_1,\rho_2 \in \PPP(E) \qquad W_1(\rho_1, \rho_2) =
\sup_{\|\varphi
\|_{Lip}\le 1}\int_{E}\varphi(x) \, (\rho_1(dx) - \rho_2(dx)), 
\eeqn
where $\| \varphi\|_{Lip} := \sup_{x \neq y} \frac{|\varphi(x) -
\varphi(y)|}{d_E(x,y)}$ is the Lipschitz semi-norm  relatively to the distance
$d_E$. This semi-norm is closely related to the usual Lipschitz semi-norm since
it satisfies 
\beqn \label{eq:dist_Lip}
\frac12 \left( \| \nabla \varphi \|_\infty + \| \varphi  -
\varphi(0) \|_\infty \right) \leq \| \varphi \|_{Lip} \leq 2 \left( \| \nabla
\varphi \|_\infty + \| \varphi\|_\infty\right) =: 2\, \|
\varphi\|_{W^{1,\infty}}.
\eeqn 
 It implies  that $W_1$ is equivalent to the
$(W^{1,\infty})'$-distance, denoted by $D_{W^{1,\infty}}$,
$$
D_{W^{1,\infty}} (\rho_1,  \rho_2) := \sup_{\|\varphi\|_{W^{1,\infty}}\le 1}\int_{E}\varphi(x) \,
(\rho_1(dx) - \rho_2(dx)), 
$$
and more precisely
\beqn \label{eq:EquivDistLip}
  \frac12\, D_{W^{1,\infty}}  \le W_1 \le 2\,
D_{W^{1,\infty}}.
\eeqn

\smallskip \noindent
- The distance induced by the $H^{-s}$ norm for $s> \frac d2$ : for any $\rho, \,
\eta \in \PPP(E)$
$$
\| \rho - \eta \|_{H^{-s}}^2 := \int_{\R^d}  |\hat \rho(\xi) - \hat \eta(\xi) |^2 
  \frac{d\xi}{\langle  \xi \rangle^{2s}}
$$
where $\hat \rho$ denotes the Fourier transform of $\rho$ (which may always be
seen as a measure on the whole $\R^d$), and $ \langle  \xi \rangle = \sqrt{1+
|\xi|^2}$. 

\smallskip \noindent $\bullet$
We will often restrict ourself to the spaces $\PPP_k(E)$ of probability measures with
finite moment of order $k > 0$ defined by 
$$
\PPP_k(E) := \{ \rho \in \PPP(E) \text{ s.t. } \; M_k(\rho) := \int_E
\langle v \rangle^k \,\rho(dv) < + \infty \}. 
$$

\smallskip
\paragraph{\bf The probability measures space $\PPP(E^N)$, its marginals spaces
$\PPP(E^j)$, and $\PPE$.}
The next step is to consider probability measures on the configuration spaces.

\smallskip \noindent
$\bullet$ The space $\PPP(E^N)$ will be endowed with two distances 
\begin{itemize}
\renewcommand{\labelitemi}{$-$}
\item $W_1$ the MKW distance on $\PPP(E^N)$ associated to $d_{E^N}$ and $p=1$, which 
has the same properties as the one constructed on  $\PPP(E)$ and satisfies in particular
the Kantorovich-Rubinstein formulation~\eqref{def:W1KR}.  
\item  $W_2$ the MKW distance associated to  the normalized square distance
$|\cdot|_2$ defined above.
\end{itemize}
 Remark that we will only work on the subspace $\Ps(E^N)$ of borelian
probability measures which are invariant under coordinates permutations. 

\smallskip \noindent $\bullet$ On the probability measures space $\PPE$, we can define  
different distances thanks to the Monge-Kantorovich-Wasserstein construction.
We will use three of them
\begin{itemize}
\renewcommand{\labelitemi}{$-$}
\item $\WW_1$, the MKW distance induced by the cost function $W_1$  on
$\PPP(E)$. In short
\beqn\label{def:WW1}
\WW_1(\alpha_1,\alpha_2) = W_{W_1,1} (\alpha_1,\alpha_2) := \inf_{\pi \in \Pi(\alpha_1,\alpha_2)}   
\int_{\PPP(E) \times \PPP(E)} W_1(\rho_1,\rho_2) \, \pi(d\rho_1, d\rho_2),
\eeqn
\item  $\WW_2$, the MKW distance induced by the cost function $W^2_2$ on
$\PPP(E)$. In short 
$$
\WW_2(\alpha_1,\alpha_2)^2= W_{W_2,2} (\alpha_1,\alpha_2)^2 := \inf_{\pi \in \Pi(\alpha_1,\alpha_2)}   
\int_{\PPP(E) \times \PPP(E)} W_2^2(\rho_1,\rho_2) \, \pi(d\rho_1, d\rho_2),
$$
\item $\WWs$, the MKW distance induced by the cost function $\| \cdot
\|_{H^{-s}}^2$ on $\PPP(E)$. In short 
$$
\WWs(\alpha_1,\alpha_2)^2 = \WW_{\| \cdot \|_{H^{-s}},2}(\alpha_1,\alpha_2)^2:= \inf_{\pi \in
\Pi(\alpha_1,\alpha_2)}  
\int_{\PPP(E) \times \PPP(E)} \| \rho_1 - \rho_2 \|_{H^{-s}}^2 \, \pi(d\rho_1,
d\rho_2).
$$
\end{itemize}

\smallskip \noindent $\bullet$ Remark that the application ''empirical
measure`` $p_N$ allows to define by push-forward a canonical map between $\PPP(E^N)$ 
and $\PPE$.  For $G^N \in \PPP(E^N)$ we denote its image under the application
$p_N$ by  $\hat G^N \in \PPE$ : $ \hat G^N := G^N_\# p_N$. In other words, $ \hat G^N$ is
the unique probability measure in $\PPE$ which satisfies the duality relation
\beqn\label{def:hatFN}
\forall \, \Phi \in C_b(\PPP(E)) \qquad 
\langle \hat G^N, \Phi \rangle = \int_{E^N} \Phi(\mu^N_X) \, G^N(dX).
\eeqn

\medskip \noindent
{\bf More properties of the space $\PPE$.}

$\bullet$ {\bf Marginals of probability measures on $\PPE$.} We can define a mapping
form $\PPE$ onto $\PPP(E^j)$ in the following way. For any $\alpha \in \PPE$ we
define the projection $\alpha_j \in \PEj$ thanks to the relation 
$$
\alpha_j := \int \rho^{\otimes j} \,d\alpha(\rho). 
$$

It may also be restated using \emph{polynomial fonctions} : 
for any $\varphi \in C_b(E^j)$ we define the 
monomial  (of order $j$) function $R_\varphi \in C_b(\PPP(E))$ by 
\beqn\label{def:polyproba}
\forall \, \rho \in \PPP(E) \qquad R_\varphi (\rho) := \int_{E^j} \varphi (X) \, \rho^{\otimes j} (dX).
\eeqn
 We remark that the monomial functions of all orders generate an
algebra of continous fonction (for the weak convergence of measures) that
are called polynomials. When $E$ is compact so that $\PPP(E)$ is also compact, they form a dense
subset of $C_b(\PPP(E))$ thanks to the Stone-Weierstrass theorem. 

In terms of polynomial fonctions, the marginal $\alpha_j$ may be defined by 
$$
\forall \, \varphi \in C_b(E^j) \qquad \langle \alpha_j , \varphi \rangle := \langle \alpha, R_\varphi \rangle.
$$

$\bullet$ Starting from $G^N \in \Ps(E^N)$, we can define its push-forward
$\hat G^N$ and then for any $1 \le j \le N$ the marginals of the push-forward
$\hat G^N_j := (\hat G^N)_j \in \Ps(E^j)$. They satisfy the duality relation
\beqn\label{def:FNj}
\forall \, \varphi \in C_b(E^j) \qquad \langle \hat G^N_j , \varphi \rangle := \int_{E^N} R_\varphi(\mu^N_X) \, G^N(dX).
\eeqn
 We emphasize that it is not equal to $G^N_j$ the $j$-th marginal of $G^N$,
but we will see later that the two probability measures $G^N_j$ and $\hat G^N_j$ are
close (a precise version is recalled in Lemma~\ref{lem:HSEquivQuant}).  

\medskip
\paragraph{ \bf Different quantities mesuring chaoticity.}
Now that everything has been defined, we introduce the quantities that we will
use to quantify the chaoticity of a sequence $G^N \in \Ps(E^N)$ of symmetric
probability measure with respect to a profil $f \in \PPP(E)$:
\begin{itemize}
\renewcommand{\labelitemi}{$-$}
\item The chaoticity can be mesured on $E^j$ for $j \geq 2$. For any $1 \leq j
\leq N$, we set
$$
\Omega_j (G^N; f) := W_1 (G^N_j,f^{\otimes j}),
$$
\item and also on $\PPP(E)$ by
$$
\Omega_\infty(G^N; f) := \WW_1 (\hat G^N, \delta_f) =  \int_{E^N}W_1(\mu^N_X,f) \,
G^N(dX),
$$
since there is only one transference plan $\alpha \otimes \delta_f$ in $\Pi
(\alpha,\delta_f)$.
\end{itemize}

\subsection{Equivalence of distances on $\PPP(E)$, $\PPP_{sym}(E^N)$ and
$\PPE$.} . 

To quantify the equivalence between the distances defined above 
on $\PPP(E)$, we will need some assumption on the moments. The metrics $W_1 $,
$W_2$ and $\| . \|_{H^{-s}}$ are uniformly topologicaly equivalent
in $\PPP_k(E)$ for any $k> 0$.
More precisely, we have 

\begin{lem}\label{lem:ComparDistances}
Choose $f,g \in \PPP(E)$. For any $k>0$, denote
$ \MMM_k := M_k(f) + M_k(g)$.

(i) For any $k >0$ and $s \ge 1$, 
there exists \mbox{$C :=
C(d) \, \Bigl[ 1 + \bigl( \frac{s-1}2 \bigr)^{\frac{s-1}2}
\Bigr]$}, such that 
there holds
\begin{equation}\label{estim:W1Hs}
W_1(f,g) \le  C \, \MMM_k^{\frac d{d+2ks}} \,  \| f- g
\|_{H^{-s}}^{\frac{2k}{d+2ks}}.
\end{equation}
(ii) For any $k>2$, there holds
\begin{equation}\label{estim:W2W1}
W_2(f,g) \le 2^{\frac32} \, \MMM_k^{1/k}\, W_1
(f,g)^{1/2-1/k} . 
\end{equation}
(iii) Without moment assumptions and for any $s> \frac{d+1}2$, there exists a constant $C = C(s,d)$
such that 
there holds
$$
W_1(f,g) \leq W_2(f,g), \qquad  \| f- g \|_{H^{-s}} \leq C \,
W_1(f,g)^{\frac12} . 
$$
 \end{lem}
 
We remark that we have kept the explicit dependance on $s$ of the constant
appearing in $(i)$ in order to be able to perform some optimization on $s$
later. The important point is that the constant may be choosen independant of
$s$ if $s$ varies in a compact set. 

 \medskip\noindent
{\sc Proof of Lemma~\ref{lem:ComparDistances}. } The proof is a mere
adaptation of classical results on comparison of distances in probability measures spaces
as it can be found in \cite{BookRachev,coursCT,MM-KacProgram} for instance. 
We nevertheless sketch  it for the  sake of completness.

\smallskip \noindent {\sl Proof of i). }
We consider a truncation sequence $\chi_R(x) = \chi(x/R)$, $R>0$, with $\chi \in
C^\infty_c(\R^d)$, $\|\nabla \chi \|_\infty \le 1$, $0 \le \chi \le 1$,
$\chi \equiv 1$ on
$B(0,1)$,  and the sequence of mollifiers  $\gamma_\eps(x) = \eps^{-d} \,
\gamma(x/\eps)$, $\eps>0$,  with $\gamma(x) =
(2\pi)^{-d/2} \, \exp(-|x|^2/2)$, so that $\hat \gamma_\eps (\xi) =  \exp (-
\eps^2 \, |\xi|^2/2)$. 
In view of the equivalence of distance~\eqref{eq:EquivDistLip}, we choose a
$\varphi \in W^{1,\infty}(\R^d)$ such that $\|
\varphi\|_{W^{1,\infty}}  \le 1$,
we define $\varphi_R := \varphi \, \chi_R$, $\varphi_{R,\eps} =
\varphi_R \ast \gamma_\eps$ and we write
\[
 \int \varphi \, (df-dg) = \int \varphi_{R,\eps} \, (df-dg) +
 \int \left(\varphi_R - \varphi_{R,\eps}\right) \, (df-dg) + 
 \int \left(\varphi - \varphi_R \right) \, (df-dg).
\]
For the last term, we have
\[
\quad \forall \, R > 0 \qquad
  \left|  \int  (\varphi_R- \varphi) \, (df-dg) \right| 
 \le \int_{B_R^c} \| \varphi \|_\infty \, {|x|^{k} \over R^k} \, (df+dg) \le 
  \frac{M_{k}}{R^k}.
\]
For the second term, we observe that 
$$
\| \varphi_R - \varphi_{R,\eps} \|_\infty \le \| \nabla
\varphi_R \|_\infty \int_{\R^d} \gamma_\eps (x) \, |x| \, dx \le C(d) \,
\eps,
$$
and we get
\[
\left| \int  \left(\varphi_R - \varphi_{R,\eps} \right) \, (df  -
  dg) \right| \le C(d) \, \eps.
\]
Finally, the first term can be estimated by 
\[
\left| \int \varphi_{R,\eps} \, (df-dg) \right| 
 \le 
\| \varphi_{R,\eps} \|_{H^s} \, \| f - g \|_{H^{-s}}, 
\]
with for any  $R \ge 1$ and $\eps \in (0,1]$
\bean
  \|  \varphi_{R,\eps} \|_{H^s} 
  &=&  \left(  \int \langle \xi \rangle^2 \, |\widehat{\varphi \, \chi_R}|^2 \,
\langle \xi \rangle^{2(s-1)} \, |\hat\gamma_\eps|^2 \, d\xi \right)^{1/2}
  \\ 
  &\le& \| \varphi \, \chi_R \|_{H^1} \, \| \langle \xi \rangle^{s-1} \, \hat 
\gamma_\eps (\xi) \|_{L^\infty}
 \le  C(d) \, R^{d/2} \| \langle \xi \rangle^{s-1} \, \hat 
\gamma_\eps (\xi) \|_{L^\infty} \,  
\eean
 The infinite norm is finite and a simple optimization leads
to
$$
\| \langle \xi \rangle^{s-1} \, \hat \gamma_\eps (\xi) \|_{L^\infty} \le 
\Bigl( \frac{s-1}2 \Bigr)^{\frac{s-1}2}\eps^{-(s-1)_+},
$$
with the natural convention $0^0=1$. All in
all, we have
$$ 
W_1(f,g) \le   C(d) \left[ 1 + \Bigl( \frac{s-1}2 \Bigr)^{\frac{s-1}2}
\right]  \left( \eps + \frac{M_{k}}{R^k} +
  R^{ \frac d2} \eps^{-(s -1)}  \, \| f - g \|_{H^{-s}} \right) . 
$$
This yields to \eqref{estim:W1Hs} by optimizing the paramater $\eps$ and $R$
with
$$
R = M_k^{\frac {2s}{d+2ks}} \,\| f - g \|_{H^{-s}}^{-\frac 2{d+2ks}}, \quad \text{and } \eps = M_k^{\frac d{d+2ks}} \| f -
g \|_{H^{-s}}^{\frac{2k}{d+2ks}}.
$$

\smallskip \noindent {\sl Proof of ii).}
We have for any $R \ge 1$ the inequality
$$
\forall \, x,y \in E, \quad |x-y|^2 \le R^2 \, d_E(x,y) + {2^k \over R^{k-2}} \,
(|x|^k + |y|^k)
$$
from which we deduce
\bean
W_2(f,g) ^2 
&\le& R^2 \inf_{\pi \in \Pi(f,g)}\int_{E \times E}
d_E(x,y) \, \pi(dx,dy)
\\
&& \hspace{1cm} +  {2^k \over  R^{k-2}} \sup_{\pi \in \Pi(f,g)} \int_{E
\times E} ( |x_i|^k + |y_i|^k)  \, \pi(dx,dy)
\\
&\le& R^2 \, W_1 (f,g)  +  {2^k \over  R^{k-2}} \, (M_k(f) +
M_k(g)), 
\eean
and then we get with $(R/2)^k = \MMM_k / W_1$
\beqn \label{ineq:tildeW2<W1}
 W_2(f,g) \le 2^{3/2} \MMM_k^{1/k}\, W_1 (f,g)^{1/2-1/k} .
\eeqn

\smallskip \noindent {\sl Proof of iii).} The first point is classical. The
second relies on the fact that \\ \mbox{$\| \delta_x - \delta_y \|_{H^{-s}}^2
\leq C d_E(x,y)$}.
\qed

\medskip
There is also a similar result on $E^N$, where the $H^{-s}$ norm is less
usefull.
\begin{lem} \label{lem:ComparDistEN}
Choose $F^N,G^N \in \PPP_{sym}(E^N)$. For any $k>0$, denotes 
$$ \MMM_k := M_k(F^N_1) + M_k(G^N_1).$$
For any $k>2$, it holds that
\begin{equation}\label{estim:W2W1EN}
W_2(F^N,G^N) \le 2^{\frac32} \, \MMM_k^{1/k}\, W_1
(F^N,G^N)^{1/2-1/k}.
\end{equation}
It also holds without moment assumptions that 
$W_1(F^N,G^N) \leq W_2(F^N,G^N)$.
\end{lem}

\medskip\noindent
{\sc Proof of Lemma~\ref{lem:ComparDistEN}. } 
The proof is a simple generalization of \eqref{estim:W2W1} to the case of $N$
variables. We skip it.  
\qed

\medskip
The inequalities of Lemma~\ref{lem:ComparDistances} also sum well on $\PPE$ in
order to get
\begin{lem} \label{lem:ComparDistPPE}
  Choose $\alpha, \beta \in  \PPE$, and for $k > 0$ define
$$
\quad \MMM_k := M_k(\alpha) + M_k(\beta) := \int M_k(\rho) \,[\alpha+
\beta](d\rho) = M_k(\alpha_1) + M_k(\beta_1).
$$
(i) For any $s \ge 1$ and with the same constant $C(d,s)$ as in point $(i)$
of Lemma~\ref{lem:ComparDistances} we have  for any $k >0$, 
\begin{equation}\label{estim:W1HsPPE}
\WW_1(\alpha,\beta) \le  C \, \MMM_k^{\frac d{d+2ks}} \, 
\WWs(\alpha,\beta)^{\frac {2k}{d+2ks}}.
\end{equation}
ii) For any $k>2$, it also holds 
\begin{equation}\label{estim:W2W1PPE}
(ii) \qquad 
\WW_2(\alpha,\beta) \le  2^{\frac32} \, \MMM_k ^{\frac 1k} \, 
\WW_1(\alpha,\beta)^{\frac12 - \frac1k},  
\end{equation}
iii) It holds without moment assumption that $\WW_1 \leq \WW_2$ and $\WWs \leq
C \WW_1^{\frac12}$ for $s>\frac{d+1}2$ with a  constant $C=C(s,d)$.
\end{lem}

\medskip\noindent
{\sc Proof of Lemma~\ref{lem:ComparDistPPE}.} 
All the above estimates are simple summations of the corresponding estimate of
Lemma~\ref{lem:ComparDistances}. We only prove $i)$.
\bean
\WW_1 (\alpha, \beta) 
&=& \inf_{\Pi \in \Pi(\alpha, \beta)} \int  W_1(\rho, \eta) \, \Pi(d\rho,d\eta)
\\
&\le& C \inf_{\Pi \in \Pi(\alpha, \beta)} \int [M_k(\rho) +
M_k(\eta)]^{\frac d{d+2ks}} \|\rho - \eta \|_{H^{-s}}^{\frac{2k}{d+2ks}} \,
\Pi(d\rho,d\eta) \\
&\le& C \left(\int  M_k(\rho)
\, [\alpha + \beta](d\rho) \right)^{\frac d{d+2ks}} \left(  \inf_{\Pi \in
\Pi(\alpha, \beta)} \int \|\rho - \eta \|_{H^{-s}}^{\frac1s} \, \Pi(d\rho,d\eta)
\right)^{\frac{2ks}{d+2ks}}\\
&\le& C [M_k(\alpha) + M_k(\beta)]^{\frac d{d+2ks}} \left( 
\inf_{\Pi \in \Pi(\alpha, \beta)} \int \|\rho - \eta \|_{H^{-s}}^2 \,
\Pi(d\rho,d\eta) \right)^{\frac k{d+2ks}}\\
&\le& C \MMM_k^{\frac d{d+2ks}}  \WWs(\alpha, \beta)
^{\frac{2k}{d+2ks}}
\eean
where we have successively used the inequality \eqref{estim:W1Hs}, H\"older
inequality, the definition of the moment of $\alpha$ and $\beta$, and Jensen
inequality.
\qed


\subsection{Quantified equivalence of chaos.}

This section is devoted to the proof of Theorem~\ref{Th:EquivChaos}, or more precisely, to the proof of the following accurate version of Theorem~\ref{Th:EquivChaos}. 

\begin{theo} \label{Th:EquivChaosAccurate}
  For any $G^N \in \Ps(E^N)$
and $f \in \PPP(E)$, there holds
 \bear
(i) \qquad \forall \, 1 \le j \le \ell \le N \qquad \Omega_j(G^N;f)  & \le& 2 \,  
\Omega_\ell(G^N;f),  \label{eq:EquivChaosPart1}\\
(ii) \qquad \forall \, 1 \le j \le N \qquad \Omega_j(G^N;f)  &\le &   
\Omega_\infty(G^N;f) + \frac{j^2}N. \label{eq:EquivChaosPart2}
\eear
 
For any $k>0$ and any $0 < \gamma <  \frac1{d+1 + \frac d k}$, there exists a
explicit constant  $C :=C(d,\gamma,k)$ such that 
\beqn \label{eq:EquivChaosPart3}
(iii) \qquad \Omega_\infty(G^N;f) \le C \, \MMM_k^{\frac 1k} \, \left(   \Omega_2(G^N;f) + {1
\over N} \right)^{\gamma},
\eeqn
where as usual $ \MMM_k := M_k(f) + M_k(G^N_1).$

For any $k>0$ and any $0 < \gamma <  \frac1{d' + \frac{d'}k}$, with $d' =
\max(d,2)$, there exists a constant $C := C(d,\gamma,k )$ such that 
\beqn \label{eq:EquivChaosPart4}
 (iv) \qquad | \Omega_N(G^N;f) - \Omega_\infty(G^N;f) |  \le C \, {  M_k(f)^{1/k}  \over N^\gamma}. 
\eeqn
\end{theo}

Let us make some remarks about the above statement. Roughly speaking, the two first inequalities  are in the good sense: the measure of chaos for a certain number of particles is bounded by the measure of chaos with more particles, and even in the sense of empirical measure (i.e. with $\Omega_\infty$). Let us however observe that 
the second inequality is meaningful only when the number $j$ of particles in the
left hand side is not too high, typically $j =o(\sqrt N)$.
The third inequality is in the "bad sense" and it is maybe the most important one, since it provides  an estimate of the measure of chaos in the sense of empirical measures by the measure of chaos for two particles only. It is for instance a key ingredient in \cite{MM-KacProgram}. See also corollary~\ref{cor:EquivChaos3} for versions adapted to probability measures with compact support or with exponential moment. The last inequality compares the measure of chaos at $N$ particles to its measure in the sense of empirical distribution.  It seems new and it will be a key argument in the next sections in order to make links between the  {\it Kac's chaos}, the 
 {\it entropy chaos} and the  {\it Fisher information chaos}.

\smallskip

\begin{rem} In the inequality~\eqref{eq:EquivChaosPart3},  
the $\Omega_2$ term in the right hand side may be replaced by any
$\Omega_\ell$ for $\ell \ge 2$, but it cannot be replaced by $\Omega_1$, which
does not measures chaoticity, as it is well known. We give a counter-example for
the sake of completeness. We choose $g$ and $h$ two distinct probability measures on
$E$, and take $f := {1
\over 2} (\, g + \, h)$. We consider the probability measure $G \in \PPE$,
and its associated sequence $(G^N)$ of marginal probability measures on
$\PPP(E^N)$ defined by 
$$
G =  {1 \over 2}( \delta_g +  \delta_h), \quad 
G^N := {1 \over 2} \, g^{\otimes N} + {1 \over 2} \, h^{\otimes N}.
$$
As $G_1 = f$, $\Omega_1(G^N,f) =0$ for all $N$, 
inequality~\eqref{eq:EquivChaosPart3} with $\Omega_2$ replaced by $\Omega_1$ will 
imply that $\Omega_\infty(G^N,f)$ goes to zero. But from
inequality~\eqref{eq:EquivChaosPart2} of 
Theorem~\ref{Th:EquivChaosAccurate}  
$$
 W_1(G^2,f^{\otimes 2}) = \Omega_2(G^N,f) \leq \Omega_\infty(G^N,f) + 
\frac C N.
$$
There is a contradiction since $G^2 \neq f^{\otimes 2}$ except if $g=h$.
\end{rem}

\medskip
We begin with some probably well known elementary inequalities and identities 
concerning Monge-Kantorovich-Wasserstein distances in space product.
For the sake of completeness we will nevertheless sketch the proofs of them.
  Remark that the two first formulas are particularly simple thanks to the
choice of the normalization \eqref{def:dEj}, and that they remains valid
if we replace $d_j$ by the normalized $l^1$-distance $\frac1j \sum_i 
|x_i-y_i|$.

\begin{prop}\label{prop:W1F1F} 

a) - For any $F^N,G^N \in \PPP_{\!sym}(E^N)$ and $1 \le j \le N$,  there hods
\beqn\label{prop:W1F1F-1}
 W_1(F^N_j,G^N_j) \le \Bigl( { j \over N} \, \Bigl[ {N \over j} \Bigr]
\Bigr)^{-1} \, W_1(F^N,G^N) \le  2  \, W_1(F^N,G^N).
\eeqn

\noindent
b) - For any $f,g \in P(E)$, there holds
\beqn\label{prop:W1F1F-2}
W_1(f^{\otimes N},g^{\otimes N}) =   W_1(f,g).
\eeqn

\noindent
c) - For any $f,g,h  \in P(E)$, there holds
\beqn\label{prop:W1F1F-3}
2 \, W_1(f \otimes h,g \otimes h) =   W_1(f,g).
\eeqn

\end{prop}

As a immediate corollary of  \eqref{prop:W1F1F-1} with $N:=\ell$,  $F^\ell :=
f^{\otimes \ell}$ and $G^\ell := G^N_\ell$, we obtain the first
inequality~\eqref{eq:EquivChaosPart1} of Theorem~\ref{Th:EquivChaosAccurate}.

As can be seen in the following proof, similar results also holds for MKW distances
constructed with arbitrary distance $D$ and exponents $p$, and therefore for the
$W_2$ distance. We do not state them precisely, but they will be useful in the
proof of the next Lemma~\ref{lem:WassNleqW1}. 

\medskip
\noindent
{\sc Proof of Proposition~\ref{prop:W1F1F}. } \\
 {\sl Proof of\eqref{prop:W1F1F-1}. } 
Consider $\pi \in \Pi(F^N,G^N)$ an
optimal transference plan in \eqref{def:W1}. 
Introducing the  Euclidean division,  $N = n \, j +r$,  $0 \le r \le j-1$, and writing $X = (X_1, ..., X_n, X_0) \in E^N$, $Y = (Y_1, ..., Y_n, X_0) \in E^N$,
with $X_i,Y_i \in E^j$, $1 \le i \le n$, $X_0, Y_0 \in E^r$, we have  
 \bean
W_1(F^N,G^N)
 &=& \int_{E^{2N}} d_{E^N} (X,Y) \, \pi(dX,dY)
 \\
&=&{1 \over N} \int_{E^{2N}} \left( \sum_{i=1}^n  j \,d_{E^j}(X_i,X_i)  +r \, d_{E^r}(X_0,Y_0)\right) \, \pi(dX,dY) 
 \\
  &\ge& {j \over N} \sum_{i=1}^n \int_{E^{2j}}d_{E^j} (X_i,Y_i)\, \tilde \pi_i(dX_i,dX_i),
 \eean
 with $\tilde \pi_i \in \Pi(\tilde F_i,\tilde G_i)$, where $\tilde F_i$ and $\tilde G_i \in \PPP(E^j)$ 
 denote the marginal probability measures of  $F^N$ and $G^N$ on the i-th block of 
variables. From the symmetry hypothesis, we have  $\tilde F_i = \tilde F_1 =
F^N_j$ and  $\tilde G_i = \tilde G_1 = G_j^N$ for any $1 \le i \le n$. As a
consequence, we have 
  \bean
\int_{E^{2j}} d_{E^j} (X_i,Y_i)\, \tilde \pi_i(dX_i,dX_i) \ge W_1(F^N_j,G^N_j),  
 \eean
 and we then deduce the first inequality in \eqref{prop:W1F1F-1}. Since the 
integer portion $n := [N/j]$ is larger than $1$, we have 
 $$
 { j \over N} \, \Bigl[ {N \over j} \Bigr] = {n \, j \over n \, j + r} \ge {n \, j \over n \, j + j} \ge {1 \over 2},
$$
from which we deduce the second inequality in \eqref{prop:W1F1F-1}. 

\smallskip\noindent
 {\sl  Proof of \eqref{prop:W1F1F-2}. }  We consider $\alpha \in \Pi(f,g)$ 
an optimal transference plan for the $W_1(f,g)$ distance and  we define the
associated transference plan 
  $\bar\pi := \alpha^{\otimes N} \in \Pi (f^{\otimes N},g^{\otimes N})$ by  
 $$
 \forall \, A_i,B_i \in E \quad 
 \bar\pi(A_1 \times ... \times A_N \times B_1 \times ... \times B_N) =
\alpha(A_1 \times B_1) \times ... \times \alpha(A_N \times B_N).
 $$
 By definition of $W_1(f^{\otimes N},g^{\otimes N})$, we then have 
  \bean
W_1(f^{\otimes N},g^{\otimes N}) &\le& {1 \over N} \sum_{i=1}^N \int_{E^{2N}}
d(x_i,y_i) \, \bar\pi(dX,dY) 
 =  W_1(f,g). 
 \eean
Since the first inequality in \eqref{prop:W1F1F-1}  in the case $j=1$ implies
the reverse inequality, the above inequality is an equality. 

\smallskip\noindent
 {\sl  Proof of \eqref{prop:W1F1F-3}. }  On the one hand, from the definition of the distance $W_1$ by transference plans, we have 
 for an optimal transference plan $\pi \in \Pi (f \otimes h,g \otimes h)$ the inequality
\bean
W_1 (f \otimes h,g \otimes h) 
&=& {1 \over 2} \int_{E^4} ( d_E(x_1,y_1) + d_E(x_2,y_2) ) \, \pi(dx_1,dx_2,dy_1,dy_2)
\\
&\ge& {1 \over 2} \int_{E^4} d_E(x_1,y_1)  \, \pi_1(dx_1,dy_1) \ge {1 \over 2} \, W_1(f,g),
\eean
since the $1$-marginal $\pi_1$ defined by $\pi_1(A \times B) = \pi(A \times E
\times B \times E)$ for any $A,B \in \BBB_E$ belongs to the transference plans set
$\Pi(f,g)$.  

On the other hand, considering an optimal transference plan $\pi \in \Pi (f,g)$ for
the $W_1$ distance, 
we define the associated transference plan $\bar \pi (dx,dy) := \pi(dx_1,dy_1) \,
\otimes  \, h(dx_2)  \delta_{y_2=x_2} \in \Pi (f \otimes h,g \otimes h)$,  and
we observe that 
\bean
W_1 (f \otimes h,g \otimes h) 
&\le& {1 \over 2} \int_{E^4} ( d_E(x_1,y_1) + d_E(x_2,y_2) ) \, \bar\pi(dx_1,dx_2,dy_1,dy_2)
\\
&=& {1 \over 2} \int_{E^4} d_E(x_1,y_1)  \, \pi(dx_1,dy_1) = {1 \over 2} \,
W_1(f,g). 
\eean
We obtain \eqref{prop:W1F1F-3} by gathering these two inequalities.
 \qed
  
\medskip
We next prove another lemma that allows to compare a distance between measures
on $\PPE$ and a distance between their marginals on $E^j$, and thus to compare
$\Omega_\ell$ and $\Omega_\infty$. 
\begin{lem}\label{lem:WassNleqW1}
For any distance $D$ on $E$ and $p \ge1$, extend $D$ on $E^j$ with
$D_{j,p}(V,W)^p = \frac1j
\sum_i D(v_i,w_i)^p$,  and define the associated MKW distance $W_{D_{j,p},p}$ on
$\PPP(E^j)$ and
the MKW distance  $\WW_{W_D,p}$ on $\PPP (\PPP(E))$ associated to $W_D$ and $p$.
Let
$\alpha$ and $\beta$ be two probability measures on $\PPE$. Then, for any $j \in
\N$,
\beqn \label{eq:WassNleqW1}
W_{D_{j,p},p} (\alpha_j,\beta_j) \leq \WW_{W_{D,p},p} (\alpha,\beta)
\eeqn
That is in particular true for the MKW distances $W_1$ and $W_2$ defined in
section~\ref{subsect:def}
$$
\forall j \in \N, \quad W_2(\alpha_j,\beta_j) \leq \WW_2 (\alpha,\beta), \qquad 
W_1(\alpha_j,\beta_j) \leq \WW_1 (\alpha,\beta).
$$
\end{lem}
\noindent
{\sc Proof of lemma~\ref{lem:WassNleqW1}.}  For simplicity we denote for any
$j$, $ W_{D_{j,p},p} =W_D$. We choose any transference
plan $\Pi$ between $\alpha$ and $\beta$ and write
\bean
\left[ W_D(\alpha_j,\beta_j) \right]^p & = & \left[ W_D \left(\int \rho^{\otimes
j}\,\alpha(d\rho),
\int \rho^{\otimes j} \,\beta(d\rho) \right) \right]^p \\
& = & \left[ W_D \left(\int \rho^{\otimes j} \,\pi(d\rho,d\eta), \int
\eta^{\otimes j}
\,\pi(d\rho,d\eta) \right)  \right]^p \\
& \le & \left[\int W_D \big( \rho^{\otimes j}, \, \eta^{\otimes j} \big)
\,\pi(d\rho,d\eta) \right]^p \\
& \le & \int [W_D ( \rho,\, \eta)]^p \,\pi(d\rho,d\eta),
\eean
where we have used the convexity property of the Wasserstein distance, the
equivalent of equality~\eqref{prop:W1F1F-2} in our general case, and Jensen
inequality. By optimisation on $\pi$ we obtain the claimed inequality.
\qed

\medskip

As a consequence of a classical combinatory trick, which goes back at least to
\cite{Grunbaum}, we have 

\begin{lem}[Quantification  of the equivalence $G^N_j \sim \hat G^N_j$] \label{lem:HSEquivQuant}
For any $G^N \in \Ps(E^N)$ and any $1 \le j \le1+  N/2$, we have 
$$
\|G^N_j - \hat G^N_j \|_{TV} \le 2 \,{j(j-1) \over N} \qquad 
\text{ and }\quad  W_1(G^N_j,\hat G^N_j) = {j(j-1) \over N},
$$
and in particular the first marginals are equal: $G^N_1 = \hat G^N_1$. 
\end{lem}

\medskip\noindent
{\sc Proof of Lemma~\ref{lem:HSEquivQuant}. } 
The second inequality is a straightforward consequence of the first inequality together with the use of
$$
W_1(G^N_j,\hat G^N_j) \le \frac12 \|G^N_j - \hat G^N_j \|_{TV}.
$$
A proof of the later may be found  in \cite[Proposition 7.10]{VillaniTOT}, in a slightly different context.
Here, the better factor $1/2$ can be obtained because of the stronger assumptions of our setting (the distance $d_{E^j}$ we deal with here is bounded by $1$).

\smallskip
The first inequality is a simple and classical combinatorial computation, see for instance \cite{Grunbaum}, \cite[Proposition 2.2]{S6}, \cite[Lemma 4.2]{MMWchaos} or \cite[Lemma~3.3]{MM-KacProgram}. We briefly sketch the proof for the convenience of the reader. 

For $1 \le j \le N$, we denote by $\CC_j^N$ the set of maps from $\{1,\ldots,j\}$ into  $\{1,\ldots,N\}$, and by $\AA_j^N$ the subset of $\CC_j^N$ made of the one-to-one maps. Remark that we have
$$
\bigl| \CC_j^N \bigr| = N^j,\qquad  \bigl| \AA_j^N \bigr| = \frac{N!}{(N-j)!}.
$$
Thanks to  the symmetry assumption made on $G^N$, we may write for any $\varphi \in C_b(E^j)$

\bean
 \langle  G^N_j , \varphi \rangle = \int_{E^N} \varphi(x_1,\ldots,x_j) G^N(dX) = \frac {(N-j)!}{N!} \sum_{s \in \AA_j^N} \int_{E^N} \varphi(x_{s(1)},\ldots,x_{s(j)}) G^N(dX)
\eean
From the definition of $\hat G^N_j$ we also get
\bean
 \langle \hat G^N_j , \varphi \rangle &=&
 \int_{\PPP(E)} \left( \int \varphi(y_1,\ldots,y_j) \rho^{\otimes j}(dY^j)  \right) \hat G^N(d\rho) \\
 & = &  \int_{E^N} \left( \int \varphi(y_1,\ldots,y_j) (\mu_X^N)^{\otimes j}(dY^j)  \right)  G^N(dX^N) \\
 &= & \frac1{N^j} \sum_{s \in \CC_j^N} \int_{E^N}  \varphi(x_{s(1)},\ldots,x_{s(j)})   G^N(dX^N).
\eean
The difference is then equals to
\begin{align*}
\langle  G^N_j - \hat G^N_j , \varphi \rangle = \Bigl( \frac {(N-j)!}{N!} - \frac1{N^j} \Bigr)
&  \sum_{s \in \AA_j^N} \int_{E^N} \varphi(x_{s(1)},\ldots,x_{s(j)}) G^N(dX) \\
& - \frac1{N^j} \sum_{s \in \CC_j^N \backslash \AA_j^N} \int_{E^N}  \varphi(x_{s(1)},\ldots,x_{s(j)}),G^N(dX^N)
\end{align*}
and may be bounded by
\bean
\bigl| \langle  G^N_j - \hat G^N_j , \varphi \rangle \bigr| &\le& 
\Bigl( 1 - \frac{N!}{N^j(N-j)! } \Bigr) \| \varphi\|_{L^ \infty} + \frac1{N^j}  \bigl| \CC_j^N \backslash \AA_j^N \bigr|  \| \varphi\|_{L^ \infty} \\
& = & 2 \Bigl( 1 - \frac{N!}{N^j(N-j)! } \Bigr) \| \varphi\|_{L^ \infty}.
\eean
 For $N \ge 2(j-1)$, we can bound the right hand side thanks to
\begin{eqnarray*}
  1- \frac{N! }{(N-j)! \, N^j} 
  &=& 
  1 - \left(1 - {1 \over N}\right) \, \cdots \, 
  \left(1 - {j-1 \over N}\right) 
  = 1 - \exp \left( \sum_{i = 0}^{j-1} 
    \ln \left(1 - \frac{i}N \right) \right) \\
  &\le& 1 - \exp \left( - 2 \sum_{i = 0}^{j-1} \frac{i}N \right) 
  \le 2 \sum_{i = 0}^{j-1} \frac{i}N \le {j(j-1) \over  N},
\end{eqnarray*}
where we have used 
\[
\forall \, x \in [0,1/2], \quad \ln ( 1 - x) \ge - 2 \, x \qquad \mbox{and}
 \qquad \forall \, x \in \R, \quad e^{-x} \ge 1 - x.
\]
We eventually get for $j \le 1+ N/2$
\[
\| \  G^N_j - \hat G^N_j \|_{TV} = \sup_{\| \varphi \|_\infty \le 1} \langle  G^N_j - \hat G^N_j , \varphi \rangle
\le 2\frac{j(j-1)}N,
\]
which ends the proof. \qed

\smallskip
Applying the previous lemmas~\ref{lem:WassNleqW1} and \ref{lem:HSEquivQuant}, we
can bound $\Omega_j$ by $\Omega_\infty$ and some rest. This is the second
inequality~\eqref{eq:EquivChaosPart2} of theorem~\ref{Th:EquivChaosAccurate}.

\smallskip\noindent
{\sc Proof of inequality~\eqref{eq:EquivChaosPart2} in 
Theorem~\ref{Th:EquivChaosAccurate}.} We simply write 
\bean
\Omega_j(G^N,f) = W_1(G^N_j,f^{\otimes j})  & \leq & W_1(G^N_j, \hat G^N_j)+ W_1(\hat G^N_j,f^{\otimes j}) \\
& \leq & \frac{j^2}N + \WW_1(\hat G^N, \delta_f) =  \frac{j^2}N + \Omega_\infty(G^N,f),
\eean
thanks to the two previous lemmas~\ref{lem:WassNleqW1} and
\ref{lem:HSEquivQuant}.
\qed

\medskip

 We establish now the key estimate which will lead to the third inequality \eqref{eq:EquivChaosPart3} in
Theorem~\ref{Th:EquivChaosAccurate} where $\Omega_\infty$ is controled by $\Omega_2$. 
Following \cite[Lemma 4.2]{MM-KacProgram}, the main idea is to use as an intermediate step the 
$H^{-s}$ norm on $\PPP(E)$, rather than the Wassertsein $W_1$ distance, because it is a 
monomial function of order two on $\PPP(E)$, and thus has a nice algebraic structure. 
This fact is stated in the following elementary lemma.

\begin{lem} \label{lem:H-s_poly} For $s>d/2$, define $\Phi_s : \R^d \rightarrow \R$ by 
\beqn \label{eq:defPhi}
\forall \, z \in \R^d, \qquad \Phi_s(z) := \int_{\R^d} e^{- i \,z \cdot \xi}\, \frac{d\xi}{\langle \xi
\rangle^{2s}}.
\eeqn
The function $\Phi_s$ is radial, bounded, and furthermore  if $s> \frac{d+1}2$,  it is Lipschitz. 
For any $\rho, \, \eta\in \PPP(E)$
\beqn \label{eq:H-s}
\|\rho - \eta \|_{H^{-s}}^2  = \int_{\R^{2d}} \Phi_s(x-y) \, (\rho^{\otimes 2} - \rho \otimes \eta ) (dx,dy)  
 + \int_{\R^{2d}} \Phi_s(x-y) \,  (\eta^{\otimes 2} - \eta \otimes \rho)  (dx,dy),
\eeqn
and   for any $\rho\in \PPP(E)$
$$
\|\rho \|_{H^{-s}}^2  = \int_{\R^{2d}} \Phi(x-y) \,\rho^{\otimes 2} (dx,dy),
$$
which means that the norm $H^{-s}$ on $\PPP(E)$ is the monomial function of
order two associated to the function $(x,y) \mapsto \Phi_s(x-y)$.
\end{lem}

\noindent
{\sc Proof of Lemma~\ref{lem:H-s_poly}. }  
We obtain that $\Phi_s$ is bounded  from the fact that $\int_{\R^d}  {\langle \xi
\rangle^{-2s}} \, {d\xi} $ is finite for $s > d/2$, and that it  is Lipschitz  from the fact that
 $\int_{\R^d}{\langle \xi \rangle^{1-2s}} \,  {d\xi}$ is finite when $s > (d+1)/2$.  
We now prove \eqref{eq:H-s}. Using the Fourier transform definition of the Hilbert norm of
$H^{-s}(\R^d)$, we have for any $\rho,\eta \in  H^{-s}(\R^d)$, 
and then for any $\rho,\eta \in \PPP(E) \subset  \PPP(\R^d) \subset H^{-s}(\R^d)$, 
\bean
\|\rho - \eta \|_{H^{-s}}^2 
&=& \int_{\R^d} (\hat \rho (\xi) - \hat \eta(\xi))
\, (\overline{\hat \rho(\xi) - \hat \eta(\xi)})  \, {d\xi \over \langle \xi
\rangle^{2s}} \\
 &=& \int_{\R^{3d}}(\rho (dx) - \eta(dx) \, (\rho(dy) - \eta(dy))  \, e^{- i \, 
 (x-y) \xi}\, {d\xi \over \langle \xi \rangle^{2s}}  \\
 &=& \int_{\R^{2d}} \Phi_s(x-y) \, (\rho^{\otimes 2} - \rho \otimes \eta ) (dx,dy)  
 + \int_{\R^{2d}} \Phi_s(x-y) \,  ( \eta^{\otimes 2} - \eta \otimes \rho)  (dx,dy). 
\eean
The last identity follows from \eqref{eq:H-s} by choosing $\eta = 0$. 
\qed

\medskip

Thanks to that Lemma, we will be able to obtain the following key estimate.

\medskip
\begin{prop}~\label{prop:WWenHs}  
For any $s > \frac{d+1}2$ there exists a constant $C  = 2\| 
\Phi_s\|_{Lip} \le \frac{2^{s+1} c_d}{2s-d-1} \in (0,\infty)$ (where $c_d$ denotes the surface of the unit sphere of $\R^d$)  
such that for any $G^N \in \PPP_{\!sym}(E^N)$, $N \ge 1$, $f \in \PPP(E)$,  
there holds
\beqn\label{ineq:WWenHs}
\WWs(\hat G^N,\delta_f)
\le C \, \left[W_1 ( \hat G^N_2, f \otimes f) \right]^{\frac12}.
\eeqn
\end{prop}

\noindent
{\sc Proof of Proposition~\ref{prop:WWenHs}. }  
Because $\PPP(E) \subset
\PPP(\R^d) \subset H^{-s}(\R^d)$ for $s>\frac d2$ and \\
\mbox{$ \Pi(\hat G^N,\delta_f) = \{ \hat G^N \otimes \delta_f \}$}, we have
$$
\left[\WWs (\hat G^N,\delta_f) \right]^2 := \inf_{\pi \in \Pi(\hat
G^N,\delta_f)} I[\pi] = I(\hat G^N \otimes \delta_f),
$$
with cost functional 
$$
 I[\pi]:= \int \!\! \int_{\PPP(E) \times \PPP(E)} \|\rho - \eta \|_{H^{-s}}^2
\, \pi(d\rho,d\eta). 
$$
Using Lemma~\ref{lem:H-s_poly}, we have  
\bean
I [\hat G^N \otimes \delta_f]
 &=&   \int_{\PPP(E)} \Bigl\{ \int_{\R^{2d}} \Phi_s(x-y) \, (\rho^{\otimes 2} - \rho \otimes f ) (dx,dy)
  \Bigr\}\, \hat G^N (d\rho) \\
&& +  \int_{\PPP(E)} \Bigl\{  \int_{\R^{2d}} \Phi_s(x-y) \,  ( f^{\otimes 2} - f \otimes \rho) 
(dx,dy) \Bigr\}\, \hat G^N (d\rho) \\
  &=& \int \!\! \int_{E^2}   \Phi_s(x-y) \, [ \hat G^N_2 (dx,dy) - 
\hat G^N_1(dx)  \, f (dy) ] \\
&& +  \int \!\! \int_{E^2}   \Phi_s(x-y) \,   [f(dx) \, f(dy)  
- f(dx)\,\hat G^N_1 (dy)].
\eean
Now we may bound the cost functional as follows: 
\bean
I [\hat G^N \otimes \delta_f]
&\le& \| \Phi_s\|_{Lip} \, \Big[W_1( \hat G^N_2, \hat G^N_1 \otimes f ) +  W_1(f
\otimes f, f \otimes \hat G^N_1) \Big]
\\
& \le & \| \Phi_s\|_{Lip} \, \Big[W_1( \hat G^N_2, f \otimes f ) +  2\,  W_1(f
\otimes f, \hat G^N_1 \otimes f) \Big] \\ 
&\le& \|\Phi_s\|_{Lip} \Big[W_1( \hat G^N_2, f \otimes f ) +   W_1(f  ,  \hat G^N_1 ) \Big] 
\\
&\le&   2 \|\Phi_s\|_{Lip} \,W_1( \hat G^N_2, f \otimes f ) , 
\eean
where we have used successively the  Katorovich-Rubinstein duality formula \eqref{def:W1KR},  the triangular inequality, 
the identity \eqref{prop:W1F1F-3}, and the first  inequality in \eqref{prop:W1F1F-1} together with the fact that $(\hat G^N_2)_1 = \hat G^N_1$.
\qed

\medskip
Putting together Proposition~\ref{prop:WWenHs}, Lemma~\ref{lem:HSEquivQuant}
above and Lemma~\ref{lem:ComparDistPPE} on comparaison of distances in $\PPE$,
we may prove inequality~\eqref{eq:EquivChaosPart3} of
Theorem~\ref{Th:EquivChaosAccurate}.

\smallskip
\noindent
{\sc Proof of  inequality~\eqref{eq:EquivChaosPart3} in Theorem~\ref{Th:EquivChaosAccurate}. } We define 
$s := \frac1{2\gamma} - \frac d{2k}$. Notice that $s > \frac{d+1}2 \ge
1$ thanks to the conditions satisfied by $\gamma$ and $k$. We can thus applied
the point $i)$ of Lemma~\ref{lem:ComparDistPPE},
Proposition~\ref{prop:WWenHs} and then Lemma~\ref{lem:HSEquivQuant} in order to
get 
\bean
\Omega_\infty(G^N;f)  & := & \WW_1(\hat G^N, \delta_f) \leq C(d,s) \MMM_k^{\frac d{d+2ks}} \WWs (\hat G^N, \delta_f)^{\frac{2k}{d+2ks}} \\
& \leq & \frac{C(d,s)}{2s-d-1} \MMM_k^{\frac d{d+2ks}} W_1( \hat G^N_2,
f^{\otimes 2})^{\frac k{d+2ks}} \\
& \leq & \frac{C(d,\gamma,k)}{\gamma^{-1} - d/k -d-1} \, \MMM_k^{\frac 1k}
\left( W_1( G^N_2, f^{\otimes 2}) + \frac2N
\right)^{\gamma},
\eean
since $\gamma =\frac k{d+2ks}$. This is the claimed inequality thanks to the
definition of $\Omega_2$. It is important to notice that the constant
$C(d,\gamma,k)$ of the last line depends on $d$, $k$ and $\gamma$ via
$s$. But as explained at the end of lemma~\ref{lem:ComparDistances}, it can be
choosen independent of $k$ and $\gamma$ if $s= \frac1{2\gamma} - \frac d{2k}$
remains in a compact subset of $\R^+$. 
\qed

\medskip
With stronger moment conditions on  the probability measures $f$ and $G^N$, we may
improve the exponent in the right hand side of \eqref{eq:EquivChaosPart3} and
therefore the rate of convergence to the chaos. Introducing the exponential
moment 
\beqn \label{def:ExpMoment}
 \forall \; F \in \PPP(E), \quad 
M_{\beta,\lambda} (F) := \int_E e^{\lambda |x|^\beta} \, F(dx),
\eeqn
$E=\R^d$, $\beta, \lambda > 0$, we have the following result. 

\begin{cor} \label{cor:EquivChaos3}
(i) There exists a constant $C=C(d)$ such that if the support of $f$ and  $G^N_1$ are both contained in the ball $B(0,R)$, for a positive $R$, then 
\beqn
\Omega_\infty(G^N;f) \le C \,R \, \left(   \Omega_2(G^N;f) + {1
\over N} \right)^{\frac1{d+1}} \left|  \ln  \left(   \Omega_2(G^N;f) + {1 \over N} \right) \right|.
\eeqn 
(ii) There exists a constant $C=C(d,\beta)$ such that if the $f$ and $G^N_1$ have bounded  exponential moment of order $M_{\beta,\lambda}$ 
for $\beta,\lambda > 0$,  there holds
 \beqn\label{eq:EquivChaosExpo}
\Omega_\infty(G^N;f) \le \frac C{\lambda^{\frac1 \beta}} \, K^{2(d+1)}  \, \left(   \Omega_2(G^N;f) + {1
\over N} \right)^{\frac1{d+1}} \left|  \ln  \left(   \Omega_2(G^N;f) + {1 \over N} \right) \right|^{1+\frac1\beta}
\eeqn 
where $K := \max(M_{\beta,\lambda} (f),M_{\beta,\lambda} (G^N_1))$. 
\end{cor}

\noindent
{\sc Proof of Corollary~\ref{cor:EquivChaos3}. }\\
{\sl Step 1. The compact support case.}
Here we simply have $M_k(f) \leq R^k$ and the same for the moments of $G^N_1$. 
Applying~\eqref{eq:EquivChaosPart3} with the explicit formula for the constant
$C$, we get for any $0 < \gamma < \frac1{d+1}$ and  $k > \frac d{\gamma^{-1} - d
-1}$
$$
\Omega_\infty(G^N;f) \le \frac{C(d,\gamma,k)}{\gamma^{-1} - d\,k^{-1}-d-1}  \,
R \, \left( \Omega_2(G^N;f) + {1\over N} \right)^{\gamma}.
$$
And we use the remark at the end of the previous proof that allows to replace 
$C(d,\gamma,k)$ by $C(d)$ if $s= \frac1{2\gamma} - \frac d{2k}$ is restricted
to some compact subspace of $[1,+\infty)$. It
will be the case in the sequel since we shall choose $k$ large and $\gamma$
close to $\frac1{d+1}$. Letting $k \to +\infty$ leads to
$$
\Omega_\infty(G^N;f) \le \frac{C(d)}{\gamma^{-1} -d-1}  \, R \, \left(  
\Omega_2(G^N;f) + {1\over N} \right)^{\gamma}.
$$
Denoting $\alpha := \frac1\gamma - d- 1$ and $a= \Omega_2(G^N;f) + {1\over N}$ which we assume  smaller than $\frac12$, the r.h.s can be rewritten 
$$
\Omega_\infty(G^N;f) \le C(d) \, \frac R  \alpha \,  a^{1/(d+1+\alpha)}.
$$ 
Some optimization leads to the natural choice $\alpha = 2\frac{(d+1)^2}{|\ln a|}$.  It comes
$$
\Omega_\infty(G^N;f) \le C(d) R \,  | \ln a | \,   a^{1/(d+1)} a^{1/(d+1+\alpha)- 1/(d+1)}.
$$
Since $\frac1{d+1} - \frac1{d+1+\alpha} \le \frac{\alpha}{(d+1)^2} \le \frac1{2 \,| \ln a |}$, we deduce
$$
a^{1/(d+1+\alpha)- 1/(d+1)} \le a^{- 1/(2 | \ln a |)} = e^{\frac12}
$$
and this concludes the proof of point $(i)$.

\medskip \noindent
{\sl Step 2. The case of exponential moment.} \\
Using the elementary inequality 
$x^k \le \bigl( \frac k{\lambda\beta e}\bigr)^{k/\beta} \, e^{\lambda \,
|x|^\beta}$, we get the following bound on the $k$ moment 
$$
M_k(F)^{1/k} \leq \Bigl( \frac k{\lambda \beta e} \Bigr)^{1/\beta} M_{\beta,\lambda}(F)^{1/k},  
$$
and it implies with our notations $\MMM_k^{\frac1k} \leq  \bigl( \frac {k}
{\lambda \beta e} \bigr)^{1/\beta}  (2K)^{1/k}$.
Applying~\eqref{eq:EquivChaosPart3} with the explicit formula for the constant
$C$ and the notation $a$ of the previous step, we get for any $0 < \gamma <
\frac1{d+1}$ and  $k > \frac d{\gamma^{-1} - d -1}$
$$
\Omega_\infty(G^N;f) \le \frac{C(d)}{(\lambda \beta e)^{1/\beta}} \frac {k^{1/\beta}} {\gamma^{-1} - d\,k^{-1}-d-1}  \, K^{1/k} \, a^\gamma.
$$  
Here we cannot take the limit as $k \to \infty$, but optimizing in $k$ the second fraction of the r.h.s, we choose $k$ satisfying $ \frac1\gamma-d-1 = \frac{2d}k$ and get the bound
$$
\Omega_\infty(G^N;f) \le \frac{C(d,\beta)}{\lambda^{1/\beta}} \frac{4d} {(\gamma^{-1} -d-1)^{1+1/\beta}}  \, K^{1/k} \, a^\gamma.
$$
Still denoting $\alpha = \frac1\gamma - d- 1 = \frac{2d}k$, the choice $\alpha = 2\frac{(d+1)^2}{|\ln a|}$ leads this time to the bound
$$
\Omega_\infty(G^N;f) \le  \frac{C(d,\beta)}{\lambda^{1/\beta}} \, K^{(d+1)/\ln 2 }\,  | \ln a |^{1+1/\beta} \,   a^{1/(d+1)},
$$
which concludes the proof.
\qed

\begin{rem} Inequality~\eqref{eq:EquivChaosPart3} in Theorem~\ref{Th:EquivChaosAccurate} says in particular that  for
any  $k > 0$ and $0 < \gamma < (d+1+d/k)^{-1}$  there exists a constant $C :=
C(d,\gamma,k)$ such that for any  $f \in \PPP(E)$, there holds
 \beqn \label{rem:Omega0FNf} 
  \Omega_\infty(f^{\otimes N};f) \le  {C  \, M_k(f)^{1/k}\over N^\gamma }.
\eeqn
\end{rem} 

For such a tensor product probability measures framework, the above rate can be improved in the following way. 

\begin{theo}[\cite{MM-KacProgram,BoissardLeGouic}] \label{theo:Omega0FNf} 
1. For a  moment weight exponent $k >0$ and an  exponent 
\begin{itemize}
\item[(i)] $ \gamma = \gamma_c := (2+1/k)^{-1}$ when $d=1$, 
\item[(ii)] $ \gamma \in (0,\gamma_c)$ with  $\gamma_c := (2+2/k)^{-1}$ when $d =2$,
\item[(iii)] $ \gamma = \gamma_c := (d+d/k)^{-1}$ when $d \ge 3$,
\end{itemize}
there exists a finite constant $C := C(d,\gamma,k)$ such that \eqref{rem:Omega0FNf} holds. 

2. Moreover, for any moment weight exponents $\lambda,\beta > 0$,  there exists a finite constant $C := C(d,\lambda,\beta, M_{\beta,\lambda}(f))$ such that
\beqn \label{rem:Omega0FNfsharp}    
\Omega_\infty(f^{\otimes N};f) \le  C \, {(\ln N)^{1/\beta}   \over N^{1/2} },  \text{ if } d = 1, \qquad 
\Omega_\infty(f^{\otimes N};f) \le  C \, {(\ln N)^{1+1/\beta}   \over N^{1/d} }, \quad 
\text{ if }  d\ge 2.
\eeqn

\end{theo}

On the one hand, using similar Hilbert norm arguments as  those used in the proof of Proposition~\ref{prop:WWenHs} and 
 inequality~\eqref{eq:EquivChaosPart3} in Theorem~\ref{Th:EquivChaosAccurate},
the first point in Theorem~\ref{theo:Omega0FNf} has been proved in \cite[Lemma 4.2(iii)]{MM-KacProgram} with however the restriction $\gamma \in (0,\gamma_c)$ 
when $d \ge 1$. The  optimal rate $\OO(1/N^{(2+1/k)^{-1}})$ in the critical case $\gamma = \gamma_c$, $d=1$, is not mentioned  in \cite[Lemma 4.2(iii)]{MM-KacProgram} but follows from a careful but straightforward reading of the proof of  \cite[Lemma 4.2(iii)]{MM-KacProgram}.   
The better rate obtained in Theorem~\ref{theo:Omega0FNf}  with respect to \eqref{rem:Omega0FNf} is due to the fact that 
for a tensor product measure one can work in the Hilbert space $H^{-s}$  with  $s>d/2$ rather than with  $s>(d+1)/2$ in the general case.
The second point in Theorem~\ref{theo:Omega0FNf} follows by adapting the proof of Corollary~\ref{cor:EquivChaos3} to this tensor product measures 
framework. 

\smallskip
On the other hand,  using matching techniques,  it has been proved in \cite{DubricYukich,BoissardLeGouic} that \eqref{rem:Omega0FNf} also holds true 
for the critical exponent $\gamma_c = 1/d$  in the compact support case (or exponential moment with $\beta=1$) when $d \ge 3$   and $\gamma_c = (d+d/k)^{-1}$ in the case of finite moment of order $k$
when $d \ge 3$. 
These last results thus slightly improve the estimates available thanks to our Hilbert norms technique. It is worth mentioning that  the critical exponents  are known to be optimal, see for instance \cite{DubricYukich,BartheBordenave}. A natural question is whether the rates in  inequality~\eqref{eq:EquivChaosPart3} and in Corollary~\ref{cor:EquivChaos3} may be improved using similar arguments as in \cite{DubricYukich,BoissardLeGouic}.

\medskip
We come to the proof of the last part of Theorem~\ref{Th:EquivChaosAccurate}, which will be a consequence of the following proposition 

\begin{prop}\label{prop:equivW1-2} For $F^N, G^N \in \Ps(E^N)$, there holds
\beqn\label{eq:equivW1-2}
W_1( F^N, G^N) = \WW_1 ( \hat F^N, \hat  G^N) .
\eeqn
\end{prop} 

\medskip\noindent
{\sc Proof of Proposition~\ref{prop:equivW1-2}. } We split the proof into two steps. 

\smallskip\noindent
{\sl Step 1.  A reformulation of the problem.
}
Since we are dealing with symmetric probability measures, it is  natural to introduce the equivalence relation $\sim$ in $E^N$ by saying that  $X=(x_1,
...,x_N),Y =(y_1, ..., y_N) \in E^N$ are equivalent,
we write $X \sim Y$,  if there exists a permutation $\sigma \in \SN$ such that $Y =
X_\sigma := (x_{\sigma(1)},\ldots,x_{\sigma(n)})$. 

We also introduce on $E^N$ the "semi"-distance $w_1$
\beqn\label{def:w1}
w_1 (X, Y)   :=   
\inf_{\sigma \in \SN} d_{E^N}(X,Y_\sigma)
= \inf_{\sigma \in \SN} {1 \over N} \sum_{i=1}^N d_E(x_i,y_{\sigma(i)}),
\eeqn
which only satisfies $w_1(X,Y) =0$ iff $ X \sim Y$. We then introduce the associated MKW functionnal $W^\dagger_1$. For $F^N, G^N \in \Ps(E^N)$, 
$$
W_1^\dagger (F^N,G^N)  :=  \inf_{\pi^N \in  \Pi (F^N,G^N)} \int_{E^N \times E^N}
w_1(X,Y) \, \pi^N(dX,dY).
$$
It is in fact a distance on the space of symmetric probability measures, but this point will also be a consequence of our proof. It is a classical result (see for instance \cite[Introduction. Example: the discrete case]{VillaniTOT}) that  
\beqn\label{eq:w1=W1}
\forall \, X,Y \in E^N, \quad W_1(\mu^N_X,\mu^N_Y) 
=  w_1(X,Y),
\eeqn
 (shortly, it means than we do not need to split the small Dirac masses when we try to optimize the transport between two empirical measures). We recall the notation $p_N$ defined in section~\ref{subsect:def} for the application that sends a configuration to the associated empirical measure : $p_N(X) = \mu^N_X$.
 
 Remark that its associated push-forward mapping restricted to the symmetric probability measures 
$$
\tilde p_N :  \Ps(E^N) \to \PPP(\PP_N(E)) \subset \PPE, 
\quad  G^N  \mapsto  \hat G^N := G^N_\# p_N,
$$
is  a bijection.  Its inverse can be simply expressed thanks to a dual formulation: for  $\alpha \in \PPP(\PP_N(E))$, its inverse $\tilde \alpha = \tilde p_N^{-1} \alpha$ is the probability measure satisfying
$$
\forall  \varphi \in C_b(E^N), \quad \int_{E^N} \varphi(X)\,\tilde \alpha(dX) = \int_{\PP_N(E)}  \tilde \varphi(\rho) \alpha(d\rho) , 
$$
where $\tilde \varphi (\rho) := \frac1{N!} \sum_{\sigma \in \SN} \varphi(X_\sigma)$, for any given $X$ such that $\mu = \mu^N_X$.
Similarly, defining $\PPP_{\!s,s}(E^N \times E^N)$ the subset of $\PPP(E^N \times E^N)$ of probability measures which are invariant under permutations
on the first and second blocks of $N$ variables separately, we have that 
$$
\tilde p_N^{\otimes 2} :  \PPP_{\!s,s}(E^N \times E^N)  \to \PPP(\PP_N(E) \times \PP_N(E))  , 
\quad  \pi^N  \mapsto  \hat \pi^N := \pi^N_\# (p_N,p_N),
$$
is a bijection. 

\noindent \smallskip
The identity~\eqref{eq:w1=W1} and the bijection $\tilde p_N$ allows us to establish the identity 
\beqn \label{eq:WW1=W1dagger}
\forall F^N G^N \in \PPP(E^N), \quad W^\dagger(F^N,G^N) = \WW_1(\hat F^N, \hat G^N).
\eeqn 
Indeed, denoting $\Pi_{s,s}(F^N,G^N) = \Pi(F^N,G^N) \cap \PPP_{\! s,s}(E^N,E^N)$, we have
\bean
W_1^\dagger (F^N,G^N)  & = &  \inf_{ \pi^N \in  \Pi_{s,s} (F^N,G^N)} \int_{E^N \times E^N} w_1(X,Y) \,\pi^N(dX,dY) \\
& = &  \inf_{ \pi^N \in  \Pi_{s,s} (F^N,G^N)} \int_{E^N \times E^N} W_1(p_N(X),p_N(Y)) \,\pi^N(dX,dY) \\
& = & \inf_{ \pi^N \in  \Pi_{s,s} (F^N,G^N)} \int_{\PP_N(E)  \times \PP_N(E)}  W_1(\rho,\eta) \,\pi^N_\#(p_N,p_N)(d\rho,d\eta) \\
& = & \inf_{ \hat \pi \in  \Pi (\hat F^N, \hat G^N)} \int_{\PPP(E) \times \PPP(E)}  W_1(\rho,\eta) \,\hat \pi(d\rho,d\eta) = \WW_1(\hat F^N, \hat G^N),
\eean
where we have essentially used  the invariance $w_1(X,Y) = w_1(X_\sigma,Y_\tau)$ for any $\sigma,\tau \in \SN$ 
and the fact that  $\tilde p_N^{\otimes 2}$ is a bijection.

\smallskip\noindent
{\sl Step 2.  The equality $W^\dagger_1 = W_1$.} 
The interest of the reformulation~\eqref{eq:WW1=W1dagger} is that we can now work on one space: $E^N$.
Remark that since $w_1(X,Y) \le d_{E^N}(X,Y)$, we always have $W^\dagger_1 \le W_1$,
and the equality will hold only if one transference plan  for $W_1^\dagger$ is concentrated on the set 
$$
\CC := \left\{ (X,Y)\in E^N \times E^N \;\text{ s.t}. \quad   w_1(X,Y)  = \inf_{\sigma \in \SN} 
d_{E^N}(X,Y_\sigma) = d_{E^N}(X,Y) \right\}.
$$
We choose an optimal transference plan $\pi$ for $W^\dagger_1$. For simplicity we
will assume that $\pi$ is symmetric, i.e. unchanged by the applications
$P_\sigma :(X,Y) \mapsto (X_\sigma,Y_\sigma)$ for any $\sigma \in \SN$. If not,
we replace it by its symmetrization $\frac1{N!} \sum_\sigma \pi_\# P_\sigma$
which will still be an optimal transference plan of $F^N$ onto $G^N$. Starting
from $\pi$, we will construct a transference plan $ \pi^* \in \Pi(F^N,G^N)$
such that  

- i) $\pi^*$ is concentrated on $\CC$. 

- ii) $I_N[\pi] = \int w_1(X,Y) \,\pi(dX,dY)  = \int w_1(X,Y) \, \pi^*(dX,dY) = I_N[ \pi^*]$

\noindent
Both properties imply then that 
\bean
W^\dagger_1(F^N,G^N)  &= & \int_{E^N \times E^N} w_1(X,Y) \,\pi(dX,dY)  = \int_{E^N \times E^N} w_1(X,Y) \, \pi^*(dX,dY) \\
& = &  \int_{E^N \times E^N} d_{E^N}(X,Y) \, \pi^*(dX,dY) \geq W_1(F^N,G^N)
\eean
which is the desired inequality.

We define $\pi^\ast$ in the following way. First, we introduce for any $X,Y \in E^N$
\bean
\CC_{X;Y}  & := & \left\{ Z \in E^N; \; Z \sim Y \text{ and } d_{E^N}(X,Z) = w_1(X,Y)   
 \right\}  \subset E^N
 \\
\rho_{X;Y} & := & \frac1{N_{X;Y}} \sum_{Z \in \CC_{X;Y}} \delta_{(X,Z)} \in \PPP(E^N \times E^N), 
\quad N_{X;Y} :=  \# \CC_{X;Y} \in \N^*. 
\eean
We note that $Z \in\CC_{X;Y} $ iff $ Z \sim Y$ and $(X,Z) \in \CC$, so that  $\Supp \rho_{X;Y} \subset \CC$.  
It can be shown that $(X,Y) \mapsto N_{X;Y}$ is a borelian application (it takes finite values and its level set are closed) and that 
$E^N \times E^N \to \PPP(E^N \times E^N)$,  $(X,Y) \mapsto \rho_{X;Y}$ is also borelian if $\PPP(E^N \times E^N)$ is endowed with the weak topology of measures.  This allows us to define a  transference plan $\pi^*$ by
$$
\pi^* :=  \int_{E^N \times E^N} \rho_{X;Y} \, \pi(dX,dY) \in \PPP(E^N \times E^N),
$$
or in other words, for any $\psi \in C_b(E^N \times E^N)$, we have 
\bean
\langle \pi^*,\psi \rangle 
&=& 
\int_{E^{2N}}  \frac1{N_{X;Y}} \sum_{Z \in \CC_{X;Y}}  \int_{E^{2N}} \psi(X',Y') \, \delta_{(X,Z)} (dX',dY')  \, \pi^N(dX,dY)
\\
&=& 
\int_{E^{2N}}  \frac1{N_{X;Y}} \sum_{Z \in \CC_{X;Y}} \psi(X,Z) \, \pi^N(dX,dY).
\eean
It remains to proof that $\pi^*$ satisfy the announced  properties. Since $\rho_{X;Y}$ is supported in $\CC$ for any $(X,Y) \in E^N \times E^N$, it is also the case for $\pi^*$.  It is also not difficult to show that the transport cost for $w_1$ is preserved. Indeed, we have 
\bean
\int_{E^{2N}} d_{E^N} (X',Y') \, \pi^*(dX',dY') 
&= &
\int_{E^{2N}} \left( \frac1{N_{X;Y}} \sum_{Z \in \CC_{X;Y}} d_{E^N} (X,Z) \right) \, \pi(dX,dY) \\
&= & 
\int_{E^{2N}} \left( \frac1{N_{X;Y}} \sum_{Z \in \CC_{X;Y}} w_1(X,Y) \right) \, \pi(dX,dY) \\
& = & 
\int_{E^{2N}} w_1(X,Y) \, \pi(dX,dY).
\eean
The fact that $\pi^*$ has first marginal $F^N$ is also clear since
for any $\varphi \in C_b(E^N)$
\bean
\int_{E^{2N}} \varphi(X') \,\pi^*(dX',dY') 
& = &
\int_{E^{2N}} \left( \frac1{N_{X;Y}} \sum_{Z \in \CC_{X;Y}} \varphi(X) \right) \, \pi(dX,dY) \\
& = & 
\int_{E^{2N}} \varphi(X) \,\pi(dX,dY)  = \int_{E^N} \varphi(X) \,F^N(dX) . 
\eean

For the second marginal, we shall use the following properties of $\CC_{X;Y}$ and $N_{X;Y}$  
$$
\forall \tau \in \SN, \quad 
Z_\tau \in \CC_{X_\tau;Y_\tau} \Leftrightarrow Z \in  \CC_{X;Y},\quad \text{and thus } \quad  N_{X_\tau;Y_\tau} =  N_{X;Y}. 
$$
Thanks to the invariance by symmetry of $\pi$ and $G^N$, we can write for any $\varphi \in C_b(E^N)$
\bean
\int_{E^{2N}} \varphi(Y) \,\pi^*(dX,dY) 
& = & 
\int_{E^{2N}} \left( \frac1{N_{X;Y}} \sum_{Z \in \CC_{X;Y}} \varphi(Z) \right) \, \pi(dX,dY) \\
& = &
 \frac1{N!} \sum_{\tau \in \SN} \int_{E^{2N}} \left( \frac1{N_{X_\tau;Y_\tau}} \sum_{Z \in \CC_{X_\tau;Y_\tau}} \varphi(Z) \right) \, \pi(dX,dY) \\
& = &
 \frac1{N!} \sum_{\tau \in \SN} \int_{E^{2N}} \left( \frac1{N_{X;Y}} \sum_{Z \in \CC_{X;Y}} \varphi(Z_\tau) \right) \, \pi(dX,dY) \\
& = &
\int_{E^{2N}} \left( \frac1{N_{X;Y}} \sum_{Z \in \CC_{X;Y}}  \tilde \varphi(Z) \right) \, \pi(dX,dY) \\
& = &
  \int_{E^{2N}} \tilde \varphi(Y) \,\pi(dX,dY) 
\\ 
&=& \int_{E^{2N}}  \tilde  \varphi(Y) \,G^N(dX) =  \int_{E^{2N}}  \varphi(Y) \,G^N(dX),
\eean
where we have introduced the symmetrization of $\varphi$ defined 
 by  $ \tilde \varphi (Z) := \frac1{N!} \sum_{\sigma \in \SN}
\varphi(Z_\sigma)  $ and we have used that  $\tilde \varphi (Z) = \tilde \varphi(Y)$ 
for any $Z  \in \CC_{X;Y}$ and the fact that $G^N$ is symmetric. This
concludes the proof.
\qed

\medskip
Putting together  Proposition~\ref{prop:equivW1-2} and \eqref{rem:Omega0FNf}, we obtain the inequality~\eqref{eq:EquivChaosPart4} of Theorem~\ref{Th:EquivChaosAccurate}.

 \smallskip\noindent
 {\sc Proof of inequality~\eqref{eq:EquivChaosPart4} in~\ref{Th:EquivChaosAccurate}. } 
 We have 
\bean
| \Omega_N(G^N,f) - \Omega_\infty(G^N,f) | & =& | W_1(G^N,f^{\otimes N}) - \WW_1(\hat G^N, \delta_f) |  \\
& = & | \WW_1(\hat G^N, \widehat{f^{\otimes N}}) - \WW_1(\hat G^N, \delta_f) | \\
& \le &  \WW_1(\widehat{f^{\otimes N}}, \delta_f) = \Omega_\infty(f^{\otimes N};f) \\ 
& \le &  {C \,M_k(f)^{1/k}\over N^\gamma },
\eean
where we have used  the definition of $\Omega_N$, $\Omega_\infty$,  the triangular inequality, Proposition~\ref{prop:equivW1-2} and \eqref{rem:Omega0FNf}.
\qed



\section{Entropy chaos and Fisher information chaos} 
\label{sec:Entropy}
\setcounter{equation}{0}
\setcounter{theo}{0}

In this section $E \subset \R^d$ stands for an open set or the adherence 
of a open space (so that the gradient of a function on $E$ is well defined).

\subsection{Entropy chaos}
The entropy of a probability measure on a compact subset of $\R^d$ with density $f
\,dx$ is well defined by the formula $\int f \ln f$. On a (possibly) unbounded
set $E$, we have to be more careful because the entropy may not be defined for
probability measure decreasing too slowly at infinity. This is a well known issue, but
we present here a rigourous definition for probability measures $F \in \PPP(E^j)$ having
 a finite moment $M_k$ for some $k>0$.  It will be usefull in the
section~\ref{sec:mixtures} where we define the level 3 entropy and  Fisher information on $\PPE$.

We emphasize that in the sequel we shall use the same notation $F$ for a
probability measure and its density $F\,dx$ with respect to the Lebesgue measure, when
the last quantity exists.
For any $k>0$ and $F \in \PPP_k(E^j) \cap L^1$, we define the
(opposite of the Boltzmann's) entropy  
\bear\label{def1EntropB}
 H_j(F) &:=& \int_{E^j} F \, \log F  
\\ \nonumber
&=&   \int_E h(F/G_k^j) \, G_k^j +
\int_E F \, \log
G_k^j \qquad(=:H_j^{(1)}(F))
\eear
with $G_k^j (V):= c^j_k \, \exp (- |v_1|^k - ... - |v_j|^k) \in \PPP(E^j)$,  $c_k$ chosen so that $G_i$ is a probability measure, and $h(s) := s \, \log
s - s + 1$. 
The RHS term is well defined in  $\R \cup \{+\infty \}$ 
as the sum of a nonnegative term and a finite real number, and it can be
checked that it is equal to the middle term, which has thus a sense. Next, we
extend the entropy functional to any  $F \in \PPP_k(E^j)$ by setting
\beqn\label{def3EntropB}
H_j(F) :=   \sup_{\phi_j \in C_b(E^j)} \Bigl\{\langle F, \phi_j \rangle - H^*
(\phi_j) \Bigr\}  + \int_E F \, \log
G_k^j \qquad(=:H^{(2)}(F))
\eeqn
where 
$$
H^* (\phi_j) :=  \int_{E^j}h^*(\phi_j) \, G^j_k
$$
and where $h^*(t) := e^t - 1$ is the Legendre transform of $h$.
Finally,  we define the normalized 
entropy functional $H$ by 
\beqn\label{defEntropBoltzPEj}
\forall \, F \in \PPP_k(E^j) \qquad
H(F) := {1 \over j} \, H_j(F).
\eeqn

We start recalling without proof  a very classical result concerning the entropy. 

\begin{lem}\label{lem:defENtropB} Let us fix $k > 0$. 
The entropy functional  $\PPP_k(E) \to \R \cup \{+\infty \}$, $\rho \mapsto
H_j(\rho)$ is well defined by the expression \eqref{def3EntropB}, 
is convex and is l.s.c. for the following notion of converging sequences: $\rho_n
\wto \rho$ in the weak sense of measures in $\PPP(E)$ 
and   $\langle \rho_n,|v|^{m} \rangle$ 
is bounded for some $m>k$ (the same holds of course for $H$). Moreover, 
$H_j(F)$ does not depend on the choice of $k$ used in the expression 
\eqref{def3EntropB}, 
$$
H(F) \ge \log c_k - M_k(F) \qquad \forall \, F \in \PPP_k(E),
$$
 and $H(F) < \infty$ iff $F \in  L^1$ 
, $F \, \log F \in L^1(E)$, and then $H(F) = H^{(1)}(F)$. 
\end{lem}

We also recall the definition of the (non-normalized) relative entropy between
two probability measures $\rho$ and $\eta$ of $\PPP(E^j)$ :
\beqn \label{def:EntropyRel}
H_j(\rho | \eta) :=  \int_{E^j} \ln \left( \frac{d\rho}{d\eta} \right)\,d\rho = 
\int_{E^j} (g \ln g +1 -g) d\eta
\eeqn
with $g = \frac{d\rho}{d\eta}$ if $\rho$ is absolutely continuous with
respect to $\eta$.  If $g$ is not defined, then $H_j(\rho|\eta) := + \infty$.
The associated normalized quantity is simply $H(\rho | \eta) := \frac1j
H_j(\rho | \eta)$. The relative entropy is
defined without moment assumption since the quantity under the last integral is
nonnegative. It can also be defined using a dual formula
similar to \eqref{def3EntropB}. For a
fixed $\eta$ it has the same properties as the entropy.

\smallskip
We now give two elementary and well known results which are fundamental for the analysis of the entropy defined on space product. 

\begin{lem}\label{lem:EntropBoltzHi&ii} On $\PPP_m(E^j)$, $m>0$, the entropy satisfies the identity 
\beqn\label{eq:EntropBoltzHi&ii}
\forall \, f \in \PPP_m(E) \qquad H(f^{\otimes j}) = H(f).
\eeqn
\end{lem}

\noindent
{\sc Proof of Lemma~\ref{lem:EntropBoltzHi&ii}. } 
If $f \in \PPP_m(E)$ is a function such that $H(f) < \infty$,  then we may use  \eqref{def1EntropB} as a definition and
$$
H (f^{\otimes j}) = {1 \over j} \int_{E^j}f^{\otimes j} \, \log
f^{\otimes j}=  \int_{E^j}  f^{\otimes j}(v_1, ..., v_j)  \, \log f (v_1) =
H_1(f).
$$
In the contrary,   $H_1(f) = \infty$ implies  $H_j(f^{\otimes j}) = \infty$. \qed

\begin{lem}\label{EntropRelative}  
(i) For any functions $f,g \in L^1_m(E) \cap \PPP(E)$, $m > 0$,  there holds
\beqn\label{ineg:fLOGffLOGg} 
H(f) := \int_E f \log f \ge \int_E f \log g, 
\quad \text{or} \quad H(f|g) := \int_E f \, \log (f/g) \geq 0,
\eeqn
with equality only if $f=g$ a.e..\\
(ii) More generally, for any nonnegative functions $f,g \in L^1_m(E)$, $m>0$,   there holds
$$
 \int_E f \log {f \over g} \ge F \log {F \over G}, \quad
 \hbox{with}\quad F := \int_E f, \,\,\, G := \int_E g.
$$
(iii) A consequence of $(i)$ is that if $F \in \PPP(E^j)$ has first marginal $f$ with $H(f) < + \infty$, 
then
$$
H(F) \ge H(f) \qquad \text{with equality only if} \,\, \,  F= f^{\otimes j} \;
\text{a.e.}.
$$
(iv) The entropy is superadditive: for any   $F \in \PPP_m(E^{i+j}) \cap
\PPP_{\!sym}(E^{i+j})$, $i, j \in \N^*$, $m > 0$, 
the following inequality holds  
\beqn\label{ineg:additiviteEntrop} 
\qquad H_{i+j} (F_{i+j}) \ge H_{i} (F_{i}) + H_{j} (F_{j}),  \quad \text{(non-normalized entropy)} ,
\eeqn
 where $F_\ell$ as usual stands for the $\ell$-th marginal of $F$. 
\end{lem}

\noindent
{\sc Proof of Lemma~\ref{EntropRelative}. } (i) 
To obtain the inequality, write $H(f|g) = \int h(f/ g) f$ and
use the fact that $h(s)= s \log s - s + 1$ is a nonnegative function. 
Next there is equality only if $h(f/g)=0$ a.e. on $\{f >0\}$. Since $h$ vanishes
only at $s=1$, it means that $f=g$ a.e. on $\{ f>0 \}$. Using that
$\int f = \int g =1$, we obtain the claimed equality. 

\smallskip\noindent (ii)  We write
$$
 \int_E f \log {f \over g} = F \int_E f/F \log {f/F \over g/G} + \int_E f \log {F \over G},
 $$
 the first term is nonnegative  thanks to \eqref{ineg:fLOGffLOGg} and the second term is the one which appears on the RHS of the 
 claimed inequality.
 
\smallskip\noindent  (iii) We use the first inequality~\eqref{ineg:fLOGffLOGg}
on $E^j$ with $F$ and $f^{\otimes j}$
$$
H(F) = \frac1j \int_{E^j} F \log F \ge \frac1j \int_{E^j} F \log f^{\otimes j} =
\int_{E^j} F(V) \log f(v_1) \,dV = H(f).
$$
Using again the point $i)$, we see that equality can occur only if $F =
f^{\otimes j}$ a.e..

\smallskip\noindent  (iv) Denote $h_\ell :=  H_\ell (F_\ell)$. 
If $h_{i+j} = +\infty$ there is nothing to prove.
Otherwise, we have $h_{i+j} < \infty$ which in turn implies
 $F \in L^1(E^{i+j})$, then $F_i \in L^1(E^i)$, $F_j \in L^1(E^j)$, so that the entropy may be defined thanks to  
\eqref{def1EntropB}.
In $\R \cup \{-\infty \}$, we compute 
\bean
h_{i+j} - h_i -  h_j 
&= & \int_{E^{i+j}} F_{i+j} \log F_{i+j} 
\\
&& - \int_{E^{i+j}} F_{i+j} \log F_{i} (v_1, .., v_i)  - \int_{E^{i+j}} F_{i+j} \log F_{j} (v_{i+1}, .., v_{i+j})
\\
&= & \int_{E^{i+j}} F_{i+j} \log F_{i+j} - \int_{E^{i+j}} F_{i+j} \log F_{i} \otimes F_j \ge 0,
\eean
thanks to  \eqref{ineg:fLOGffLOGg}.  \qed

\medskip
Our first result shows that entropy chaos is a stronger notion than Kac's chaos.

\begin{theo}[Entropy and chaos]\label{theo:ChaosEntrop}
Consider $(G^N)$ a sequence of $\PPP_{\!sym}(E^N)$ such that $\langle G^N_1, |v|^m \rangle \le a$ for any $N \ge 1$ and for some fixed $m, a > 0$
and consider  $f \in \PPP(E)$. 

\smallskip
{\bf 1)} If  $G^N_j \wto F_j$ weakly in $\PPP(E^j)$ for some given $j \ge 1$, then  
\beqn\label{ineq:HFNtoHf}
H(F_j) \le \liminf  H(G^N).
\eeqn
In particular, when $(G^N)$ is $f$-Kac's chaotic, \eqref{ineq:HFNtoHf} holds for any $j \ge 1$ with $F_j := f^{\otimes j}$. 

\smallskip\noindent
{\bf 2)}  On the other way round, if $(G^N)$ is $f$-entropy chaotic, then $(G^N)$ is $f$-Kac's chaotic.
\end{theo}

\noindent
{\sc Proof of Theorem~\ref{theo:ChaosEntrop}. } {\sl Step 1. } For any $N \ge j$ we
introduce the  Euclidean decomposition $N = n \, j + r$, $0 \le r \le j-1$, exactly as in the proof of Proposition~\ref{prop:W1F1F}.
Iterating $n$ times the superadditivity inequality \eqref{ineg:additiviteEntrop} we have
$$
H_N(F^N) \ge n \, H_j(F^N_j) +  H(F^N_r),
$$
with the convention $H(F^N_r) = 0$ when $r=0$.
We get \eqref{ineq:HFNtoHf} by passing to the limit in that inequality divided by $N$, using 
that $H$ is l.s.c. and that $H(F^N_r)$ is bounded by below thanks to
Lemma~\ref{lem:defENtropB} and the condition on the moment.

\smallskip\noindent
{\sl Step 2. }We assume that $(G^N)$ is $f$-entropy chaotic, that is 
$$
G^N_1 \wto f \,\, \hbox{weakly in} \,\, \PPP(E) 
\quad\hbox{and}\quad
H(G^N) \to H(f)  < \infty.
$$
Let us fix $j \ge 1$.  The sequence $(G^N_j)$ being bounded in  $\PPP_m(E^j)$, 
there exists  $F_j \in \PPP(E^j)$ and a subsequence
 $(G^{N'})$ such that $G^{N'}_j \wto F_j$ weakly in  $\PPP(E^j)$. Thanks to step 1, we have 
$$
H(F_j) \le \liminf H(G^{N'}_j) \le   \liminf H(G^{N'}) = H(f) = H(f^{\otimes j}).
$$
Since the first marginal of $F_j$ is $(F_j)_1 = \lim_{N \rightarrow
+\infty} G_1^N = f$, the third point of Lemma~\ref{EntropRelative}  gives that
$F_j = f^{\otimes j}$ a.e..
As a conclusion and because we have identified the limit, we have proved that the
all sequence $(G^N_j)$  weakly converges to $f^{\otimes j}$.
\qed

\subsection{Fisher chaos}

We now establish similar results for the Fisher information functional. 
For an arbitrary probability measure $G \in \PPP(E^j)$, we define the normalized Fisher
information by
\beqn\label{def:Fisher1}
I^{(1)}_j(G) := 
 \left\{ \begin{array}{ll}  {\displaystyle \int_{E^j} {|\nabla G |^2 \over G} 
=  \int_{E^j} |\nabla \ln G |^2 \, G } \in \R \cup \{+\infty \}   & \text{if }
G \in W^{1,1}({E^j}), \\
 +\infty &   \text{if } G \notin W^{1,1}({E^j}),  \end{array} \right.
\eeqn

For $G \in \PPP({E^j})$, we also give an alternative definition
\beqn\label{def:Fisher2}
I^{(2)}_j(G) :=  \sup_{\psi \in C_b^1({E^j})^d}  \langle G,  - {| \psi
|^2 \over 4} - \hbox{div}\, \psi\rangle \in \R \cup \{+\infty \} .
\eeqn

\begin{lem}\label{lem:Fisher1} For all $j \in  \N$, the identity  $I^{(1)}_j =
I^{(2)}_j$ holds 
on $\PPP(E^j)$, and we simply denoted by $I_j$ the usual (non-normalized) Fisher information
and by $I = j^{-1} \, I_j$ the normalized Fisher information. 
The functionals $I_j$ and $I$
are proper, convex, l.s.c. (in the sense of the weak convergence  of measures) 
on $\PPP(E^j)$. 
\end{lem}

\noindent
{\sc Proof of Lemma~\ref{lem:Fisher1}. } For the sake of simplicity, we only 
deal with the case $j=1$. We split the proof into two steps. 

\smallskip\noindent
{\sl Step 1. }  Assume that $f \in W^{1,1}$. Since for all $ \psi
\in C^1_b(E)^d$
$$|\nabla \ln f|^2 - \nabla \ln f \cdot \psi + \frac{|\psi |^2}4 =
\left| \nabla \ln f - \frac\psi 2\right|^2 \ge 0,
$$
we have
$$
I^{(1)}(f) =  \int_E | \nabla \ln f |^2 f \geq
\int_E \left( \nabla \ln f  \cdot \psi - \frac{|\psi |^2}4 \right)\,
f .
$$

For any sequence $(\psi_n)$ of smooth functions approximating $2\nabla
\ln f = 2\frac{\nabla f}f$, we obtain that
\bear
I^{(1)}(f) 
&=&  \lim_{n\to\infty} \int_E \left( \nabla \ln
f \cdot \psi_n - \frac{|\psi_n |^2}4 \right)\, f  \nonumber
\\
&=&  \sup_{\psi \in C^1_b(E)^d} \int_E \left( \nabla \ln
f \cdot \psi - \frac{|\psi |^2}4 \right)\, f  \nonumber
\\
& = & \sup_{ \psi \in C^1_b(E)^d} \int_E  \left[ \nabla f \cdot \psi - f \,
\frac{|\psi|^2}4 \right] \qquad  =: I^{(3)}(f) \label{def:I3}.
\eear
The remaining equality $I^{(3)} = I^{(2)}$ is just a simple integration by
parts. 
 Remark that maximizing sequences $(\psi_n)$ must converge
(up to some subsequence) pointwise  to $2 \, \nabla  \ln f$ a.e. on $ \{ f
\neq 0 \}$. We shall use that point in the sequel.

We also remark that this reformulation $I^{(2)}$ is also exactly
the one obtained when using the general Fenchel-Moreau theorem on the convex
function $(a,b) \to \frac{|b|^2}a$ (which is used in the integral defining
$I^{(1)}$).

\smallskip\noindent
{\sl Step 2. }
It remains to check that the
equality $I^{(1)}=I^{(2)}$ is also true on $ \PPP(E) \backslash W^{1,1}(E)$.
In other words that if $f \notin W^{1,1}(E)$ then $I^{(2)}(f) = +\infty$.
In what follows, we prove the contraposition : $I^{(2)}(f) < +\infty$ implies
$f \in W^{1,1}(E)$. Once it will be done, we will have $I^{(1)} =I^{(2)}$
everywhere, from what follows that $I$ is l.s.c. in the
sense of the weak convergence of measures. 

Consider $f \in \PPP(E)$ and assume $I^{(2)}(f) < \infty$. We
deduce that for any $\psi \in C^1_b(E)^d$ and any $t \in \R$
$$
 \int_{E}  f \, [- t^2 \, {|\psi|^2 \over 4} - t \, \hbox{div} \, \psi   ] \le
I^{(2)}(f),
$$
so that by optimizing in $t \in \R$ and using that  $f \in \PPP(E)$, we get
$$
\forall \, \psi \in C^1_b(E)^{d} 
\qquad \left|\int_{E}  f \,   \hbox{div}\, \psi  \right|^2\le 4 \, I^{(2)}(f) 
\,   \int_{E}  f \,  {|\psi|^2 \over 4} \le I^{(2)}(f) \,   \|\psi
\|_{L^\infty}^2.
$$
That inequality implies  $f \in BV(E)$ and $\|\nabla f \|_{TV} \le
\sqrt{I^{(2)}(f)}$. Using that $f \in BV(E)$ and making an integration by part
in the definition
of $I^{(2)}(f)$, we find  
$$
I^{(2)}(f) =  \sup_{ \psi \in C^1_b(E)^d} \int_{E}  [\nabla f \cdot \psi  -
f \, {|\psi|^2 \over 4} ]  =  I^{(3)}(f).
$$
Now, for any compact subset $K \subset E$ with zero Lebesgue measure, we may 
find a sequence $\rho_\eps \in C^1_c (E)$ such that  $0 \le \rho_\eps  \le
1$, $\rho_\eps = 1$ on $K$ and $\rho_\eps \to 0$ a.e., so that for any
$t > 0$ and using that  $f \in BV(E) \subset L^1(E)$, we get for all $\eps >0$
\begin{align*}
t \, \int_{K} |\nabla f | & \le  t \, \int_E |\nabla f | \rho_\eps \\
& \le \sup_{\psi \in C^1_c(E)^d, \| \psi \|_\infty \le 1}
   \int_{E} \nabla f \cdot \psi \,t \rho_\eps \\
& \le \sup_{\psi \in C^1_c(E)^d, \| \psi \|_\infty \le 1}
   \int_{E} \Bigl[ \nabla f \cdot \,t \psi  \rho_\eps  -
   f t^2 \frac{|\psi|^2 \rho_\eps^2}4 \bigr] + \frac{t^2}4 \int_E f \rho_\eps^2 \\
& \le I^{(3)}(f) + \frac{t^2}4 \int_E f \rho_\eps^2.
\end{align*}
Passing to the limit $\eps \to 0$ using that $f \in L^1(K)$ and then $t \to \infty$, we deduce that  $\nabla f$ vanishes on 
$K$, which precisely means that  $\nabla f$ is a measurable function. We have
proved $f \in W^{1,1}(\R^d)$.
\qed

\smallskip
Similarly, we define for two measures $\rho$ and $\eta$ on $E^j$ their 
(non-normalized) relative
Fisher information $I(\rho|\eta)$ by
\beqn \label{def:FisherRel}
I_j(\rho|\eta) :=  \int_{E^j} \frac{|\nabla g|^2} g d\eta = \int_{E^j}
\left|\nabla \ln \frac{d\rho}{d\eta}\right|^2  d\rho, 
\eeqn
where $g = \frac{d\rho}{d\eta}$ if $\rho$ is absolutely continuous with respect
to $\eta$. If not, $I_j(\rho|\eta) := + \infty$. The associated normalized
quantity is simply $I(\rho|\eta) := \frac1j I_j(\rho|\eta)$. For a fixed $\eta$,
the relative Fisher information has roughly the same properties as the Fisher
information. In particular, if $\eta$ has a derivable density, we have the
equality
\beqn \label{eq:Fisherrelativedual}
I_j(\rho|\eta)  = \sup_{\varphi \in C^1_b(E^j)^{dj}} \int_{\R^j} \left( 
 - \varphi \cdot \frac{\nabla \eta}{\eta} -  \diver \varphi - \frac{|\varphi|^2}4
\right) \,d\rho.
\eeqn

\begin{lem}\label{lem:FisherCondii} For any $f \in \PPP(E)$ there holds
$I(f^{\otimes j}) = I(f)$. 
\end{lem}

\noindent
{\sc Proof of Lemma~\ref{lem:FisherCondii}. }
If $I(f) < \infty$ then $f \in
W^{1,1}(E)$ and also $f^{\otimes j}\in W^{1,1}(E^j)$. 
The following computation is then meaningful  
\bean
I(f^{\otimes j}) =  \frac1j \int_{E^j} {|\nabla_{E^j} f^{\otimes j}|^2 \over
f^{\otimes j}} 
= \int_{E^j} {|\nabla_E f |^2 \over   f} \otimes  f^{\otimes (j-1)}  = 
I(f).
\eean
Since  $I_j(f^{\otimes j}) < \infty$ implies $f^{\otimes j} \in W^{1,1}(E^j)$ 
and then $f \in W^{1,1}(E)$, we also have $I_j(f^{\otimes j}) = j \, I(f)$ if
$I(f) = \infty$. 
\qed

\begin{lem}\label{lem:FisherCondiv}
For any $F \in P_{sym}(E^j)$ and $1 \le \ell \le j$, then holds\\
(i)   $I(F_\ell) \leq I(F)$.\\
(ii) The Fisher information is super-additive. It means that 
\beqn \label{Fishersuperadd}
I_j(F) \ge  I_\ell (F_\ell) +  I_{j-\ell}(F_{j-\ell}), \quad \text{(non-normalized Fisher information)} ,
\eeqn
with in the case $I_\ell (F_\ell) +  I_{j-\ell}(F_{j-\ell}) < + \infty$ equality  only if $F = F_\ell \otimes F_{j-\ell}$.\\
\smallskip
(iii) If $I(F_1) < +\infty$, the equality  $I(F_1) = I(F)$ holds if and only if $F =
(F_1)^{\otimes j}$.

\end{lem}

\noindent
{\sc Proof of Lemma~\ref{lem:FisherCondiv}. }

\smallskip\noindent
{\sl Proof of (i). }
If $I(F) = +\infty$ the conclusion is clear. Otherwise, thanks to the equivalent definition $I^{(3)}$ of the Fisher information and the symmetry assumption of $F$, we have 
\bean
I(F) &=&     \sup_{\psi \in C_b(E^j)^{dj}}  {1 \over j} \int_{E^j}  \Bigl( \psi(x_1,\ldots,x_j) \cdot \nabla F
- F \, {|\psi(x_1, \ldots, x_j)|^2 \over 4} \Bigr) 
\\
&=&     \sup_{\psi \in C_b(E^j)^d}  \int_{E^j}  \Bigl( \psi(x_1,\ldots,x_j) \cdot \nabla_1 F
- F \, {|\psi(x_1, \ldots, x_j)|^2 \over 4} \Bigr) 
\\
&\ge&  \sup_{\psi \in C_b(E^\ell)^d}  \int_{E^j}  \Bigl( \psi(x_1,\ldots,x_\ell) \cdot \nabla_1 F - F \, {|\psi(x_1,\ldots,x_\ell)|^2 \over 4} \Bigr) 
\\
&=&  \sup_{\psi \in C_b(E^\ell)^d}  \int_{E^\ell}  \Bigl( \psi  \cdot \nabla_1
F_\ell - F_\ell \, {|\psi |^2 \over 4} \Bigr)  
= I(F_\ell). 
\eean
 
 \noindent{\sl Proof of the superadditivity property (ii).} 
The first proof of that result seems to be the one by Carlen in
\cite[Theorem 3]{Carlen1991}. We sketch now another proof that uses the
third formulation $I^{(3)}$. We recall that in the definition of $I^{(3)}_j(F)$ the supremum is taken over the $\psi =(\psi_1,\ldots,\psi_j)$, with all 
$\psi_i : E^j \to \R^d$. We now restrict the supremum over the $\psi$ such that:

- The $\ell$ first $\psi_i$ depend only on $(x_1,\ldots,x_\ell)$, with the
notation $\psi^\ell = (\psi_1, \ldots,\psi_\ell)$.

- The $(j- \ell)$ last $\psi_i$ depend only on $(x_{\ell+1},\ldots,x_j)$,
with the notation $\psi^{j-\ell} = (\psi_{\ell+1}, \ldots,\psi_j)$.

We then have the inequality
\bean
I_j(F) & \geq &  
\sup_{\psi^\ell ,\, \psi^{j-\ell} } 
\int_{E^j} [\nabla_\ell f \cdot \psi^\ell + \nabla_{j-\ell} f \cdot
\psi^{j-\ell} - f \, \frac{|\psi^\ell|^2+ |\psi^{j-\ell}|^2}4 ] \\
& = &  \sup_{\psi^\ell \in C^1_b(E^\ell)^{\ell d}} 
\int_{E^\ell} [\nabla f_\ell \cdot \psi^\ell - f_\ell \, \frac{|\psi^\ell|^2}4
]\\ 
&& \hspace{3cm} + 
\sup_{\psi^{j-\ell} \in C^1_b(E^{j-\ell})^{(j-\ell)d}} 
\int_{E^{j-\ell}} [\nabla f_{j-\ell} \cdot \psi^{j-\ell} 
- f_{j-\ell} \, \frac{|\psi^{j-\ell}|^2}4 ] \\
& = & I_\ell(F_l) + I_{j-\ell}(F_{j-\ell})
\eean
If the inequality is an equality, we use the remark made at the end of Step 1 in the proof 
of Lemma~\ref{lem:Fisher1} : Maximizing sequences $\psi^\ell_n$
and
$\psi^{j-\ell}_n$ for respectively $I_\ell$ (resp. $I_{j-\ell}$) should converge
pointwise towards $2 \, \nabla \ln f_l$ (resp. $2 \, \nabla \ln f_{j-\ell}$) up
to some subsequence, a.e. on $\{ f_\ell \neq 0\}$ (resp.  $\{ f_{j - \ell} \neq
0 \}$). If we have equality, we also must have $(\psi^\ell_n,\psi^{j-\ell}_n)
\to 2 \,  \nabla  \ln f$ on $\{ f\neq 0\}$, a set that is included in $\{
f_\ell \neq 0\} \times \{ f_{j - \ell} \neq 0\}$ and thus
$$
\nabla \ln f = ( \nabla \ln f_\ell ,\nabla \ln f_{j-\ell}) = \nabla \ln (f_\ell
\otimes f_{j-\ell}),
$$
which implies the claimed equality since $f$ and $f_\ell
\otimes f_{j-\ell}$ are probability measures.

\smallskip\noindent{\sl The case of equality (iii).} 
Using recursively the superadditivity in that particular case, we get with the
notation $F_1=f$
$$ 
I(f) = I(F) \ge \frac {j-1}j I(F_{j-1}) + \frac 1j I(f)  \ge \frac {j-2}j
I(F_{j-2}) + \frac 2j I(f) \ge \ldots \ge I(f).
$$
Therefore, all the inequalities are equalities. We obtain that
$$
F = F_{j-1} \otimes f = F_{j-2} \otimes f \otimes f = \ldots = f^{\otimes j},
$$
by applying recursively the case of equality in \eqref{Fishersuperadd}.
\qed

\medskip\smallskip
 
It is classical and essentially a consequence of the Sobolev inequality and  the Rellich-Kondrachov Theorem  (together with very standard manipulations on the entropy functional which are similar to the 
ones presented at the end of the proof of Theorem~\ref{theo:EntropCvgceKac})  that for $(f_n)$ a sequence of  $ \PPP(E)$, the conditions 
$$
f_n \wto f \,\,\hbox{weakly in }\PPP(E), \quad M_k(f_n) \,\, \hbox{bounded}, \,\, k > 0,  \quad \hbox{and} \quad I(f_n) \le C 
$$
imply that  $H(f_n) \to H(f)$. 
A natural question is whether a similar result holds for a sequence $(F^N)$ in $\PPP(E^N)$.
 Before answering affirmatively to that question, we establish a normalized non-relative HWI inequality for 
 a large class of  sets $E \subset \R^d$. It is a variant of  the famous HWI
inequality of Otto-Villani \cite{OttoVillani} that will be the cornerstone of the argument. Let us
mention that its good behaviour in any dimension is of particular importance
here and it is due to the good (separate) behaviours of $H$, $W_2$ and $I$ with respect to  the dimension.

  \begin{prop}\label{prop:HWI}  
   Assume that  $E \subset \R^d$ is a bi-Lipschitz volume preserving deformation of a convex set of
$\R^d$, $d \ge 1$: there exists a convex subset $E_1 \subset \R^d$  and a bi-lipschitz diffeomorphism $T : E_1 \to E$ which preserves the volume (i.e. its Jacobian is always equal to $1$). Then, the  normalized non relative HWI inequality holds in $E$: there exists a constant $C_E \in [1,\infty)$  such that
\beqn\label{eq:GalHWI} 
\forall \, F^N,\, G^N \in \PPP_2(E^N) \qquad H(F^N ) \le  H(G^N) +   C_E \, W_2(F^N,G^N) \,\sqrt{ I (F^N)} .
\eeqn
More precisely, the above inequality holds with $C_E := \| \nabla T \|_\infty \, \| \nabla T^{-1} \|_\infty$ where \break
 $\| \nabla T \|_\infty := \sup_{v \in E} \sup_{|h|_2 \leq 1} | \nabla T(v) \, h  |_2$. 
 \end{prop}

Before going to the proof, remark that the class of set $E$ which are bi-Lipschitz volume preserving deformation of convex set is rather large. For instance, it is shown in~\cite[Theorem 5.4]{FonsecaParry}
 that any star-shaped bounded domain with Lipschitz boundary (and some additional assumptions) is in the previously mentioned class.

\medskip\noindent 
{\sc Proof of Proposition~\ref{prop:HWI}. } We proceed in three steps. 

\smallskip
\noindent 
{\sl Step 1. $E=\R^d$. } Let us first recall the famous HWI inequality of Otto-Villani. 
Consider $\rho = e^{-V(x)} \, dx$ a probability measure on  $\R^D$ such that $D^2 V \ge 0$.
For any probability measures $f_0, f_1 \in \PPP_2(\R^D)$, there holds
\beqn\label{ineq-HWI}
H_D(f_0 | \rho) \le  H_D(f_1 |\rho) +  \tilde W_2(f_0,f_1) \,\sqrt{
I_D (f_0 |\rho)},
\eeqn
where $H_D$ and $I_D$  stand for the non normalized relative
entropy and relative Fisher information defined in \eqref{def:EntropyRel} and
\eqref{def:FisherRel}
respectively, 
and  $ \tilde W_2$ stands for the non  normalized quadratic MKW distance in $\R^D$
based on the usual  Euclidean norm $|V| = (\sum_{i=1}^D |v_i|^2)^{1/2}$ for any 
$V = (v_1, ..., v_D) \in \R^D$. Inequality \eqref{ineq-HWI} has been proved in \cite{OttoVillani}, see also \cite{VillaniTOT,VillaniOTO&N,MR1842429,MR1846020,CorderoHWI}. 
We easily deduce  the {\it ``non relative"} inequality \eqref{eq:GalHWI}  from the 
{\it ``relative"} inequality \eqref{ineq-HWI}. In order to do so, we simply apply the  HWI inequality  \eqref{ineq-HWI} 
in $\R^D$, $D = dN$, with respect to the Gaussian
$\gamma_\lambda(v) := (2\pi\lambda)^{-D/2} e^{- |v|^2/2\lambda}$,  and we get 
$$
H_D(F^N | \gamma_\lambda) \le  H_D(G^N |\gamma_\lambda) +  \tilde
W_2(F^N,G^N) \,\sqrt{ I_D (F^N |\gamma_\lambda)}.
$$
We write  the relative entropy and the relative Fisher information in terms of the
non-relative ones, and we get 
$$
H_D(F^N | \gamma_\lambda) = H_D(F^N) - \int F^N \ln (\gamma_\lambda) =
 H_D(F^N) + \frac D2 \log( 2 \pi \lambda ) + \frac{M_2(F^N)}{2\lambda},
$$
\bean
I_D (F^N |\gamma_\lambda)  &=& \int F^N \left|  \nabla \ln F^N + \frac v
\lambda \right|^2  = I_D (f_0) + \frac 2\lambda \int v \cdot \nabla f_0 +
\frac{M_2(f_0)}{\lambda^2}\\
& = &  I_D (f_0) - \frac {2D}\lambda + \frac{M_2(f_0)}{\lambda^2}.
\eean
Inserting this in the relative HWI inequality, simplifying the terms involving $ \log( 2 \pi \lambda )$, letting $\lambda \to +
\infty$ and dividing the resulting limit by $N$, we obtain the claimed result. 
 
 \smallskip
\noindent 
{\sl Step 2. $E\subset\R^d$ is convex. } The proof is the same as in the case $E = \R^d$ using that the 
HWI inequality \eqref{ineq-HWI} holds in a convex set.   We have no precise reference for that last result but all the
necessary arguments can be find  in \cite{VillaniOTO&N}. More precisely, \cite[Chapter  20]{VillaniOTO&N} explains that
 the HWI inequality \eqref{ineq-HWI} holds  when the entropy is displacement convex, while 
  it is proved in \cite[Chapters 16 and  17]{VillaniOTO&N} that the entropy on a convex set $E$ is displacement convex, exactly as on $\R^d$. 
 
  \smallskip
\noindent 
{\sl Step 3. General case. }  
We choose two absolutely continuous probability measures $F^N$ and $G^N$ on $E^N$, and defined the corresponding  probability measures $F_1^N$ and $G^N_1$ on $E_1^N$ by
$$
F_1^N(v_1,\ldots,v_N) := F^N( T(v_1),\ldots,T(v_N)) = F^N \circ T^{\otimes N} (V),
$$
and the same formula for $G^N_1$. It can be checked that $\nabla_{v_j} F_1^N =    {}^t\nabla T(v_j) \nabla_{v_j} F^N \circ  T^{\otimes N} $, so that $|\nabla_{v_j} F_1^N| \leq   \| \nabla T \|_\infty \, |\nabla_{v_j} F^N \circ T^{\otimes N}|$.  Turning to Fisher information, it comes
$$
I(F^N_1) := \int_{E_1^N} \frac{| \nabla F_1^N |^2}{F^N_1} \,dV  \leq  \| \nabla T \|_\infty^2 
\int_{E_1^N}  \frac{| \nabla F^N \circ T^{\otimes N} |^2}{F^N \circ T^{\otimes N} } \,dV = \| \nabla T \|_\infty^2 \, I(F^N),
$$
where we have used the fact that $T$ preserves the volume.

For the MKW distance, remark that $|(T^{-1})^{\otimes N} (V) - (T^{-1})^{\otimes N} (V' )| \leq \| \nabla T^{-1} \|_\infty \, |V-V'|$. Therefore,
\bean
W_2(F^N_1,G^N_1)^2 & = & \inf_{\pi_1 \in \Pi(F^N_1,G^N_1)} \int |V- V'|^2 \,\pi_1(dV,dV') \\
& = & \inf_{\pi \in \Pi(F^N,G^N)} \int |(T^{-1})^{\otimes N} (V) - (T^{-1})^{\otimes N} (V' ) |^2 \,\pi(dV,dV') \\
& \le & \| \nabla T^{-1} \|_\infty^2 \, \inf_{\pi \in \Pi(F^N,G^N)} \int |V- V'| \,\pi(dV,dV') \\
& =& \| \nabla T^{-1} \|_\infty^2 \, W_2(F^N,G^N)^2.
\eean
For the entropy, the preservation of volume ensures the equality $H(F^N_1) = H(F^N)$, and a similar one for $G^N$. Finally, using the HWI inequality in $E_1$ 
proved in step 2 and the above properties, we get
\bean
H(F^N) &=& H(F^N_1) \leq H(G^N_1) +  \sqrt{I(F^N_1)} \, W_2(F^N_1,G^N_1) \\
&  \leq & H(G^N) + \| \nabla T \|_\infty \, \| \nabla T^{-1} \|_\infty \, \sqrt{I(F^N)} \, W_2(F^N,G^N),
\eean 
which is exactly the claimed result.
 \qed

\medskip
Let us finally prove now our main result Theorem~\ref{theo:LesDiffChaos}
which is a
consequence of the characterization of the Kac's chaos in Theorem~\ref{Th:EquivChaosAccurate}
together with Proposition~\ref{prop:HWI}.

 \smallskip
\noindent
{\sc Proof of Theorem~\ref{theo:LesDiffChaos}. }  We recall that the implication (iii) $\Rightarrow$ (iv) has been
yet proven in Theorem~\ref{theo:ChaosEntrop}. We split the proof into two steps. 

\smallskip\noindent
 {\sl Step 1. (i) $\Rightarrow$ (ii). } Fix a $j \in \N$, there exists a
subsequence of $(G^N)$, still denoted by $(G^N)$, and some compatible 
and symmetric probability measures $F_j \in \PPP(E^j)$,
such that $G^N_j \to F_j$ weakly in $\PPP(E^j)$. In particular
$F_1 = f$. As a consequence of Lemma~\ref{lem:Fisher1} and
Lemma~\ref{lem:FisherCondiv} point (i), we have
$$
I(f) \le I(F_j) \le \liminf  I(G^N_j) \le \liminf I(G^N)  = I(f).
$$
Using now the third point of Lemma~\ref{lem:FisherCondiv} we deduce $F_j =
f^{\otimes j}$. The uniqueness of the limit implies that the whole sequence $G^N$
is in fact $f$-Kac's chaotic.

\smallskip\noindent
 {\sl Step 2. (ii) $\Rightarrow$ (iii). } We write twice the normalized non relative HWI inequality of
Proposition~\ref{prop:HWI}, and get  
$$
|H(G^N) - H(f^{\otimes N})| \le 
  C_E \,  W_2 (G^N,f^{\otimes N}) \, \Bigl( \sqrt{I(G^N)} + \sqrt{I ( f^{\otimes N})} \Bigr). 
$$
Using the previous inequalities together with the inequality of the
Lemma~\ref{lem:ComparDistEN} 
$$
W_2(G^N,f^{\otimes N}) \le C_E \, 2^{\frac32} \, [M_k(G^N_1) +M_k(f)]^{1/k}\, W_1
(G^N,f^{\otimes N})^{1/2-1/k}  $$
we get  \eqref{eq:HGN-Hf<WN} since $M_k(f) \leq \sup M_k(G^N_1)$ and $I(f) \leq \sup I(G^N)$ .
\qed


\bigskip
\section{Probability measures on the ``Kac's spheres"}

\label{sec:Sphere}
\setcounter{equation}{0}
\setcounter{theo}{0}



We generalize the preceding two sections to the important case of probability measures with
support on the ``Kac's spheres"
$$
\KK\SS_N :=  \{V = (v_1, ..., v_N) \in \R^N,  \,\,v_1^2 + ... + v_N^2 = N \}. 
$$
We refer to \cite{Carrapatoso1} where similar results are obtained to the (even
more important) case of probability measures with support on the ``Boltzmann's spheres" 
$$
\BB\SS_N :=  \{V = (v_1, ..., v_N) \in (\R^3)^N,  \,\,|v_1|^2 + ... + |v_N|^2 = N, \,\, v_1 + ... + v_N = 0 \}. 
$$

\bigskip

\subsection{On uniform probability measures on the Kac's spheres as $N\to\infty$}
\label{sec:H&S-6}

\begin{defin} For any $N \in \N^*$ and $r > 0$, we denote by  $\sigma^{N,r}$ the
uniform probability measure of $\R^{N}$ carried by the sphere  $S^{N-1}_r$ 
defined by 
$$
S^{N-1}_r := \{V \in \R^N; \,\,\, |V|^2 = r^2 \}.
$$
We define $\sigma^N \in \PPP(E^N)$, $E=\R$, the sequence  $\sigma^N :=
\sigma^{N,\sqrt{N}}$ of probability measures uniform on the Kac's spheres
$$
\KK\SS_N := S^{N-1}_{\sqrt{N}}  := \{V \in \R^N; \,\,\, |V|^2 = N \}.
$$
\end{defin}

We begin with a classical and elementary lemma that we will use several times  in the sequel. 

\begin{lem}\label{lem:sigmaNell} (i) For any $1 \le \ell \le N-1$, there holds
$$
\sigma^N_\ell (V) = \Bigl(  1 - { |V|^2 \over N}\Bigr)_{\!+}^{{N-\ell-2 \over 2}} \, { |S_1^{N-\ell-1} | \over N^{\ell/2} \, |S_1^{N-1}|} , 
$$
where we recall that $ |S^{k-1}_1| = 2 \, \pi^{k/2} / \Gamma(k/2)$.  

(ii) For any fixed $\ell$, the sequence  $(\sigma^N_\ell)_{N \ge N_\ell}$
is bounded in $L^\infty$ (with $N_\ell = \ell+4$),   in  $H^s$ for any $s \ge 0$ (with $N_\ell = N(\ell,k)$ large enough)
and the exponential moment $M_{2,1/4} (\sigma^N_1)$ defined in~\eqref{def:ExpMoment} is bounded (uniformly in $N$).

(iii) For any function $\varphi \in C_b(\R^N)$, any  $r > 0$ and $1 \le
\ell \le N-1$, there holds
$$
\int_{S^{N-1}_r} \varphi(V,V') \, d\sigma^{N}_r(V,V') 
= \int_{B^\ell(r)} { |S^{N-\ell-1}_{\sqrt{r^2-V^2}} | \over | S^{N-1}_r |}
\left\{ 
\int_{S^{N -\ell-1}_{\sqrt{r^2-V^2}}}  \varphi(V,V') \,
d\sigma^{N-\ell}_{\sqrt{r^2-V^2} }(V') 
\right\} dV,
$$
where $V \in  \R^\ell$ and $V' \in \R^{N-\ell}$. This precisely means that
$$
\sigma^N (dV,dV') = \sigma^N_\ell (dV) \, \sigma^{N-\ell}_{\sqrt{N-|V|^2}} (dV').
$$
\end{lem}

\noindent
{\sc Proof of Lemma~\ref{lem:sigmaNell}. } (i) One possible definition of $\sigma^{N,r}$ is
$$
  \sigma^{N,r}  := {1 \over r^{N-1} \, |S_1^{N-1}|} \lim_{h \to 0} {1 \over h}
\, \Bigl( {\bf 1}_{B^N(r+h)}-  {\bf 1}_{B^N(r+h)} \Bigr),
  \qquad B^N(\rho) := \{V \in \R^N; \,\, |V| \le \rho \}, 
  $$
where the surface $r^{N-1} \, |S_1^{N-1}|$ of the Sphere $S^{N-1}_r$ stands
for the normalization constant such that $\sigma^{N,r}$ is a probability
measure. For any $\varphi \in C_b(E^\ell)$, $1 \le \ell \le N-1$, we compute
\bean
\Bigl\langle {\bf 1}_{B(\rho)} , \varphi \otimes {\bf 1}^{N-\ell} \Bigr\rangle
&=& \int_{\R^\ell} {\bf 1}_{|V|^2 \le \rho^2} \, \varphi(V)  
\left\{\int_{\R^{N-\ell}}  {\bf 1}_{x_{\ell+1}^2 + ... + x_N^2 \le \rho^2-|V|^2} 
\, dx_{\ell+1} \, ... \, dx_N \right\}dV
\\
&=& \int_{\R^\ell} \varphi(V)  \, \omega^{N-\ell} \, (\rho^2 - |V|^2)^{N-\ell
\over 2}_+ \, dV,
\eean 
where $\omega^k = |B^k(1)|$ is the volume of the unit ball of  $\R^k$. We deduce
$$
  \sigma^N_\ell(r) = {1 \over Z_{N,r}} \, {d \over dr} \left[ \omega^{N-\ell} \,
(r^2 - |V|^2)^{N-\ell \over 2}_+ \right] = 
  {\omega^{N-\ell} \, (N-\ell) \over r^{N-1} \, |S^{N-1}|}\, r \,  (r^2 -
|V|^2)^{N-\ell -2 \over 2}_+.
  $$
  We conclude using the relation $ |S_1^{k-1}| = k \, \omega^k$.

\medskip
\noindent
(ii) The estimates on  $\sigma^N_\ell$ are deduced from its explicit
expression after some tedious but easy  calculations.  
We only prove the last one which will be a key argument in the proof of
the accurate rate of chaoticity in Theorem~\ref{theo:EmpiricalMeasure}. 
For any $k \ge 1$ and introducing  $n := (N-4)/2$, we easily estimate
\bean
\int_{\R^2} |v_1|^{2k} \, \sigma^N_2 (dv)
&=& {1 \over 2\pi} \, {N-2 \over N}  \int_{\R^2}   |v_1|^K  \, \Bigl(  1 - { |v|^2 \over N} \Bigr)_{\!+}^{{N-4 \over 2}}  \, dv
\\
&\le&    \int_{0}^{\sqrt{N}}   r^{K+1}  \, \Bigl(  1 - { r^2 \over N} \Bigr)^{{N-4 \over 2}}  \, dr
\\
&=&   N^{k+1} \int_{0}^{1}   s^{k}  \, \Bigl(  1 -  s  \Bigr)^{n}  \, ds.
\eean
Thanks to $k+1$  integrations by parts, we deduce
\bean
\int_{\R^2} |v_1|^{2k} \, \sigma^N_2 (dv)
&\le& N^{k+1} \int_{0}^1 (1-z)^k  \, z^{n} \, dz
\\
&=& N^{k+1}\, {k \over n+1} \int_{0}^1 (1-z)^{k-1}  \, z^{n+1} \, dv
\\
&=& N^{k+1}\, {k \over n+1} \, ... \, {2 \over n+k-1} \, {1 \over n+k} \, {1 \over n+k+1} ,
\eean
and then
\bean
\int_{\R^2} e^{|v|^2/4} \, \sigma^N_1 (v) \, dv 
&\le& \sum_{k=0}^\infty {1 \over k! \, 4^k}   \int_{\R^2} \, |v_1|^{2k} \, \sigma^N_2 (dv)
\\
&=& \sum_{k=0}^\infty {1 \over  4^k}   \, {(2n+4)^{k+1} \over (n+1) \, ... ( n+k+1)}  
\\
&  \le & 2 \, \sum_{k=0}^\infty {1 \over  2^k}   \, {(n+2) \over (n+1)}   \le 6.
\eean

\medskip
\noindent
(iii) We come back to the proof of (i). We set $m=\ell$ and $n=N-\ell$ and we
write
\bean
\langle \sigma^{N,r}, \varphi \rangle
&=& {1 \over Z_{N,r}} 
\lim_{h\to0} \, {1 \over h}  \left[\int_{B^N(r+h)} \varphi - \int_{B^N(r)}
\varphi \, \right]
\\
&=& {1 \over Z_{N,r}} 
\lim_{h\to0} \, {1 \over h}  \left[\int_{|v|\le r+h} \int_{|v'| \le \sqrt{
(r+h)^2 - |v|^2}}  \varphi - \int_{|v|\le r} \int_{|v'| \le \sqrt{ (r+h)^2 -
|v|^2}}   \varphi \, \right]
\\
&&+ {1 \over Z_{N,r}} 
\lim_{h\to0} \, {1 \over h}  \left[\int_{|v|\le r} \int_{|v'| \le \sqrt{ (r+h)^2
- |v|^2}}  \varphi -  \int_{|v|\le r} \int_{|v'| \le \sqrt{ r^2 - |v|^2}}  
\varphi \, \right]
\\
&=& {1 \over Z_{N,r}} 
\lim_{h\to0} \, {1 \over h}   \int_{r \le |v|\le r+h} \int_{|v'| \le \sqrt{
(r+h)^2 - |v|^2}}  \varphi  
\\
&&+ {1 \over Z_{N,r}} \int_{B^m_r} \lim_{h\to0} \, {1 \over h}  \left[ 
\int_{B^n(\sqrt{ (r+h)^2 - |v|^2}) }  \varphi - \int_{B^n(\sqrt{ r^2 - |v|^2})
}  \varphi \, \right].
\eean
We invert the integral and the limit on the last line using dominated
convergence, since the integral on $v'$ are bounded by $\|\varphi\|_\infty /
{\sqrt{r^2 - |v|^2}}$. The first term is bounded (for any $0 < h \le r$) by
$$
{1 \over Z_{N,r}} 
\lim_{h\to0} \, {1 \over h}   \int_{r \le |v|\le r+h} \int_{|v'| \le \sqrt{ 3 \,
r \, h }}  | \varphi |\le C_{N,r} \, \| \varphi \|_{L^\infty}  \,\lim_{h\to0} 
\sqrt{h} = 0,
$$
and the second term converges to 
$$
\int_{B^m(r)} { Z_{n,\sqrt{r^2-|v|^2}}  \over Z_{m+n,r} } \left\{ 
\int_{S^{n-1}_{\sqrt{r^2-|v|^2}}}  \varphi(v,v') \, d\sigma^n_{\sqrt{r^2-|v|^2}
}(v') 
\right\}dv,
$$
which is exactly the claimed identity.  
\qed

\smallskip
Let us recall the following classical result. 

\begin{theo} \label{theo:Diaconis-Freedman}
The sequence  $\sigma^N$ is $\gamma$-chaotic, where $\gamma$ still stands for
the gaussian distribution $\gamma (dx) = {(2 \pi)^{-1/2}} \, e^{-x^2/2} \,dx$ on  $\R$, 
and more precisely 
\beqn\label{estim:Poincar√©1}
  \| \sigma^N_\ell - \gamma^{\otimes \ell} \|_{L^1}\leq  \quad  2 \, { \ell+3
\over N-\ell-3} \quad \hbox{pour tout} \quad 1 \le \ell \le N-4.
\eeqn
\end{theo}

The fact that $\sigma^N$ is $\gamma$-chaotic is sometime called ``Poincar\'e's Lemma". 
In fact, it should go back to Mehler \cite{Mehler} in 1866. Anyway, we refer to 
\cite{DiaconisFreedman1987,CCLLV} for a bibliographic discussion about this important result,
and to  \cite{DiaconisFreedman1987} for a proof of estimate \eqref{estim:Poincar√©1}. 
We give now a different quantitative version of the ``Poincar\'e's Lemma".
 
\begin{theo}\label{theo:ChaosQuantifPoincare} There exists a
numerical constant $C \in (0,\infty)$ such that 
\beqn  \label{estim:Poincar√©2}
\Omega_N(\sigma^N;\gamma) := W_1(\sigma^N,\gamma^{\otimes N}) \le \frac C {\sqrt N}.
\eeqn
\end{theo}

\begin{rem}\label{rem:ChaosQuantifPoincare} It is worth observing that it is not
clear that one can deduce  \eqref{estim:Poincar√©2} from \eqref{estim:Poincar√©1}
or that the reverse implication holds. In particular,
using \eqref{estim:Poincar√©1} and Theorem~\ref{Th:EquivChaosAccurate}  we obtain an
estimate on $W_1(\sigma^N,\gamma^{\otimes N})$ which is weaker than
\eqref{estim:Poincar√©2}. 
\end{rem}

\noindent
{\sc Proof of Theorem~\ref{theo:ChaosQuantifPoincare}. }
There is a simple transport map from $\gamma^{\otimes N}$ onto $\sigma^N$ which
is given by the radial projection $P : V \mapsto \frac{V}{|V|_2}$ with
the notation $|V|_k = ( N^{-1} \sum_i |v_i|^k)^{1/k}$ for any $k>0$ for the
normalized distance of order $k$. The fact it is an admissible map comes from
the invariance by rotation of $\gamma^{\otimes N}$ and
$\sigma^N$. Is it optimal?  It is not obvious because $P(V)$ is not
necessary the point of $\KK\SS_N$ wich is the closest to $V \in \R^N$, for the
$|\cdot|_1$ distance (for which it costs less to displace in the direction of
the axis). However, it may still be optimal for rotationnal symmetry reasons,
but it is less obvious. Nevertheless, it will be sufficient for our estimate.
Since,
$$
|P(V) - V|_1 = \left| {\textstyle \frac{1}{|V|_2}} - 1 \right| |V|_1
$$
we get as all our distances are normalized
\bean
W_1(\gamma^{\otimes N},\sigma^N) & \leq & \int_{\R^N} |P(V) - V|_1
\gamma^{\otimes N}(dV)
\\
& = & \int_{\R^N} \left| {\textstyle \frac1{|V|_2}} - 1
\right| |V|_1 \gamma^{\otimes
N}(dV)
\\
& = & \left(\int_0^{+\infty} \left|
{\textstyle \frac{\sqrt N}R} - 1 \right| R^N e^{-R^2/2} \,dR \right) 
 \frac{|S^{N-1}|}{(2\pi)^{N/2}} \left(
   \int_{S^{N-1}_1} |V|_1 \,d \sigma^{N,1}\right).
\eean
 Using that $|V|_1 \leq |V|_2$ because of the normalization, we may
bound the last integral by 
\bean
\int_{S^{N-1}_1} |V|_1 \,d \sigma^{N} \leq 
\int_{S^{N-1}_1 } |V|_2 \,d \sigma^{N} = N^{-1/2}.
\eean
Remark that this integral is also equal to $\frac1{\sqrt N} M_1(\sigma^N)$
which can be explicited thanks to the formula for $\sigma_1^N$ of
Lemma~\ref{lem:sigmaNell}.
Using this in the previous inequality
and performing the change of variable $R = \sqrt N R'$, we get
\bean
W_1(\gamma^{\otimes N},\sigma^N) & \leq & \frac{|S^{N-1}|}{\sqrt N (2\pi)^{N/2}}
\int_0^{+\infty} | \sqrt N - R| R^{N-1} e^{-R^2/2} \,dR \\
& \leq & \frac{|S^{N-1}|N^{\frac N2}}{(2\pi)^{N/2}}
\int_0^{+\infty} | 1 - R'| (R')^{N-1} e^{-N R'^2/2} \,dR'.
\eean
We can simplify the prefactor, using the formula for $|S^{N-1}|$ and Stirling's
formula 
\bean
\frac{|S^{N-1}|\,N^{\frac N2}}{(2\pi)^{N/2}} & = &
\frac{N^{\frac N2}}{\Gamma(\frac N2) 2^{N/2-1} } \\
& = & \frac{\sqrt N e^{N/2}}{\sqrt \pi}[1 + O(1/N)].
\eean
Turning back to the transportation cost, we get
$$
W_1(\gamma^{\otimes N},\sigma^N)  \leq  \frac{e \,\sqrt N }{\sqrt \pi}[1 +
O(1/N)]  \int_0^{+\infty} \left(R e^{(1- R^2)/2}\right)^{N-1} e^{-
R^2/2} | 1 - R| \,dR.
$$
After studying the function $g(r) = r e^{(1- r^2)/2}$, we remark that it is
strictly increasing form $0$ to $1$, then strictly decreasing from $1$ to
$+\infty$, that its maximum in $1$ is $1$, and that $g(1 + \eps) = 1 - \eps^2+
O(\eps^3)$. We shall also use the less sharp but exact bound
$$
g(1+ \eps) \leq 1 - \frac{\eps^2}4, \quad \text{ for } \eps \in [ -\frac12 ,
\sqrt2 -1 ].
$$
We can  now cut the previous integral in three parts $\ds  \int_0^{1/2} +
\int_{1/2}^{\sqrt 2} + \int_{\sqrt 2}^{+\infty}$. We bound the first part by
$$
\int_0^{1/2} \ldots \leq \frac12 g\left({\textstyle \frac12}\right)^{N-1},
$$
and the third part by
$$
\int_{\sqrt 2}^{+\infty}  \ldots \leq g(\sqrt 2)^{N-1} \int_{\sqrt 2}^{+\infty} 
e^{- r^2/2} r \,dr = g(\sqrt 2)^{N-1} e^{-1}.
$$
For the last part, we perform the change of variable $r = 1 + u/\sqrt N$. It comes
\bean
\int_{1/2}^{\sqrt 2} \ldots  & = & \frac1N \int_{-\sqrt N / 2}^{(\sqrt 2-1)\sqrt
N} g\left( 1 + {\textstyle \frac u{\sqrt N}} \right)^{N-1} |u| e^{-u/\sqrt N -
u^2/2 N} \,du \\
& \le & \frac1N \int_{-\sqrt N / 2}^{(\sqrt 2-1)\sqrt
N} \left( 1 -  {\textstyle \frac {u^2}{4 N}} \right)^{N-1} |u| \,du \\
& \le & \frac1N (1 + O(N^{-1})) \int_{-\infty}^{+\infty} 
 e^{-\frac {u^2}4 } |u| \,du \\
& \le & \frac4N (1 + O(N^{-1}))
\eean
Putting all together, we finally get
$$
W_1(\gamma^{\otimes N},\sigma^N)  \leq \frac C{\sqrt N} (1+ O(N^{-1})) + C
\sqrt N \lambda^N,
$$
with $\lambda = \max(g(\sqrt 2),g(1/2)) <0.86$. This implies the claimed 
inequality. \qed

\medskip
\noindent
{\sc Proof of \eqref{eq:ChaosEstimSigmaN} in Theorem~\ref{theo:EmpiricalMeasure}. }    
The proof of the last estimate in \eqref{eq:ChaosEstimSigmaN} follows from \eqref{estim:Poincar√©2}
and Lemma~\ref{lem:sigmaNell}-(ii)
together with \eqref{eq:EquivChaosExpo}. \qed


\subsection{Conditioned tensor products on the Kac's spheres}
 \label{subsec:condtensprod}

We begin with a sharp version of the local central limit theorem (local CLT) or
Berry-Esseen
type theorem which will be the cornerstone argument in this section. 

\begin{theo}\label{theo:localTCL} Consider $g \in \PPP_3(\R^D) \cap L^p(\R^D)$,
$p \in (1,\infty]$, such that 
\beqn\label{hypo:localTCL}
  \int_{\R^D} x \, g(x) \, dx = 0, \quad
\int_{\R^D} x \otimes x \, g(x) \, dx = Id, \quad \int_{\R^D} |x|^3 \, g(x) \,
dx =: M_3. 
\eeqn
We define the iterated and renormalized convolution by   
\beqn\label{estim:localTCLinfty}
g_N(x) := \sqrt{N} \, g^{(*N)} (\sqrt{N} \, x) . 
\eeqn 
There exists an integer $N(p)$ and a constant $C_{BE}= C(p,k,M_3(g),\|g
\|_{L^p})$ such that 
\beqn\label{estim:localTCL2}
\forall \, N \ge N(p) \qquad
\|  g_N - \gamma \|_{L^\infty} \le \frac{C_{BE}}{\sqrt N }. 
\eeqn

\end{theo} 

\begin{rem}
Theorem~\ref{theo:localTCL} is a sharper but less general version of 
\cite[Proposition 26]{CCLLV}. 
The proof follows the proof of \cite[Proposition 26]{CCLLV} and uses an argument
from
\cite[Proposition 26]{CCLLV}, see also \cite{LionsToscani1995}. 
The first  local CLT have been established in  the pioneer works by 
A.~C.~Berry \cite{Berry1941} and C.-G. Esseen \cite{Esseen1942} who proved the
convergence in 
$\OO(1/\sqrt{N})$ uniformly on the distribution fonction in dimension
$D=1$, see for instance \cite[Theorem 5.1, Chapter XVI]{FellerBookII}.
Since that time, many variants of the local CLT have been established
corresponding to different regularity assumption 
made on the probability measure $g$,  we refer the interested reader to the recent
works  \cite{Rio2009}, \cite{Bobkov-BEth},  \cite{Barthe&co2004}
and the references therein. 
\end{rem}

The proof of Theorem~\ref{theo:localTCL} use the following technical lemma which
proof is postponed after the proof of the Theorem. 
  
\begin{lem}\label{lem:borneFourierg} (i)  Consider $g \in
\PPP_3(\R^D)$ satisfying  \eqref{hypo:localTCL}. There exists $\delta \in (0,1)$
such that
$$
\forall \, \xi \in B(0,\delta) \qquad |\hat g (\xi)| \le e^{- |\xi|^2/4}.
$$
(ii) Consider $g \in \PPP(\R^D) \cap L^p(\R^D)$, $p \in (1,\infty]$. For any
$\delta > 0$ there exists  $\kappa = \kappa(M_3(g),\| g \|_{L^p},\delta) \in
(0,1)$ such that 
\beqn\label{estim:hatgGRANDEfrequance}
\sup_{|\xi| \ge \delta} |\hat g (\xi) | \le \kappa(\delta).
\eeqn
\end{lem}

\noindent
{\sc Proof of Theorem~\ref{theo:localTCL}. } We follow closely the proof of
\cite[Theorem 27]{CCLLV} which is more general but less precise, and we use a
trick
that we found in the proof of \cite[Theorem 1]{GoudonJT2002}.  We observe that 
$$
\hat g_N(\xi) =   (\hat g (\xi/\sqrt{N}))^N,
\qquad
\hat \gamma (\xi) =  (\hat \gamma (\xi/\sqrt{N}))^N.
$$
Because $g \in L^1 \cap L^p$, the Hausdorff-Young inequality implies 
 $\hat g \in L^{p'} \cap L^\infty$ with $p' \in [1,\infty)$, and then 
 $\hat g_N (\xi) = (\hat g(\xi/\sqrt{N}) )^N \in L^1$ for any $N \ge p'$. As a
consequence we may write 
$$
|g_N(x) - \gamma(x)| = (2\pi)^D \, \left| \int_{\R^D} (\hat g_N (\xi) -\hat
\gamma (\xi)) \, e^{i \, \xi \cdot x}\, d\xi \right|
\le (2\pi)^D \, \int_{\R^D}  |\hat g_N - \hat \gamma| \, d\xi.
$$
We split the above integral between low and high frequencies  
\bean
\| g_N - \gamma \|_{L^\infty} 
&\le& \int_{|\xi| \ge \sqrt{N} \, \delta} \left|\hat g_N   \right|\, d\xi   + 
\int_{|\xi| \ge \sqrt{N} \, \delta} \left| \hat \gamma    \right| \, d\xi  
\\
&&+ \int_{|\xi| < \sqrt{N} \, \delta}   \left|\hat g_N  - \hat \gamma \right|
 \, d\xi  \qquad (=: T_1 + T_2 + T_3 ).  
\eean
For the first term, we have 
\bean
T_1  
&\le&  \int_{|\xi| \ge \sqrt{N} \, \delta} \left|\hat g \left( {\xi \over
\sqrt{N}} \right)  \right|^N \, d\xi 
= N^{d/2} \int_{|\eta| \ge   \delta} \left| \hat g \left( \eta
\right)\right|^{N}    \, d\eta 
\\
&\le& \left(\sup_{|\eta| \ge \delta}|\hat g(\eta)|\right)^{N-p'} \, N^{d/2} 
\int_{\eta >  \delta} \left| \hat g \left( \eta  \right)\right|^{p'}   \, d\eta 
\\
&\le& \kappa(\delta)^{N-p'} \, N^{d/2}  \, C_p \, \|g \|_{L^p}^p
\eean
with $\delta \in (0,1)$ given by point (i) of
Lemma~\ref{lem:borneFourierg}, $\kappa(\delta)$ given by  point (ii) of
Lemma~~\ref{lem:borneFourierg}  and $N \ge p'$. The second term may be estimated
in the same way, and we clearly obtain that there exists a constant  $C_1 =
C_1(D,p,\|g \|_{L^p})$ such that  
\beqn\label{bdd:T1+T_1}
T_1 + T_2 \le {C_1 \over \sqrt{N}}.
\eeqn
Concerning the third term, we write 
\bean
T_3
&=&  
\int_{|\xi| \le \sqrt{N} \, \delta} {\left|\hat g_N (\xi) - \hat   \gamma_N 
(\xi)  \right|  \over |\xi|^3} \, |\xi|^3   \, d\xi  , 
\eean
with 
\bean
{\left|\hat g_N (\xi) - \hat   \gamma_N  (\xi)  \right|  \over |\xi|^3}  
&=& {1 \over N^{3/2}}\, {\left|\hat g (\xi/\sqrt{N}) ^N - \hat  
\gamma(\xi/\sqrt{N}) ^N   \right| \over |\xi/\sqrt{N}|^3}  
\\
&=&   
 {1 \over N^{3/2}}\, {\left|\hat g (\xi/\sqrt{N}) - \hat   \gamma(\xi/\sqrt{N})
 \right|  \over |\xi/\sqrt{N}|^3}  
\\
&& \times \,  \left|\sum_{k=0}^{N-1} \hat g (\xi/\sqrt{N})^k \, \hat  
\gamma(\xi/\sqrt{N}) ^{N-k-1}   \right|.
\eean
Estimate (i) of Lemma~\ref{lem:borneFourierg}  implies  
\bean
 &&\left|\sum_{k=0}^{N-1} \hat g (\xi/\sqrt{N})^k \, \hat \gamma(\xi/\sqrt{N})
^{N-k-1}   \right| 
 \\
 &&\le 
\sum_{k=0}^{N-1} e^{ - {|\xi|^2 \over 4 \, N } \, k} \, e^{ - {|\xi|^2 \over 2
\, N } \, (N-k-1)}  \le N \, e^{ - {|\xi|^2 \over 4 } \,  {  N-1 \over N} }  \le
N   \, e^{ - {|\xi|^2 \over 8 } }. 
 \eean
We deduce 
\bean
T_3
&=&  {1 \over N^{3/2}} \, \left( \sup_{\eta} {{ |\hat g (\eta) - \hat  
\gamma(\eta)  |  \over |\eta|^3}}  \right) \, 
\int_{\R^D} N   \, e^{ - {|\xi|^2 \over 8 } } \,    |\xi|^3 \, d\xi 
\\
&\le&  {1 \over N^{1/2}} \, (M_3(g) + M_3(\gamma))  \, C_{k,d}.
\eean
We conclude by gathering the estimates on each term.  \qed 

\medskip 
\noindent{\sc Proof of Lemma~\ref{lem:borneFourierg}. } Thanks to a Taylor
expansion, we have 
\bean
\hat g (\xi) &=& 1 - {\xi^2 \over 2} + \OO (M_3(g) \, |\xi|^3)
\\
\hat \omega (\xi) &=& 1 - {\xi^2 \over 4} + \OO ( |\xi|^3),  \qquad \omega(x) :=
{1 \over \sqrt{\pi}} \, e^{-x^2 }, 
\eean
from which we deduce that there exists $\delta = \delta (M_3(g))  \in
(0,1)$ small enough such that 
$$
\forall \, \xi \in B_\delta \quad
|\hat g(\xi)|\le 1 - {3 \over 8} \, \xi^2 \le \hat\omega(\xi) , \quad
\hat\omega(\xi):=  \, e^{-\xi^2/4}.
$$
That is nothing but (i). On the other hand, (ii) is a consequence of 
\cite[Proposition 26, (iii)]{CCLLV}. \qed 
 
  \bigskip

For a given ``smooth enough" probability measure  $f \in \PPP(E)$, $E = \R$,  we define 
$$
Z_N(r) := \int_{S^{N-1}(r)} f^{\otimes N}  \, d\sigma^{N,r}, \qquad 
Z'_N(r) := \int_{S^{N-1}(r)}  {f^{\otimes N} \over \gamma^{\otimes N}} \,
d\sigma^{N,r} = { Z_N(r) \over \gamma^{\otimes N} (r)}.
$$

We give a sharp estimate on the asymptotic behavior of $Z'_N$ as $N\to\infty$. 

\begin{theo}\label{theo:asymptotZN} Consider $f \in \PPP_6(\R) \cap L^p (\R)$,
$p \in (1,\infty]$,  satisfying 
\beqn\label{hypo:localTCL2}
\int_\R f \, v \, dv = 0,
\eeqn 
and define
\beqn\label{hypo:localTCL3}
E := \int_\R f \, |v|^2 \, dv, \qquad
\Sigma := \left( \int_\R (v^2-E)^2 \, f(v) \, dv \right)^{1/2}.
\eeqn
Then $Z_N(r), Z_N'(r)$ are well defined for all $r>0$ and there holds with the
above notations
\beqn\label{eq:asympZprim}
Z'_N(r)  \, \alpha_N(r^2) =  { \sqrt{2} \over \Sigma}\, \alpha_N(N) \ \left(
\exp \left\{ -  \left( {r^2-N \, E \over  \sqrt{N} \, \Sigma} \right)^2 \right\}
+ \frac{R_N(r)}{\sqrt N} \right)
\eeqn
where
$$
 \alpha_N(s) = s^{{N \over 2} -1}\, e^{-{s \over 2}}
 \quad\text{and}\quad  \| R_N \|_\infty \leq   C(p,\|f\|_p,M_6(f))
$$
\smallskip
As a particular case, there holds
\beqn\label{eq:asympZprimN}
Z'_N := Z'_N(\sqrt{E \, N}) =  {\sqrt{2} \over \Sigma} \, \left( 1  +  \OO
\Bigl( N^{-\frac12}  \Bigr) \right).
\eeqn
\end{theo}

\noindent
{\sc Proof of Theorem~\ref{theo:asymptotZN}. } We follow the proof
of \cite[Theorem 14]{CCLLV} but using the sharper estimate
proved in Theorem~\ref{theo:localTCL} (instead of \cite[Theorem 27]{CCLLV}).

Before going on, let us remark that it is not obvious that $Z_N(f;r)$ is
well defined for all $r>0$ under our assumption on $f$ which is not necessarily
continuous, since we are restricting $f^{\otimes N}$ to surfaces of $\R^N$. 
But, in fact the product structure of $f^{\otimes N}$ makes it possible. To see
this, take $f$ and $g$ two measurable functions equal almost everywhere, and
call $\NN$ the negligible set on which they differ. Then the tensor products
$f^{\otimes N}$ and $g^{\otimes
N}$ differs only on the negligible set $\bar \NN = \cup_{i} \; \R^{\otimes(i-1)}
\times \NN \times \R^{\otimes (N-i)}$. It is not difficult to see that
because of the particular structure of $\bar \NN$, the intersection of $\bar \NN
\cap S^{N-1}_r$ is also $\sigma^N_r$-negligible for all $r>0$. Therefore
$f^{\otimes N}$ and $g^{\otimes N}$ are equal $\sigma^N_r$-almost everywhere on
$S^{N-1}_r$, and there is no ambiguity in the definition of $Z_N(f,r)$ for all
$r>0$.

\medskip
We now define the law $g$ of $v^2$ under $f$
\beqn\label{def:hFROMf}
h(u) := {1 \over 2 \, \sqrt{u}} \, (f(\sqrt{u}) + f(-\sqrt{u})) \, {\bf 1}_{u >
0}, 
\eeqn
remarking that  $h \in \PPP_3(\R) \cap L^q(\R)$ with $q > 1$ as it has been
shown in the proof of \cite[Theorem 14]{CCLLV}. Consider $(\VV_j)$ a sequence
of random variables which is i.i.d. according to $f$. On the one hand, the law
 $s_N(du)$ of the random variable 
$$
S_N := \sum_{j=1}^N  |\VV_j|^2 
$$
can be computed by writing 
\bean 
\Ee (\varphi (S_N))  
&=& \int_0^\infty \varphi(r^2) \, |S^{N-1}_1|\, r^{N-1} \Bigl(\int_{S^{N-1}_r}
f^{\otimes N}(V) \, \sigma^{N,r}(dV) \Bigr)  \, dr 
\\
&=& \int_0^\infty \varphi(u) \, |S^{N-1}_1| \, u^{{N-1\over 2}} \Bigl( \int_{
S^{N-1}_{\sqrt{u}} }  f^{\otimes N}(V) \, \sigma^{N,\sqrt{u}}(dV)  \Bigr)  \,
{du \over 2 \sqrt{u}},
\eean
which implies  
$$
s_N(du) =  {1 \over 2}  \, |S^{N-1}_1|\, u^{{N\over 2}- 1} \, Z_N(\sqrt{u}). 
$$
On the other hand, we have $s_N = h^{(*N)}$. Gathering these two identities, we
get 
\bear \nonumber
h^{(*N)} (r^2) 
&=& {1 \over 2}  \, |S^{N-1}_1|\, r^{N-2} \, Z_N(r)  = 
 {\pi^{N/2} \over \Gamma(N/2)}\, r^{N-2} \, Z'_N(r) \, {e^{-r^2/2} \over
(2\pi)^{N/2}} 
\\ \label{eq:asymptZN1}
&=& {\alpha_N(r^2)  \over \Gamma(N/2)}\, { Z'_N(r) \,  \over 2^{N/2} } . 
\eear
Let us define  $g(u) := \Sigma  \, h(E + \Sigma \, u)$, so that $g  \in
\PPP_3(\R) \cap L^q(\R)$ and
$$
\int_{\R } g(y) \, y \, dy = 0, \qquad\int_{\R } g(y) \, |y|^2 \, dy = 1.
$$
Applying Theorem~\ref{theo:localTCL} to $g$ and using the identity $g^{(*N)}(u)
= \Sigma  \, h^{(*N)}(N \, E + \Sigma \, u)$, we obtain 
\beqn\label{eq:asymptZN2}
\sup_{r \ge 0} \left| h^{(*N)}(r^2) - {1 \over \sqrt{N} \, \Sigma} \gamma \left(
{r^2-N \, E \over  \sqrt{N} \, \Sigma} \right)\right| \le {C_{BE} \over 
N \, \Sigma},
\eeqn
where  $C_{BE}$ is the constant given in Theorem~\ref{theo:localTCL} and
associated to $g$. Gathering the  Stirling formula 
\beqn\label{eq:asymptZN3}
\Gamma (N/2)  = { \sqrt{\pi N} \, \alpha_N(N) \, 2^{-{N\over2} + 1}} \, \left( 1
+   \OO(N^{-1/2})  \right),
\eeqn 
with \eqref{eq:asymptZN1}, \eqref{eq:asymptZN2}, we obtain
$$
\forall \, r > 0 \quad
\left| {\alpha_N(r^2)   Z'_N(r) \,  \over  \sqrt{\pi N} \, \alpha_N(N) \, 2 \, (
1 + \OO(N^{-1/2})) }
 - {1 \over \sqrt{N} \, \Sigma \, \sqrt{2\pi} } \exp \left( \left( {r^2-N \, E
\over  \sqrt{N} \, \Sigma} \right)^2/2 \right) \right| \le {C_{BE} \over 
N \, \Sigma}. 
 $$  
Estimate  \eqref{eq:asympZprim} readily follows.  \qed

\medskip
For a given $f \in \PPP_6(\R)  \cap L^p(\R)$, $p > 1$, we define the
corresponding sequence of  ``conditioned product measures" (according to the
Kac's spheres 
$ \KK\SS_N$), we write $F^N := [f^{\otimes N}]_{\KK\SS_N}$, by
\beqn\label{def:CondProd}
F^N := {1 \over Z_N(f;\sqrt{N})} \, f^{\otimes N} \, \sigma^N.
\eeqn

We show that  $(F^N)$ is  well defined for $N$ large enough and
is $f$-chaotic. 

\begin{theo}\label{theo:[fotimesN]fchaos} Consider $f \in \PPP_6(\R) \cap
L^p(\R)$, $p > 1$, satisfying 
\beqn \label{hyp:[fotimesN]fchaos} 
\int_\R f \, v \,dv =0 \qquad \text{and} \quad \int_\R f \,  v^2 \,dv =1.
\eeqn
The sequence $(F^N)$ of
corresponding conditioned product measure is $f$-chaotic, more precisely 
$$
\Omega_\ell(F^N,f) := W_1(F^N_\ell, f^{\otimes \ell} ) \le {1 \over 2}\, \|
F^N_\ell - f^{\otimes
\ell} \|_1 \le  {C \, \ell^2 \over \sqrt{N}},
$$
for some constant $C = C(f) \in (0,\infty)$. 
\end{theo} 

\begin{rem} The $f$-Kac's chaoticity property of the sequence $F^N = [f^{\otimes
N}]_{\KK\SS_N}$ 
is stated and proved for smooth densities $f$ in the seminal article by  M. Kac
\cite{Kac1956}. 
Next, the same chaoticity property is proved with large generality (on $f$) in 
\cite{CCLLV}. 
Theorem~\ref{theo:[fotimesN]fchaos} is a {\it ``quantified"} version of 
\cite[Theorems 4 \& 9]{CCLLV} 
and  \cite[paragraph 5]{Kac1956}.
\end{rem}

\noindent
{\sc Proof of Theorem~\ref{theo:[fotimesN]fchaos}. } 
As in Theorem~\ref{theo:asymptotZN}, it is not obvious that $F^N$ is
well defined under our assumption on $f$ which is not necessarily continuous,
since we are restricting $f^{\otimes N}$ to a surface of $\R^N$.  But the
argument given at the beginning of the proof of Theorem~\ref{theo:asymptotZN}
 shows in fact that the restriction of $f^{\otimes N}$ to $\KK\SS_N$ is
unambiguously defined. Since, Theorem~\ref{theo:asymptotZN}
implies that $Z_N(f,\sqrt N)$ is finite and non zero for $N$ large enough, 
we deduce that $F^N$ is well defined for $N$ large enough.

\smallskip
Let us fix $\ell \ge 1$ and
$N \ge \ell + 1$.  
Denoting $V = (V_\ell, V_{\ell,N})$, with $V_\ell = (v_j)_{1 \le j \le
\ell}$, $V_{\ell,N} = (v_j)_{\ell+1 \le j \le N}$,  
we write thanks to the equality (iii) of Lemma~\ref{lem:sigmaNell}
$$
F^N (dV) = \left( { f \over \gamma}\right) ^{\otimes \ell} \!\!\! (V_\ell) \, {1
\over Z'_N(\sqrt{N}) } \,  \left( { f \over \gamma} \right) ^{\otimes N-\ell} 
\!\!\!  (V_{\ell,N}) \,   
\sigma^{N-\ell,\sqrt{N-|V_\ell |^2}} (dV_{\ell,N})  \, \, \sigma_\ell^N (V_\ell)
\, dV_\ell,
 $$
so that, coming back to the notation $V = V_\ell = (v_j)_{1 \le j \le \ell} \in
\R^\ell$, we have 
$$
  F^N_\ell (V)  = \left( \prod_{j=1}^\ell { f(v_j) \over \gamma(v_j) }\right)
{ Z'_{N-\ell}(\sqrt{N-|V|^2}) \over Z'_N(\sqrt{N}) } 
 \,  \sigma^N_\ell (V) = \left(  \prod_{j=1}^\ell  f(v_j) \right)
\theta_{N,\ell}(V),
 $$
if we define the quantity $\theta_{N,\ell}$ by
\beqn \label{eq:deftheta}
\theta_{N,\ell}(V)   :=  (2\pi)^{\frac \ell 2} \, e^{\frac{|V|^2}2} {
Z'_{N-\ell}(\sqrt{N-|V|^2}) \over Z'_N(\sqrt{N}) }  \,  \sigma^N_\ell (V).
\eeqn
The key point is now to prove that $\theta_{N,\ell}$ goes to $1$.
Recalling the Stirling formula $\Gamma(k) = \sqrt{\frac{2\pi}k}\left(
\frac k e \right)^k (1+ \OO(k^{-1}))$, we write $\sigma^N_\ell$ as 
\bean
\sigma^N_\ell (V) 
&=& { |S^{N-\ell-1}_1| \over |S^{N-1}_1| }  \, { (N-|V|^2)_+^{N-\ell-2 \over 2} 
\over N^{N-2 \over 2} }
\\
&= & { \alpha_{N-\ell} (N-|V|^2)\over N^{-\frac l 2} \alpha_N(N) }\, 
\frac{e^{- {|V|^2 \over 2}}}{(2\pi)^{ \frac l2}} \, {\bf 1}_{ |V|\le \sqrt{N} }
(1 + \OO(\frac{\ell^2}N)) , 
\eean
from which we deduce 
\bear
\theta_{N,\ell}(V) 
&=& 
  { 
Z'_{N-\ell}(\sqrt{N-|V|^2}) \over Z'_N(\sqrt{N}) } \, 
{ \alpha_{N-\ell} (N-|V|^2)\over N^{-\frac l 2} \alpha_N(N) } \, {\bf 1}_{
|V|\le \sqrt{N} }
(1 + \OO(\frac{\ell^2}N))  \nonumber
\\
&=&   \frac{\alpha_{N-\ell} (N-\ell)}{N^{-\frac l 2} \alpha_N(N)}
\frac{ e^{- \left( {\ell- |V|^2 \over \sqrt{N-\ell} \, \Sigma }
\right)^2 } + \OO((N-\ell)^{-1/2})} 
{ 1 +\OO (N^{-1/2})}
\,  {\bf 1}_{ |V|\le \sqrt{N} } (1 + \OO(\frac{\ell^2}N)) \nonumber
\\
&=& 
 \underbrace{\left( e^{- \left( {\ell- |V|^2 \over
\sqrt{N-\ell} \, \Sigma } \right)^2 } + \OO((N-\ell)^{-1/2}) \right)
}_{\theta_{N,\ell}^1 (V)}
 \, \underbrace{(1 + \OO(\frac{\ell^2}N)) \, {\bf 1}_{ |V|\le \sqrt{N}
}}_{\theta_{N,\ell}^2} \label{eq:Theta12}
\eear
where we have successively used  \eqref{eq:asympZprim},   \eqref{eq:asympZprimN}
 the definition of $\alpha_{N-\ell}(N-\ell)$, and a calculation yielding 
$$
\frac{\alpha_{N-\ell} (N-\ell)}{N^{-\frac l
2} \alpha_N(N)} = 1 + \OO(\ell^2/N).
$$  
It implies in particular the two following estimates on  $\theta_{N,\ell}$ which
will also be very useful in the proof of the next theorems
\beqn \label{estim:BoundsTheta}
\theta_{N,\ell}(v) \leq C \, {\bf 1}_{ |V|\le \sqrt{N}}, \quad  |
\theta_{N,\ell} (V) - 1 | \le 
{C \ell^2 \over N^{1/2}}  + C \,  {|V|^4
\over N^{1/2}} {\bf 1}_{|V| \ge N^{1/8}}.
\eeqn
Once they are proven, the conclusion follows since from the second one
\bean
\|F^N_\ell - f^{\otimes \ell} \|_1 
&=& \|(\theta_{N,\ell}  - 1) \,  f^{\otimes \ell} \|_1
\\
&\le& {C \, \ell^2 \over N} \, \| f  \|_1  +   {C  \over N^{1/2}} \, \| v^6 f 
\|_1.
\eean

\medskip
It only remains to prove the estimates~\eqref{estim:BoundsTheta}. The first uniform
estimate in \eqref{estim:BoundsTheta} is clear from~\ref{eq:Theta12} since $\| \theta^1_{N,\ell}
\|_\infty$ and $\theta^2_{N,\ell}$ are also uniformly bounded. For the second
estimate, we first control
\bean
| \theta^1_{N,\ell} (V) - 1 |
&=& | \theta^1_{N,\ell} (V) - 1 | \, {\bf 1}_{|V| \le N^{1/8}} + |
\theta^1_{N,\ell} (V) - 1 | \, {\bf 1}_{|V| \ge N^{1/8}} 
\\ 
&\le& |  2 \left( {\ell- |V|^2 \over \sqrt{N-\ell} \, \Sigma } \right)^2  +
\OO(N^{-1/2})  | \, {\bf 1}_{|V| \le N^{1/8}} + C \,  {|V|^4 \over N^{1/2}}
{\bf 1}_{|V| \ge N^{1/8}} 
\\
&\le& {C \ell^2 \over N^{1/2}} \, {\bf 1}_{|V| \le N^{1/8}} + C \,  {|V|^4
\over N^{1/2}} {\bf 1}_{|V| \ge N^{1/8}},
\eean
which implies a similar bound for $\theta_{N,\ell}$ since
\bean
| \theta_{N,\ell} (V) - 1 |
& \le &  | \theta^2_{N,\ell}(V)| \, | \theta^1_{N,\ell} (V) - 1 | +  |
\theta^2_{N,\ell}(V) - 1 | 
\\
& \le & C \, | \theta^1_{N,\ell} (V) - 1 |  + C \frac{\ell^2}N
\\
&\le& {C \ell^2 \over N^{1/2}} + C \,  { |V|^4
\over N^{1/2}} {\bf 1}_{|V| \ge N^{1/8}}.
\eean
This concludes the proof.
 \qed

\medskip
\noindent
{\sc Proof of \eqref{eq:ChaosEstimFN} in Theorem~\ref{theo:EmpiricalMeasure}. }    
The proof of the two last estimates in \eqref{eq:ChaosEstimFN} follows from Theorem~\ref{theo:[fotimesN]fchaos}
together with \eqref{eq:EquivChaosPart3} and  \eqref{eq:EquivChaosPart4}. \qed


\subsection{Improved chaos for conditioned tensor products on the Kac's spheres.}
\label{sec:ImprovedChaosKac}

In this section, we aim to prove rate of chaoticity for stronger notions of chaos for the sequence $(F^N)$ defined in the preceding section. 
Let us first recall   the notion of {\it entropy chaos} and {\it Fisher information chaos }
  in the context of the ``Kac's spheres" as they have been yet defined in the introduction. For $f \in \PPP(E)$ smooth enough, we define the usual relative entropy and usual relative Fisher information 
$$
H(f|\gamma) :=\int_E  u \, \log u \, \gamma \, dv, \quad
I (f|\gamma)  := \int_E {|\nabla u|^2 \over  u}  \, \gamma \, dv, \quad u := f/\gamma,
$$
and similarly for  $G^N \in \Ps(\KK\SS_N)$, we define the (normalized) relative entropy and relative Fisher information 
$$
H(G^N|\sigma^N) := {1 \over N} \int_{\KK\SS_N} g^N \log g^N \, d\sigma^N, \quad
I (G^N|\sigma^N)  := {1 \over N} \int_{\KK\SS_N}  {|\nabla g^N|^2 \over  g^N}  \, d\sigma^N,
$$
where $g^N := {dG^N \over d\sigma^N}$ stands for the Radon-Nikodym derivative of $G^N$ with respect to $\sigma^N$.

\begin{defin} We say that a sequence  $(G^N)$ of $\PPP(\KK\SS_N)$ is  

i) $f$-entropy chaotic if  $G^N_1 \wto f$ and
$$
H(G^N|\sigma^N) \to H(f|\gamma),
$$

ii)  $f$-Fisher information chaotic if  $G^N_1 \wto f$ and
$$
I(G^N|\sigma^N) \to I(f|\gamma).
$$
\end{defin}
It is worth emphasizing again that our definition is slightly different (weaker)
that the corrseponding definition in \cite{CCLLV}. But they are in fact
equivalent as we shall see in next section
(Theorem~\ref{theo:LesDiffChaosSKN}).

\begin{theo}\label{theo:EntropCvgceKac} For any  $f \in \PPP_6(\R) \cap
L^p(\R)$, $p > 1$, satisfying the moment
assumptions~\eqref{hyp:[fotimesN]fchaos} of Theorem~\ref{theo:[fotimesN]fchaos},
 the corresponding conditioned product sequence of measures
$(F^N)$ defined by \eqref{def:CondProd} is $f$-entropy chaotic. More precisely,
there exists $C = C(p,\|f \|_{L^p}, M_6(f))$ such that  
\beqn\label{ineq:EntropCvgce1}
|H(F^N|\sigma^N) - H(f|\gamma) |\le {C \over \sqrt{N}} .
\eeqn
\end{theo}

\medskip\noindent
{\sc Proof of Theorem~\ref{theo:EntropCvgceKac}. } With the notation
$F^N :=  [f ^{\otimes N}]_{\KK\SS_N}$, we write for any $N\ge 1$
\bean 
H( F^N | \sigma^N) 
&=& {1 \over N}\int_{\KK\SS_N} \Bigl( \log  {f^{ \otimes N}\over Z'_N(f) \,
\gamma^{\otimes N } }  \Bigr) \, dF^N
\\
&=&\int_\R \Bigl( \log  {f \over \gamma}  \Bigr) \, F^N_1 - {1 \over N}\log
Z'_N(f).
\eean	
Thanks to the bound \eqref{eq:asympZprim} on $Z'_N(f)$ which implies that
$(Z'_N(f))$ is bounded, we deduce    
\bean
H(F^N|\sigma^N) 
= \int_{\R} F^N_1 \left( \log {f \over \gamma} \right)   + \OO (1/N).
\eean
Recalling the notation  $\theta_N := \theta_{N,1}$ defined in
\eqref{eq:deftheta}
and the estimates~\eqref{estim:BoundsTheta} it satisfies, 
we may then write 
\bean
H(F^N|\sigma^N) 
=  H(f|\gamma) + \underbrace{\int_{\R} (\theta_N - 1) \, f\, \left( \log {f
\over \gamma} \right)}_{=: T}   + \OO (1/N),
\eean
with 
$$
|T| \le \underbrace{C \, \int_{\R} |\theta_N - 1| \, f \,  (1 + |v|^2) \,
dv}_{=:
T_1} +  \underbrace{\int_{\R} |\theta_N - 1| \, f \, | \log f |\, dv}_{=: T_2}. 
$$
In order to deal with $T_1$, we use the second estimate
of~\eqref{estim:BoundsTheta} and get
$$
T_1 
\le {C \over N^{1/2}} \int_{\R^d}   \langle v \rangle^2 \, f \, dv +  {C \over
N^{1/2}} \int_{\R^d}   \langle v \rangle^6 \, f \, dv = {C \over N^{1/2}}.  
$$
In order to deal with $T_2$, we make the more sophisticated (but standard)
splitting: for any  $N, R, M \ge 1$, we write 
\bean
T_2 
&\le& \int_{B_R} |\theta_N-1|\, f   |\log f | + C_\theta \int_{B_R^c} f \, |\log
f| 
 \\
&\le& \sup_{B_R} |\theta_N - 1| \, C_f
+ C_\theta \int_{B_R^c} f  (\log f)_+ \, {\bf 1}_{f \ge M}  + C_\theta
\int_{B_R^c} f ( \log f )_+ \, {\bf 1}_{M \ge f \ge 1}\\
&& +  C_\theta \int_{B_R^c} f ( \log f)_- \, {\bf 1}_{1 \ge f \ge e^{-|v|^{2}}}
+ C_\theta \int_{B_R^c} f  ( \log f)_-  \, {\bf 1}_{e^{-|v|^{2}} \ge f \ge 0}.
\eean
For the second term, we write $ f \, ( \log f )_+ \le  f^{(1+p)/2}  \le 
f^p/M^{(p-1)/2}$ on  $\{f \ge M\}$. 
For the third term, we write  $ f \, ( \log f )_+ \le f \, \log M \le f \, (\log
M) \, |v|^6/R^6$ on $\{f \le M, \,\, |v|\ge R\}$.
For the fourth term, we write $\log f \ge - |v|^{2}$ on $\{f \ge \exp(-|v|^{2})
\}$, and thus $ f (\log f)_-  \le f \, |v|^2 \le f \, |v|^6 / R^4$ on 
$\{1 \ge f \ge e^{-|v|^2}, \,\, |v|\ge R\}$.
For the last term, we write $ f \, ( \log f )_- \le 4 \, \sqrt{f}$ on $\{0 \le f
\le 1 \}$, and thus $ f \, ( \log f )_- \le 4 \, e^{-|v|^2/2}$ 
on  $\{  e^{-|v|^2} \ge f \ge 0, \,\, |v| \ge R \}$. 
We deduce 
\bean
T_2 
&\le& C_f \sup_{B_R} |\theta_N - 1| + C_\theta \left( {1 \over
M^{(p-1)/2}} + {(\log M)_+ \over R^6} + {1 \over R^4} + e^{-R} \right) 
\\
&\le& {C(\| f \|_p, M_6(f)) \over N^{1/2}},
\eean
with the choice $R=N^{1/8}$ (which allows to use the second estimate
of~\eqref{estim:BoundsTheta}), and then $M^{(p-1)/2} = R^6$. 
\qed

\medskip
Before stating a similar result with the Fisher information, we introduce a
notation: the gradient on the Kac's spheres $\KK\SS_N$  will be denoted by
$\nabla_\sigma$
$$
\nabla_\sigma F(V)  := P_{V^\perp} \nabla F(V) = \left( Id -
\frac{ V \otimes V}{|V|^2}\right) \nabla F(V) = \nabla F(V) - \frac{V
\cdot \nabla F(V) }{N} V,
$$ 
if $F$ is a smooth function on $\R^N$. $P_{V^\perp}$ stands for the projection
on the hyperplan perpendicular to $V$. We will use many times that 
\beqn\label{eq:nablasigma}
\nabla \left[ F\left( \frac  V{|V|} \right)\right] = \frac1{|V|}
P_{V^\perp} \nabla F \left( \frac V{|V|} \right) = \frac1{|V|}
 \nabla_\sigma F \left( \frac V{|V|} \right). 
\eeqn

\smallskip
\begin{theo} \label{FNFisherchaos}
For any  $f \in \PPP_6(\R)$, satisfying the moment
assumptions~\eqref{hyp:[fotimesN]fchaos} of Theorem~\ref{theo:[fotimesN]fchaos},
 the corresponding conditioned product sequence of measures
$(F^N)$ defined by \eqref{def:CondProd} satisfies
$$\sup_{N \in \N} I(F^N | \sigma_N) < +\infty
$$
if $I(f)< +\infty$. If moreover 
$$
 \int_{\R} \frac{ f'(v)^2}{f(v)} \langle v \rangle^2  \,dv < + \infty,
 $$ 
the sequence $F^N$ is Fisher information chaotic.
\end{theo}

{\sc Proof of Theorem~\ref{FNFisherchaos}.}
We only proof the second point. The first point (boundedness of the Fisher
information)  can be deduced from the above proof. It suffices in fact to use
the simple bound $| \nabla_\sigma G | \le | \nabla G|$ instead of
equality~\eqref{eq:NormNablaSigma}. 

Remark also that the bound on the Fisher
information implies that $f$ is continuous and uniformly bounded since $E=\R$. Therefore, the $L^p$  (for $p>1$)
assumption which is necessary in theorem~\ref{theo:[fotimesN]fchaos} 
is implied by our bound on the Fisher information.  We can therefore apply the
estimates~\eqref{estim:BoundsTheta} on the quantity $\theta_{N,i}$ for $i=1,2$ 
defined in~\eqref{eq:deftheta}. They imply in particular that $\| \theta_{N,i}
\|_\infty$ is uniformly bounded and that $\theta_{N,i}$ converges point-wise
towards $1$. We start with the formula
$$
I(F^N | \sigma^N)  
= \frac1N \int_{\KK\SS_N} | \nabla_\sigma \ln \left( \frac{f^{\otimes
N}}{\gamma^{\otimes N}}\right) |^2 \,F^N(dV) .
$$
As $\nabla_\sigma$ is the projection on the Kac's spheres of the usual gradient,
we have from \eqref{eq:nablasigma} for any function  $G$ on $\R^N$ 
\beqn \label{eq:NormNablaSigma}
\left| \nabla_\sigma G(V) \right|^2 
 = \left|  \nabla G(V)  \right|^2 - \frac1N \left| V \cdot \nabla G(V)
\right|^2 . 
\eeqn
Using this with $G =  \ln  \left( \frac {f^{\otimes N}} {\gamma^{\otimes N}} \right)
$ in the Fisher information formula, it comes %
$$
I(F^N | \sigma^N)  =  \frac1N \int_{\KK\SS_N} \left|\nabla 
 \ln  \left( \frac {f^{\otimes N}} {\gamma^{\otimes N}} \right)
  \right|^2 \,F^N(dV)  - \frac1{N^2} \int_{\KK\SS_N}
\left| V \cdot \nabla  \ln  \left( \frac {f^{\otimes N}} {\gamma^{\otimes N}} \right) \right|^2 
\,F^N(dV) .
$$
Recalling that $F^N_1 = f \, \theta_{N,1}$ from \eqref{eq:deftheta}, 
by symmetry, the first term in the right hand side is equal to
\beqn \label{eq:FishSigmaDecomp}
\int_\R  \left| \partial_v  \ln \frac {f(v)}{\gamma(v)} \right|^2 F^N_1(dv) =
I(f | \gamma) + \int_\R \left| \frac {\nabla f(v)}{f(v)} + v  \right|^2 (
\theta_{N,1}(v) -1) f(v) \,dv.
\eeqn
The last term goes to zero from the hypothesis on $f$, the uniform bound
$|\theta_{N,1}| \leq C$ and the pointwise convergence of $\theta_{N,1}$ to $1$. 
To handle the second term in the RHS of~\eqref{eq:FishSigmaDecomp}, we compute
\bean
\frac1{N^2}\left| V \cdot \nabla \ln \left( \frac f \gamma \right)^{\otimes N}
\right|^2 & = & \frac1{N^2}\left( \sum_{i=1}^N  v_i \left[\ln \frac f \gamma
\right]' (v_i)\right)^2 \\
& = &  \frac1{N^2} \sum_{i=1}^N  v_i^2 \left(\left[\ln \frac f \gamma \right]'
(v_i)\right)^2 + \frac1{N^2} 
\sum_{i \neq j}^N  v_i v_j \left[\ln \frac f \gamma \right]' \!\!\!\!(v_i)
\left[\ln \frac f \gamma \right]' \!\!\!\!(v_j).
\eean
After integration, it comes thanks to the symmetry of $F^N$
\[ \begin{split}
\frac1{N^2} \int_{\KK\SS_N}  \left| V \cdot \nabla \ln \left( \frac f \gamma
\right)^{\otimes N} \right|^2   \,F^N(dV) 
& =  \frac1N \int_{\R}   v^2 \left(\left[\ln \frac f \gamma \right]'
(v)\right)^2  \, F^N_1(dv) \\
 & + \frac{N-1}N \int_{\R^2}   v_1 v_2 \left[\ln \frac f \gamma \right]'
\!\!\!\!(v_1) \left[\ln \frac f \gamma \right]' \!\!\!\!(v_2)   \,
F^N_2(dv_1,dv_2).
\end{split} \]
Using the uniform bound $F^N_1(v) = \theta_{N,1}(v)\, f(v)  \leq C f(v)$, and
the hypothesis on $f$, we obtain that the first term of the r.h.s. is bounded by
$\frac C N$. The second term denoted by $R_2(N)$ is equal to
\bean
R_2(N) & = & \frac{N-1}N \int_{\R^2}   v_1 v_2 \left[\ln \frac f \gamma \right]'
\!\!\!\!(v_1) \left[\ln \frac f \gamma \right]' \!\!\!\!(v_2)   \, f(v_1)f(v_2)
\,  dv_1dv_2 \\
& &  + \frac{N-1}N \int_{\R^2}   v_1 v_2 \left[\ln \frac f \gamma \right]'
\!\!\!\!(v_1) \left[\ln \frac f \gamma \right]' \!\!\!\!(v_2)  
(\theta_{N,2}(v_1,v_2) - 1) \, f(v_1)f(v_2) \,  dv_1dv_2 \\
& =& \frac{N-1}N  \left( \int_{\R}   \left( v f'(v)  + v^2  f(v) \right)   \, 
dv \right)^2 + R_3(N) = R_3(N),
\eean
after an integration by parts and because of the equality $\int v^2 \,f(dv) =1$.
The term $R_3(N)$ goes to zero by dominated convergence since
\[ \begin{split}
\int_{\R^2} v_1 v_2 \left|  \left[\ln \frac f \gamma \right]' \!\!\!\!(v_1)
\left[\ln \frac f \gamma \right]'   \!\!\!\!(v_2)  \right|  & f(v_1)f(v_2) \, 
dv_1dv_2  =  \left( \int_\R   v \left| \left[\ln \frac f \gamma \right]'  
\!\!\!\!(v) \right| f(v)\,  dv \right)^2 \\
& \le   I(f | \gamma ) \int_\R v^2\, df  .
\end{split} \]
This concludes the proof.
\qed


\subsection{Chaos for arbitrary sequence of probability measures on the Kac's spheres.}
\label{sec:Entropie-Gibbs}
In that last section, we  aim to  present the relationship between Kac's chaos, entropy chaos and Fisher information chaos
in the Kac's spheres framework.

We begin  with a result which is the analogous for probability measures on the
Kac's spheres to the lower semi continuity of the Entropy and Fisher information yet established
on product spaces.

\begin{theo}\label{theo:EntropGammaCvgce} 
For any  sequence $(G^N)$ of $\PPP(\KK\SS_N)$  
such that  $G^N_j \wto G_j$ weakly in  $\PPP(E^j)$, there holds 
$$
H(G_j|\gamma^{\otimes j}) \le \liminf H(G^N|\sigma^N), 
\qquad 
I(G_j|\gamma^{\otimes j}) \le \liminf I(G^N|\sigma^N) .
$$
\end{theo}

For the proof, we shall need the following integration by parts formula on the
Kac' spheres, which proof is postponed to the end the proof of
Theorem~\ref{theo:EntropGammaCvgce} .

\begin{lem} \label{lem:IPPKac}
 Assume that $F$ (resp. $\Phi$) is a function (resp. vector field in $\R^N$)
on the Kac's spheres $\KK\SS_N$ with integrable gradient. Then the following
integration by part formula holds
\beqn \label{eq:IPPKac}
\int_{\KK\SS_N} \left[ \nabla_\sigma F(V) \cdot \Phi(V)  + F(V)
\diver_\sigma \Phi(V)  - \frac{N-1}N F(V) \Phi(V) \cdot V \right] d\sigma^N(V)
= 0 
\eeqn
where $\diver_\sigma$ stands for the divergence on the sphere, given by
$$
\diver_\sigma \Phi(V) :=  \sum_{i=1}^N \nabla_\sigma \Phi_i(V) \cdot e_i
=  \diver \Phi(V) - \sum_{i=1}^N \frac{V \cdot \nabla \Phi_i(V)}{|V|^2} v_i 
$$
where the  last formula is useful only if $\Phi$ is defined on a
neighborhood of the sphere.
\end{lem}

\noindent
{\sc Proof of Theorem~\ref{theo:EntropGammaCvgce}.}
We refer to \cite[Theorem 17]{CCLLV} for a proof of the inequality involving the
entropy and we give only the proof of the second inequality, which in fact relies
 on the characterization $I^{(3)}$ of the Fisher
information. Precisely, the previous Lemma~\ref{lem:IPPKac} can be used to get a
reformulation of the Fisher information relative to
$\sigma^N$ on the sphere
\bear
I_N (G^N | \sigma^N) &:=& \int_{\KK\SS_N} |\nabla_\sigma \ln G^N |^2 \,G^N(dV) 
= \sup_{\Phi \in C^1_b(\R^N)^N} \int_{\KK\SS_N} \left( \nabla \ln
 G^N \cdot \Phi - \frac{|\Phi|^2}4 \right)
\,G^N  \nonumber \\
& \hspace{-1cm}= & \hspace{-1cm} \sup_{\Phi \in C^1_b(\R^N)^N} \int_{\KK\SS_N} \left(
\frac{N-1}N \Phi(V) \cdot V - \diver_\sigma \Phi(V) - \frac{|\Phi(V)|^2}4
\right)
\,G^N(dV) \label{FisherSigmaN}.
\eear
Next applying the equality \eqref{eq:Fisherrelativedual} to the probability measure
$\gamma^{\otimes j}$, we get that for or any $\eps >0$, we can choose a $\varphi
\in C^1_b(\R^j)^j$ such that 
$$
\frac1j I_j(F^j|\gamma^{\otimes j})- \eps \le  \frac1j \int_{\R^j} \left( 
\varphi \cdot V_j -  \mathrm{div} \varphi - \frac{|\varphi|^2}4
\right) \,F^j(dV_j).
$$
Remark that the r.h.s. is quite similar to~\eqref{FisherSigmaN}.
With the notation $N = nj + r, \; 0 \le r < j$ and $V_N=
(V_{j,1},\ldots,V_{j,n},V_r)$, we define 
$$
\Phi(V_N) := (\varphi(V_{j,1},\ldots,\varphi(V_{j,n}),0) \in C^1_b(\R^N)^N,
$$ 
and use it in the equality~\eqref{FisherSigmaN}. We get
\bean
\frac1N I(G^N|\sigma^N) & \ge & \frac1N \int_{\KK\SS_N} \left(
\frac{N-1}N \Phi(V_N) \cdot V_N - \diver_\sigma \Phi(V_N) -
\frac{|\Phi(V_N)|^2}4 \right) \,G^N(dV_N) \\
& \ge & \frac nN \int_{\R^j} \left( \frac{N-1}N \varphi(V_j) \cdot V_j - 
\diver \varphi(V_j) - \frac{|\varphi(V_j)|^2}4 \right)
\,G^N_j(dV_j)  + \frac{R_\varphi(N)}N,
\eean 
where 
\bean
R_\varphi(N) & = & \frac1N \int \left(\sum_{i=1}^N [V \cdot \nabla \Phi_i(V) ]
v_i
\right) \,G^N(dV_N) 
= \frac1N \int \left( \sum_{i,\ell}^N \frac{\partial \Phi_i}{\partial v_\ell} v_i
v_\ell \right) \,G^N(dV_N) \\
&=& \frac nN \int \left( \sum_{i,\ell}^j \frac{\partial \varphi_i}{\partial v_\ell}
v_i v_\ell \right) \,G^N_j(dV_j)  = O(1),
\eean
if $\nabla \varphi$ decrease sufficiently quickly at infinity. Passing to the
limit, we get
$$
\liminf_{N \to +\infty} I(G^N | \sigma^N) \ge \frac1j \int_{\R^j} \left( 
\varphi \cdot V_j -  \mathrm{div} \varphi - \frac{|\varphi|^2}4
\right) \,F^j(dV_j) \ge I(F^j|\gamma^{\otimes j})- \eps
$$
which concludes the proof. \qed
 
\medskip \noindent
{\sc Proof of Lemma~\ref{lem:IPPKac}} As before, we will use the normalized norm
 $| V|_2 := \sqrt{ \frac1N \sum v_i^2}$. 
Choosing any smooth function $q$ on $(0,+\infty)$ with compact support, we
define
$$
w(V) :=  q(|V|_2) \, F\left( \frac V{|V|_2}\right) \, \Phi \left( \frac
V{|V|_2}\right).
$$
Its divergence is given by
\bean
\diver w &=& \frac{q'(|V|_2)}N F\left( \frac V{|V|_2}\right) \, \Phi \left(
\frac V{|
V|}\right) \cdot \frac V{|V|_2} 
+ \frac{q(|V|_2)}{|V|_2} \nabla_\sigma F\left( \frac
V{|V|_2}\right) \, \Phi \left( \frac V{| V|}\right) \cdot \frac V{|V|_2} \\
&& + \frac{q(|V|_2)}{|V|_2} \, F\left( \frac V{|V|_2}\right) \, \diver_\sigma
\Phi
\left( \frac V{|V|_2}\right).
\eean
Integrating this equality, and using polar coordinate, we get
\bean
0 &=&  \left( \int_{\KK\SS_N} \big[ \nabla_\sigma F(V) \cdot \Phi (V) + F(V) 
\diver_\sigma \Phi (V) \big]\,\sigma^{N}(dV) \right) \left( \int_0^\infty
q(r) r^{N-2}\,dr \right) \\
&& \hspace{1cm} +  \frac1N \left( \int_{\KK\SS_N}  F(V) \, \Phi (V) \cdot V
\,\sigma^{N}(dV) \right) \left( \int_0^\infty
q'(r) r^{N-1}\,dr \right).
\eean
Since $\int_0^\infty q'(r) r^{N-1}\,dr = -(N-1) \int_0^\infty
q(r) r^{N-2}\,dr$, we obtain
$$
\int_{{\KK\SS_N}} \left[ \nabla_\sigma F(V) \cdot \Phi(V)  + F(V)
\diver_\sigma \Phi(V)  - \frac{N-1}N F(V) \cdot\Phi(V) \cdot V \right] d\sigma^N(V)
= 0,
$$
which is the claimed result.
\qed

\medskip
 The next theorem will be the key estimate in the proof of the variant  of
Theorem~\ref{theo:LesDiffChaos} adapted to the Kac's spheres. It relies on the HWI
inequality on the Kac's spheres, which allows to quantify the convergence of the
relative entropy.

\begin{theo}\label{theo:Kac+FisherToEntropy}
Consider $(G^N)$ a sequence of  $\PPP(\KK\SS_N)$ which is  $f$-chaotic, $f \in
\PPP(E)$. Assume furthermore that 
$$
M_k(G^N)^{\frac1k} \le K \;\text{for } k \ge 6, \quad 
\text{and}\quad I(G^N|\sigma^N) \le K.
$$
Then $f$ satisfies  $M_k(f) < \infty$, $I(f) < \infty$, and $(G^N)$ is 
$f$-entropy chaotic. More precisely, there exists  $C_1 := C_1(K)$ 
and for
any $\gamma_2  < \frac18 \frac{k-2}{k+1}$ a constant  $C_2(\gamma_2)$ such that
$$
|H(G^N|\sigma^N) - H(f|\gamma)|\le C_1 \, \Bigl( W_1(G^N,f^{\otimes
N})^{\gamma_1} + C_2 N^{-\gamma_2} \Bigr),
$$
with $\gamma_1 := 1/2-1/k$. 
\end{theo}

The proof uses the following estimate

\begin{theo}\label{theo:CarlenLiebBarthe} {\bf (\cite[Theorem 1.2]{CarlenLL2004}, \cite[Theorem 2]{BartheCEM2006}). } For any sequence $(G^N)$
of $\PPP(\KK\SS_N)$, there hold for all $1 \le k \le N$
\begin{equation} \label{eq:add-HI-SN}
\begin{array}{l}
\ds  H(G^N_k|\sigma^N_k) \le 2  \, \frac Nk \Bigl[  \frac N k \Bigr]^{-1} \, H(G^N|\sigma^N), \\[1em]
\ds \int    \Bigl( \bigl| \nabla \ln g^N_k \bigr|^2 -  \frac1N \bigl[ \nabla \ln g^N_k  \cdot V_k \bigr]^2  \Bigr) \, G^N_k(dV_k) \le 2 \, \frac Nk \Bigl[  \frac N k \Bigr]^{-1} \, I(G^N|\sigma^N),
\end{array}
\end{equation}
where $gN_k := \frac{d G^N_k}{\sigma^N_k}$ stands for the Radon-Nikodym derivative of $G^N_k$ with respect to $\sigma^N_k$. 
\end{theo}

The first  inequality on the the entropy in \eqref{eq:add-HI-SN} was first proved in~\cite{CarlenLL2004} with $k=1$. 
It was generalized in~\cite{BartheCEM2006} to the case $k \ge 2$ and to the Fisher information. Remark that in both references, the result are stated in a somewhat different formulation, involving marginal with $k$ variables still defined on the Kac's sphere $\KK\SS^N$. We choose to translate it in the formalism we use here. It is worth emphaszing that in the l.h.s of the second  inequality, it is not exaclty the Fisher information $I(G^N_k | \sigma^N_k)$ that appears: there is a additional factor, with weight $N^{-1}$, due to the fact that the Fisher Information on the sphere is defined thanks to projected gradient.

\medskip
\noindent
{\sc Proof of Theorem~\ref{theo:Kac+FisherToEntropy}. Step  1. } 
Applying the second inequality of~\eqref{eq:add-HI-SN} with $k=1$, we get
\beqn \label{eq:Barthe1}
\int    \Bigl( 1 - \frac{v^2}N \Bigr) \Bigl[  \bigl(\ln g^N_1\bigr)' \Big]^2   \, G^N_1(dv) \le 2 \, I(G^N|\sigma^N),
\eeqn
with the notation $g^N_1:= \frac{dG^N_1}{d \sigma_1^N}$.
If view of the definition of Fisher information by duality~\eqref{eq:Fisherrelativedual}, we can choose, a smooth function $\varphi$ with compact support such that
$$
I(f|\gamma) \le \int_\R \Bigl( - \varphi \bigl( \ln \gamma \bigr)' - \varphi' - \frac{\varphi^2}4 \Bigr) df+ 1.
$$
Next  the trick used in order to obtain~\eqref{eq:Fisherrelativedual}   also leads to
\begin{multline*}
\int \Bigl( 1 - \frac{v^2}N \Bigr) \Bigl[ - \varphi \bigl( \ln \sigma^N_1 \bigr)' - \varphi' - \frac{\varphi^2}4  \Bigr]^2 \, G^N_1(dv) 
+ \frac 2N \int v \varphi \, G^N_1(dv)  \\
 \le \int    \Bigl( 1 - \frac{v^2}N \Bigr) \Bigl[  \bigl(\ln g^N_1\bigr)' \Big]^2   \, G^N_1(dv). 
\end{multline*}
Using the strong convergence of $\sigma_1^N$ to $\gamma$ stated
in~\ref{lem:sigmaNell}, the weak convergence of $G^N_1$ to $f$, and the  inequality~\eqref{eq:Barthe1},  we may pass to the (inferior) limit and we get  
$$
I(f| \gamma) \le  1 + \liminf_{N \to +\infty} \int    \Bigl( 1 - \frac{v^2}N \Bigr) \Bigl[  \bigl(\ln g^N_1\bigr)' \Big]^2   \, G^N_1(dv)  \le 1 + 2K, $$
which leads to $I(f) \le 2K$ after a simple calculation. 
\Black
%
%
%
%
%
Then, introducing the restriction $F^N = f^{\otimes N}/ Z(\sqrt N) \sigma^N$ of
$f^{\otimes N}$
to $\KK\SS_N$  defined in \eqref{def:CondProd} 
and using point $i)$ of Theorem~\ref{FNFisherchaos}, we get
$$
\sup_N I(F^N|\sigma^N)  \le C_2.
$$
\smallskip
\noindent
{\sl  Step 2. } Because the Ricci curvature of the metric space $\KK\SS_N$ is
positive (it is $K:= (N-1)/N$) we may use the HWI inequality in weak
$CD(K,\infty)$ geodesic space (see \cite[Theorem 30.21]{VillaniOTO&N}) which
generalizes the standard HWI inequality  \eqref{ineq-HWI}  quoted in
Proposition~\ref{prop:HWI}.  However, we have to be careful, because it is now
valid with $\tilde W_2$ replaced by the MKW distance constructed with the
geodesic distance on the sphere, and not with the distance induced by the square norm
of $\R^N$. Fortunately, both distances are equivalent, and if we add a constant
$\frac\pi 2$ in the right hand side, we can still write the HWI inequality with
our usual distance $W_2$.
We then have 
$$
H(F^N|\sigma^N) -  H(G^N|\sigma^N)  \le  \frac\pi 2  \sqrt{I(F^N|\sigma^N)}  \, 
W_2(F^N,G^N),
$$
and
$$
 H(G^N|\sigma^N) - H(F^N|\sigma^N) \le  \frac\pi 2   \sqrt{I(G^N|\sigma^N)}  
\, W_2(F^N,G^N), 
$$
so that
$$
|H(F^N|\sigma^N) - H(G^N|\sigma^N)| \le  C_2 \, W_2(F^N,G^N) . 
$$
We rewrite it under the form
\bean
|H(G^N|\sigma^N) - H(f|\gamma)|
&\le&  C_3 \left[ W_2(G^N,f^{\otimes N}) +  W_2(F^N,f^{\otimes N}) \right] 
+  |H(F^N|\sigma^N) - H(f|\gamma)|.
\eean
For the first term, we have using inequality of Lemma~\ref{lem:ComparDistEN}
$$
W_2(G^N,f^{\otimes N}) \le 4 \,K \, W_1(G^N,f^{\otimes N})^{1/2-1/k}.
$$
For the second term, we have for any $\eps >0$
\bean
W_2(F^N,f^{\otimes N}) 
&\le& 4 \,K \, \Omega_N(F^N;f)^{1/2-1/k}
\\
&\le& 4 \,K \,  \Bigl( \Omega_\infty(F^N;f) + C_\eps \,  N^{- {1 \over
2+\eps+2/k}} \Bigr)^{1/2-1/k }
\\
&\le& C_\eps  \,  \Bigl( \Omega_2(F^N;f)^{{1 \over 2+\eps+1/k}}  + C_\eps \, 
N^{- {1 \over 2+\eps+2/k}} \Bigr)^{1/2-1/k }
\\
&\le& C_\eps  \,  N^{ -  {1/4-1/2k \over 2+\eps+2/k}},
\eean
where we have successively used Lemma~\ref{lem:ComparDistEN}, the inequality
\eqref{eq:EquivChaosPart3}, \eqref{eq:EquivChaosPart4} and 
Theorem~\ref{theo:[fotimesN]fchaos} in the case $d=1$ (and then $d'=
\max(d,2)=2$). The third and last term is bounded by $C \, N^{-1/2}$ thanks to
Theorem~\ref{theo:EntropCvgceKac}.
\qed 

\medskip
The lower semi continuity properties of Theorem~\ref{theo:EntropGammaCvgce} and
Theorem~\ref{theo:Kac+FisherToEntropy} allow us
 to give a variant of Theorem~\ref{theo:LesDiffChaos} in the framework of
probability measures with support on the Kac's spheres.

\begin{theo}\label{theo:LesDiffChaosSKN} Consider $(G^N)$ a sequence of
$\PPP_{\! sym}(\KK\SS_N)$ such that  $M_6(G^N_1)$ is bounded and 
 $G^N_1 \wto f$ weakly in  $\PPP(\R)$. 

\noindent
In the  list of assertions below, each one implies the assertion which follows:

(i) $(G^N)$ is $f$-Fisher information chaotic, i.e. $I(G^N | \sigma^N) \to I(f|\gamma)$, $I(f) < \infty$;

(ii) $(G^N)$ is  $f$-Kac's chaotic and $I(G^N |\sigma^N)$ is bounded;  

(iii) $(G^N)$ is $f$-entropy chaotic,   that is $H(G^N | \sigma^N) \to H(f|\gamma)$, $H(f) < \infty$;

(iv) $(G^N)$ is  $f$-Kac's chaotic. 

\end{theo}
 
 \noindent
{\sc Proof of Theorem~\ref{theo:LesDiffChaosSKN}.}
The proof is very similar to the one of Theorem~\ref{theo:LesDiffChaos}. $i) \Leftrightarrow ii)$ and 
 $iii) \Leftrightarrow iv)$ relies on the l.s.c. properties of Theorem~\ref{theo:EntropGammaCvgce}. 
 And $ii) \Leftrightarrow iii)$ uses Theorem~\ref{theo:Kac+FisherToEntropy}. We omit the details.
\qed

 \medskip
 We finally conclude this section with the proof of Theorem~\ref{theo:OpenPb}.

 \smallskip\noindent
 {\sc Proof of Theorem~\ref{theo:OpenPb}. } We only deal with the case $j=1$, but the general case $j \ge 1$ can be managed
 in a very similar way because we already know that $G^N_j \wto f^{\otimes j}$ weakly in $\PPP(E^j)$ thanks to Theorem~\ref{theo:Kac+FisherToEntropy}
 and Theorem~\ref{theo:LesDiffChaosSKN}. 
 With the notations of Theorem~\ref{theo:OpenPb}, we have to prove
 $$
 H(G^N_1 |f) = \int_E  \log (G^N_1/f) \,  G^N_1   \, \to \, 0 
 \quad\hbox{as}\quad N \to \infty.
 $$
First, we observe that since $G^N$ is symmetric and has support on the Kac's spheres, $M_2(G^N) =1$.
Moreover,
\bean
I(G^N_1 | \sigma^N_1) 
&=& \int_E | \nabla \log G^N_1 - \nabla \log \sigma^N_1|^2 \, G^N_1
\\
&=& I(G^N_1) + \int_E [ 2 \, \Delta  \log \sigma^N_1 + |\nabla  \log \sigma^N_1|^2 ]  \, G^N_1,
\eean
so that 
\bean
 I(G^N_1) 
 &\le& I(G^N_1 | \sigma^N_1) +   \int_E ( 2 \, \Delta  \log \sigma^N_1 + |\nabla  \log \sigma^N_1|^2 )_-  \, G^N_1.
\eean
We easily compute
\bean
&&2 \, \Delta  \log \sigma^N_1 + |\nabla  \log \sigma^N_1|^2 
=
\\
&&\qquad
= {N-3 \over 2} \, \left\{  2 {  (2 \, v)^2 /N^2 \over (1 - v^2/N)^2 } - 2 { 2 /N \over (1 - v^2/N) }  +  {   (2 \, v/N)^2  \over (1 - v^2/N)^2 }  \right\} 
\, {\bf 1}_{v^2 \le N}
\eean
and then 
\bean
(2 \, \Delta  \log \sigma^N_1 + |\nabla  \log \sigma^N_1|^2)_- &= & 
2 \,  {N-3 \over N} \,     { (4 \, v^2/N - 1)_-  \over   (1 - v^2/N)^2 }   
\, {\bf 1}_{v^2 \le N/4}
\\
&\le& 2 \,     { 1 \over   (1 - 1/4)^2 }   = {32 \over 9}. 
\eean
Thanks to the boundedness assumption \eqref{eq:MkGN1&IGN} we get that $I(G^N_1) \le C$ for some constant $C \in (0,\infty)$, 
and then  $I(G^N_1|\gamma) \le 2[ I(G^N_1) +  M_2(G^N_1) ] \le C$. 

\smallskip  
Next, we introduce the splitting  
\bean
 H(G^N_1 | f)
&=& \underbrace{ H(G^N_1|\gamma) - H(f|\gamma) }_{=: T_1} + \underbrace{\int_E (  f - G^N_1) \, \log \frac f \gamma  }_{=: T_2} 
\eean
and we show that $T_i \to 0$ for any  $i=1,2 $. For the first term $T_1$, using twice the HWI inequality we have 
$$
|T_1| \le \left( \sqrt{I(G^N_1|\gamma)} + \sqrt{I(f|\gamma)} \right) \, W_2 (G^N_1,f) \to 0
$$
because of the uniform bound on the Fisher information and of the convergence property 
$W_2 (G^N_1,f) \to 0$. That last convergence is a consequence of  \cite[Theorem 7.2 (iii) $\Rightarrow$ (i)]{VillaniTOT},  $G^N_1 \wto f$ weakly when $N \to \infty$ and $\langle G^N_1, v^2 \rangle = \langle f, v^2 \rangle $ for any $N \ge 1$ when $k= 2$, and it is a is a consequence of  \cite[Theorem 7.2 (ii) $\Rightarrow$ (i)]{VillaniTOT},  $G^N_1 \wto f$ weakly   as $N \to \infty$ and $M_k( G^N_1) \le C$ for any $N \ge 1$ when $k> 2$.

Before dealing with the last term, we remark that the bound on the Fisher
information of $f$ implies some regularity, precisely that
$\sqrt f$ and then $f$ are $\frac12$-H\"older.
Therefore $\ln \frac f \gamma$ is continuous and satisfies from the assumption~\eqref{eq:f>exp} the bound
$$
\left| \ln \frac f \gamma \right| \leq \ln \| f \|_\infty +  \alpha |v|^{k'} + |\beta|  + \frac{v^2}2 \leq C \langle v \rangle^{\max(k',2)}.
$$
We then conclude that $T_2 \to 0$ by using  \cite[Theorem 7.2 (iii) $\Rightarrow$ (iv)]{VillaniTOT} when $k=2$ and 
\cite[Theorem 7.2 (ii) $\Rightarrow$ (iv)]{VillaniTOT} when $k>2$.  
\qed


\section{On mixtures according to De Finetti, Hewitt and Savage}
\label{sec:mixtures}
\setcounter{equation}{0}
\setcounter{theo}{0}

In this section we develop a quantitative and qualitative approach concerning
the sequence of probability measures of $\Ps(E^N)$, $E \subset \R^d$,  in the general framework of 
convergence to ``mixture of probability measures" (here we do not assume chaos property).

\smallskip
Depending on the result, we will need some hypothesis on the set $E$ that we will make precise
in each statement. While in the first and second sections the results hold with great generality only
assuming that 

- $E$ is a Borel set of $\R^d$;

\noindent
we shall assume in the third and fourth sections that 

  - $E=\R^d$ or $E$ is an  open  set of $\R^d$ with smooth boundary in order that the strong maximum principle and the Hopf lemma hold (that we furthermore assume to be bounded in the third section); 
 
 \noindent
and we shall also assume in the  fourth section that 

- the normalized non relative HWI inequality \eqref{eq:GalHWI} holds in $E$ (e.g. it satisfies the assumptions of Proposition~\ref{prop:HWI}).

\subsection{The De Finetti, Hewitt and Savage theorem and weak convergence in $\PPP(E^N)$}
\label{subsec:DeFHS}

We begin by recalling the famous De Finetti, Hewitt and Savage theorem
\cite{deFinetti1937,HewittSavage} for which we state a quantified version that is maybe new. 

\begin{theo}\label{theo:H&S} Assume $E \subset \R^d$ is a Borel set. 
Consider a sequence $(\pi^j)$ of symmetric and compatible probability measures of
$\PPP(E^j)$, 
that is $\pi^j \in \Ps(E^j)$ and $(\pi^j)_{|E^\ell} = \pi^\ell$ for any $1 \le
\ell \le j$, and consider $(\hat \pi^j)$ the associated sequence of empirical distribution 
in $\PPE$ defined according to \eqref{def:hatFN}. For any $s>\frac d2$,
the sequence $(\hat \pi^j)$ is a Cauchy sequence for the distance $\WWs$, and
precisely
\beqn
\left[\WWs (\hat \pi^N, \hat \pi^M) \right]^2  \leq  \, 2 \| \Phi_s\|_\infty 
\left( \frac1M + \frac1N \right),
\eeqn
where $\Phi_s$ is the function introduced in Lemma~\ref{lem:H-s_poly}. In
particular, the sequence $(\hat \pi^j)$ converges towards some $\pi \in
\PPP(\PPP(E))$ with the speed $\WWs(\hat \pi^j,\pi) 
\leq   \frac C{\sqrt j}$. The limit $\pi$ is characterized by the relations
\beqn\label{eq:H&S1}
\forall  \, j \ge 1, \quad 
\pi^j = \pi_j := \int_{\PPP(E)}   \rho^{\otimes j} \,\pi(d\rho) \quad \hbox{in}\quad \PPP_{\! sym}(E^j), 
\eeqn
or in other words, with the notations of section~\ref{subsect:def}
\beqn\label{eq:H&S2}
\forall \, \varphi \in C_b(E^j) \quad
\langle \pi^j,\varphi \rangle = \int_{\PPP(E)} R_\varphi (\rho) \, \pi(d\rho).
\eeqn
Reciprocally, for any mixture of probability measures $\pi \in \PPP(\PPP(E))$, the sequence
$(\pi_j)$ of probability measures in $\PPP(E^j)$ defined by the second identity in \eqref{eq:H&S1}
is such that the $\pi_j$ are symmetric and compatible. 
\end{theo}

\noindent
{\sc Proof of Theorem~\ref{theo:H&S}.} We split the proof into two steps. 

\smallskip\noindent
{\sl Step 1. }  In order to estimate the distance between $\hat \pi^N$ and $\hat \pi^M$ we shall use as in the 
proof of Proposition~\ref{prop:WWenHs} the fact that $\|
\cdot\|^2_{H^{-s}}$ is a polynomial on $\PPP(E)$, but we have to choose a good
transference plan. Fortunately, their is at least one simple choice. The compatibility and symmetry
conditions on $(\pi^N)$ tell us that $\pi^{N+M}$ is an admissible transference
between $\pi^N$ and $\pi^M$. Using the symmetry of $\pi^{N+M}$ and the isometry
between ($E^N/\SN,w_1)$ and $(\PPP_N(E),W_1)$ stated in step 1 in the proof of 
Proposition~\ref{prop:equivW1-2}, we will interpret it as a transference plan $\tilde
\pi^{N+M}$ on $\PP_N(E) \times \PP_M(E)$ between  $\hat \pi^N$ and $\hat
\pi^M$. 
More precisely,  $\tilde \pi^{N+M}\in \PPP(\PPP(E) \times \PPP(E))$ is defined as the probability measure satisfying
$$
\forall \, \Phi \in C_b(\PPP(E) \times \PPP(E)) \quad \langle \tilde \pi^{N+M},
\Phi \rangle = \int_{E^N \times E^M} \Phi(\mu^N_X,\mu^N_Y) \, \pi^{N+M} (dX,dY).
$$
With that transference plane we have
\begin{align*}
\left[ \WWs (\hat \pi^N, \hat \pi^M) \right]^2  & \leq   \int_{\PPP(E) \times \PPP(E)}  \| \rho - \eta
\|^2_{H^{-s}} \,\tilde \pi^{N+M}(d\rho,d\eta) 
\\
& \le  \int_{\PPP(E) \times \PPP(E)} \biggl( \int_{\R^{2d}} \Phi_s(x-y) \, [ (\rho^{\otimes 2} -  \rho
\otimes  \eta) \cdot  \\
& \hspace{4cm}
+ (\eta^{\otimes 2} - \eta \otimes \rho)] (dx,dy) \biggr)  \,\tilde \pi^{N+M}(d\rho,d\eta),
\end{align*}
with the help of \eqref{eq:H-s}.  We can then compute
\bean
&&\!\!\!\!\!\!\!\!\! \left[ \WWs (\hat \pi^N, \hat \pi^M) \right]^2   \leq
\\
& \le  & \int  \left( \int_{\R^{2d}} \Phi_s(x-y) \, [(\mu^N_X)^{\otimes 2} - 
\mu^N_X \otimes  \mu^M_Y  ] (dx,dy) \right)  \,\pi^{N+M}(dX,dY) 
\\
& & +  \int  \left( \int_{\R^{2d}} \Phi_s(x-y) \, [ 
(\mu^M_Y)^{\otimes 2} - \mu_Y^M \otimes \mu_X^N ] (dx,dy) \right)  \,\pi^{N+M}(dX,dY) 
\\
& \le  & \int  \left( \frac1{N^2} \sum_{i,j=1}^N \Phi_s(x_i-x_j) - \frac1{NM}
\sum_{i,j'=1}^M \Phi_s(x_i-y_{j'})  \right)  \,\pi^{N+M}(dX,dY) 
\\
& & + \int  \left( \frac1{M^2} \sum_{i',j'=1}^M \Phi_s(y_{i'}-y_{j'}) - 
\frac1{NM} \sum_{i',j=1}^M \Phi_s(x_{i'}-y_j)
\right)  \,\pi^{N+M}(dX,dY) 
\\
& \le & \frac{\Phi_s(0)}N +
\frac{N-1}N \int \Phi_s(x-y) \,\pi^2(dx,dy) - \int \Phi_s(x-y) \,\pi^2(dx,dy) 
\\
&\ &  + \frac{\Phi_s(0)}M + \frac{M-1}M \int \Phi_s(x-y) \,\pi^2(dx,dy) -
\int \Phi_s(x-y) \,\pi^2(dx,dy),
\eean
and we conclude with 
\bean
\left[ \WWs (\hat \pi^N, \hat \pi^M) \right]^2
 & \le & \left( \frac1M + \frac1N \right) \left(  \Phi_s(0) -  \int \Phi_s(x-y) \,\pi^2(dx,dy) \right) 
\\
 & \le & 2 \| \Phi_s\|_\infty  \left( \frac1M + \frac1N \right).
\eean
The existence of the limit $\pi$ is due to the completeness  of $\PPP(\PPP(E))$. 

\medskip\noindent
{\sl Step 2. }
Now it remains to characterize the limit $\pi$.
We fix $j \in \N$, we denote by $\pi_j$ its $j$-th marginal defined thanks to 
the second identity in \eqref{eq:H&S1} and by $\hat \pi^N_j = (\hat\pi^N)_j$
the $j$-th marginal of the empirical probability measure $\hat\pi^N$ as defined in \eqref{def:FNj}. 
We easily compute 
\bean
\| \hat \pi^N_j -\pi_j\|^2_{H^{-s}} &= & \left\| \int_{\PPP(E)} \rho^{\otimes j}
\,\hat \pi^N(d\rho) -  \int_{\PPP(E)} \rho^{\otimes j}\, \pi(d\rho) \right\|^2_{H^{-s}} 
\\
& = & \inf_{\Pi \in \Pi(\hat \pi^N,\pi)} \left\|   \int_{\PPP(E)}  [\rho^{\otimes j} -
\eta^{\otimes j}]\, \Pi(d\rho,d\eta) \right\|^2_{H^{-s}} 
\\
& \leq &  \inf_{\Pi \in \Pi(\hat \pi^N,\pi)}   \int_{\PPP(E)}  \| \rho^{\otimes j} -
\eta^{\otimes j} \|^2_{H^{-s}} \, \Pi(d\rho,d\eta) 
\\
& = & \left[ \WWs (\hat \pi^N, \pi) \right]^2 \leq \frac C N.
\eean
Next we fix $s>\frac{jd}2$, so that using Sobolev embeddings on $\R^{jd}$, $\|
\varphi \|_\infty \leq C \| \varphi \|_{H^s}$ for any $\varphi \in H^s(\R^{jd})$,
which implies by duality that $\| \rho \|_{H^{-s}} \leq C \| \rho \|_{TV}$ for
any $\rho \in \PPP(\R^{jd})$. Using the Grunbaum lemma~\ref{lem:HSEquivQuant}
and the compatibility assumption $\pi^N_j = \pi^j$,
we get the
inequality
$$
\| \pi^j - \hat \pi^N_j \|_{H^{-s}} = \| \pi^N_j - \hat \pi^N_j \|_{H^{-s}}
\leq C \, \| \pi^N_j - \hat \pi^N_j \|_{TV} \leq \frac{C j^2}N.
$$
Combining the two previous inequalities leads to
$$
\| \pi^j - \pi_j \|_{H^{-s}} \leq \| \pi^j - \hat \pi^N_j \|_{H^{-s}} + \|
\hat \pi^N_j - \pi_j\|_{H^{-s}} \leq \frac C{\sqrt N} + \frac{C j^2}N,
$$
which implies the claimed equality in the limit $N \to +\infty$.
\qed

\bigskip
Let us now introduce some definitions. 
For $k > 0$, we define 
$$
\PPP_k(\PPP(E)) := \{\pi \in \PPP(\PPP(E)); \,\, M_k(\pi) := M_k(\pi_1) < \infty
\}
$$
and for $k,a > 0$, we define
$$
\BB\PPP_{k,a}(E^N) := \{ F \in \PPP(E^N); \,\, M_k(F_1) \le a \}.
$$

\begin{defin} \label{def:CVG_PPP_k}
For given sequences $(F^N)_N$ of $\Ps(E^N)$,  $(\pi_n)_n$ of $\PPP(\PPP(E))$ and $\pi \in \PPP(\PPP(E))$, we say that

\smallskip
- $(F^N)$ is bounded in $\PPP_k(E^N)$ if there exists $a > 0$ such
that $M_k(F^N_1) \le a$; 

\smallskip
- $(\pi_n)$ is bounded in $\PPP_k(\PPP(E))$ if there exists  $a > 0$
such that $M_k(\pi_{n,1}) \le a$; 

\smallskip
- $(F^N)$ weakly converges to $\pi$ in $\PPP_k(E^j)_{\forall \, j}$, we write $F^N \wto \pi$ weakly in $\PPP_k(E^j)_{\forall \, j}$,
 if  $(F^N)$ is bounded in $\PPP_k(E^N)$ and $F^N_j \wto \pi_j$ weakly in 
$\PPP(E^j)$ for any $j \ge 1$; 

\smallskip
- $(\pi_n)$ weakly converges to $\pi$ in $\PPP_k(\PPP(E))$ if  $(\pi_n)$ is
bounded in $\PPP_k(\PPP(E))$ and $\pi_n \wto \pi$ weakly in 
$\PPP(\PPP(E))$.  
\end{defin}

With that (not conventional) definitions, any bounded sequence in 
$\PPP_k(\PPP(E))$ is weakly compact in  $\PPP_k(\PPP(E))$, 
and for any sequence $(F^N)$ of probability measures of $\Ps(E^N)$ which is bounded in
$\PPP_k(E^N)$, $k > 0$, there exists a subsequence $(F^{N'})$ and a mixture of 
probability measures $\pi \in \PPP_k(\PPP(E))$ such that   $F^{N'} \wto \pi$ in $\PPP(E^j)_{\forall \, j}$. 

\smallskip
We now present a result about the equivalence of convergences for sequence of
$\PPP_{\!sym}(E^N)$, $N\to\infty$, without any chaos hypothesis. 

\begin{theo}\label{theo:CvgePENPPEquiv}  Assume $E \subset \R^d$ is a Borel set. 

\noindent
{\bf (1) }Consider $(F^N)$ a sequence of $\PPP_{\!sym}(E^N)$ 
and  $\pi \in \PPP(\PPP(E))$. The three following assertions are equivalent: 

\begin{itemize}

\item[(i)] $F^N \wto \pi$ in $\PPP(E^j)_{\forall \, j}$, that is $F^N_j \wto \pi_j$ weakly in $\PPP(E^j)$ for any $j \ge 1$; 

\item[(ii)] $\hat F^N \wto \pi$ weakly in $\PPP(\PPP(E))$; 

\item[(iii)] $W_1(F^N,\pi_N) \to 0$. 

\end{itemize}

\noindent
{\bf (2) }  For any $\gamma \in [\frac1{2d'},\frac1{d'})$  (recall that $d'=\max(d,2)$),  and any $k> \frac {d'}{\gamma^{-1} - d'}  \ge 1$, there exists a constant
$C = C(\gamma,d,k)$ such that the following estimate holds
\beqn\label{ineq:CvgePENPPEquiv}
\forall \, N \ge 1 \qquad  | W_1 (F^N,\pi_N) - \WW_1(\hat F^N,\pi)| \leq 
{C \, M_k(\pi_1)^{1/k} \over N^\gamma}.
\eeqn

\noindent
{\bf (3) } With the same notations as  in the second point, we have for any mixture of probability measures $\alpha,\beta \in
\PPP(\PPP(E))$
\beqn\label{ineq:MixtureW1WW1}
\WW_1(\hat \alpha_j,\alpha) \le {C M_k(\alpha_1)^{1/k} \over j^\gamma},
\eeqn
where $\hat \alpha_j$ is empirical probability distribution in $\PPE$ associated to the $j$-th marginal $\alpha_j \in \PPP(E^j)$, 
as well as 
\beqn\label{ineq:MixtureW1WW1-2}
\WW_1(\alpha ,\beta) - {C (M_k(\alpha_1)^{\frac1k} + M_k(\beta_1)^{\frac1k}) \over j^\gamma} \le
W_1(\alpha_j,\beta_j) \le \WW_1(\alpha ,\beta).
\eeqn

\end{theo}

\medskip\noindent
{\sc Proof of Theorem~\ref{theo:CvgePENPPEquiv}. } 
{\sl Step 1. }   Equivalence between (i) and (ii) is
classical. Let us just sketch the proof. For any $\varphi \in C_b(E^j)$ we have
from the Grunbaum lemma  recalled  in  Lemma~\ref{lem:HSEquivQuant} that
\bean
\langle \hat F^N, R_\varphi \rangle 
&=& \langle F^N, \widetilde{ \varphi \otimes {\bf 1}^{\otimes N-j} }\rangle +
\OO ( j^2/N) 
\\
&=& \langle F^N_j,   \varphi  \rangle + \OO (j^2/N) .
\eean
We deduce that the convergence
$\langle \hat F^N, R_\varphi \rangle \to \langle \pi, R_\varphi \rangle$
is equivalent to the convergence $\langle F^N_j,   \varphi  \rangle \to \langle \pi_j,   \varphi  \rangle$
since that $ \langle \pi, R_\varphi \rangle =  \langle \pi_j,   \varphi  \rangle$ thanks to Theorem~\ref{theo:H&S}. 

Therefore, $i)$ is equivalent to the convergence $\langle \hat F^N, \Phi \rangle \to \langle \pi, \Phi \rangle$ for any polynomial function $\Phi \in C_b(\PPP(E))$. 
But now, the family of probability measures $\hat F^N$ (and $\pi$) belongs to the compact subset of
$\PPE$
$$
\mathcal K := \{ \alpha  \in \PPE, \; \text{s.t. } \alpha_1 = F_1 \},
$$
and also any converging subsequence $\hat F^{N'}$ should converge weakly towards
a probability measure $\tilde \pi$ having the same marginals as $\pi$. Since by
Theorem~\ref{theo:H&S} marginals uniquely characterize a probability measure on $\PPE$,
it implies $\tilde \pi = \pi$ and then weak convergence against polynomial
function implies the standard weak convergence of probability measures ii).

\smallskip\noindent
It is classical that the MKW distance is a metrization of  the weak convergence of measures. 
Even in that "abstract"  case, (ii) is equivalent to $\WW_1(\hat F^N,\pi) \to 0$ (recall that the distance chosen in order
to define $\WW_1$ is bounded).
Thus, for sequences having a bounded moment $M_k(F^N_1)$ for some $k>0$, the equivalence
between (ii) and (iii) will be a consequence of 
\eqref{ineq:CvgePENPPEquiv}. For sequences for which no moment $M_k$ is bounded, 
the same conclusion is true. The correct argument still relies on a version of  inequality
\eqref{ineq:CvgePENPPEquiv}, with a slower and less explicit rate of
convergence, which can be obtained from an adaptation of
Lemma~\ref{lem:ComparDistances}.

\smallskip\noindent
{\sl Step 2. } We now prove \eqref{ineq:CvgePENPPEquiv}. 
For $\hat \pi_N$ we have the following representation:  
\beqn\label{eq:relationHatpiN}
\hat \pi_N = \widehat{\int \rho^{\otimes N} \, \pi(d\rho)} =  \int \widehat{\rho^{\otimes N}} \, \pi(d\rho).
\eeqn
Thanks to Proposition~\ref{prop:equivW1-2}, we may compute 
\bean
| W_1(F^N,\pi_N) &-& \WW_1(\hat F^N,\pi) |
= | \WW_1(\hat F^N,\hat \pi_N)  - \WW_1(\hat F^N,\pi) |
\\
&\le& \WW_1( \hat \pi_N,\pi) 
= \WW_1\left(  \int_{\PPP(E)}   \widehat{\rho^{\otimes N}}   \pi(d\rho) , \int_{\PPP(E)}  \delta_\rho \, \pi(d\rho) \right)
\\
&\le& \int_{\PPP(E)}   \WW_1\left(   \widehat{\rho^{\otimes N}} ,  \delta_\rho \right)  \pi(d\rho) 
=  \int_{\PPP(E)}  \Omega_\infty(\rho) \,   \pi(d\rho) ,  
\\
 &  \le & \frac{C(d,\gamma,k)}{N^{\gamma}} \int_{\PPP(E)} 
M_k(\rho)^{1/k} \, \pi(d\rho) \le \frac{C(d,\gamma,k)}{N^{\gamma}}
M_{k}(\pi_1)^{1/k}, 
\eean
where we have successively used the triangular inequality for the $\WW_1$ distance, the relation \eqref{eq:relationHatpiN}, 
the convexity property of the $\WW_1$ distance and the definition of the chaos measure $\Omega_\infty$. 
We also used the bound \eqref{rem:Omega0FNf} and the Jensen inequality (recall that $1/k \in (0,1]$)  in the last line.  

\smallskip\noindent
{\sl Step 3. } We now prove the third point. For the first inequality, choose $s = \frac1{2\gamma} - \frac d{2k}$. Then by our assumptions, $s>\ \max(1,\frac d2)$ and we can apply  
Lemma~\ref{lem:ComparDistPPE} on the comparison
of distances in $\PPE$  and Theorem~\ref{theo:H&S} to get
$$
\WW_1(\hat \alpha_j, \alpha) 
\leq C \ M_k(\alpha_1)^{\frac 1k} \, \WWs(\hat \alpha_j,\alpha)^{\frac{2k}{d+2ks}} 
\leq \frac{C \,  M_k(\alpha_1)^{\frac1k} }{j^\gamma} 
$$
For the first part of the second inequality \eqref{ineq:MixtureW1WW1-2} we write 
$$
\WW_1(\alpha,\beta) \le \WW_1(\alpha,\hat\alpha_j) + \WW_1(\hat\alpha_j,\hat\beta_j) + \WW_1(\hat\beta_j,\beta),
$$
we use the inequality  just proved above and the identity \eqref{prop:equivW1-2}. The second part of the second inequality  \eqref{ineq:MixtureW1WW1-2} is
a mere application of Lemma~\ref{lem:WassNleqW1}. 
\qed


\subsection{Level-3 Boltzmann entropy  functional for mixtures } 

In this section we recover some well known results on the Boltzmann entropy for mixture of 
probability measures  as stated in \cite{ArkerydCI99} and proved by Robinson and Ruelle in \cite{RR}. However our proof differs 
from the one of \cite{RR}, and in particular it does not use the abstract representation result of Choquet and Meyer \cite{ChoquetMeyer1963} but an abstract Lemma  \ref{lem:niv3abst} that we introduce for our purposes. 

\smallskip
Let us assume that  $E \subset \R^d$ is a Borel set and let us fix a real number  $m>0$. 
Then, for any $\pi \in \PPP_m(\PPP(E))$ we define 
\beqn\label{eq:defHH}
\HH(\pi) := \int_{P(E)} H(\rho) \, \pi(d\rho), 
\eeqn
where $H$ is the Boltzmann's entropy defined on $\PPP_m(E)$.

\begin{theo}
\label{th:Rob&Ruelle} \mbox{}
${\bf (1)}$
The functional $\HH
: \PPP_m(\PPP(E)) \to \R \cup \{\infty \}$ is proper, affine and l.s.c. with
respect to the weak convergence in $\PPP_m(\PPP(E))$. 
Moreover, for any  $\pi \in \PPP_m(\PPP(E))$,  there holds
\beqn\label{eq:Rob&Ruelle}
\HH(\pi)  = \sup_{j \in \N^*} H(\pi_j) 
= \lim_{j\to\infty} H(\pi_j),
\eeqn
where $\pi_j$ is the j-th marginal of $\pi$ defined in Theorem~\ref{theo:H&S} and $H$ is the normalized Boltzmann's entropy defined on $\PPP_m(E^j)$
for any $j \ge 1$.

\smallskip
${\bf (2)}$
Consider $(F^N)$ a sequence of  $\PPP_{\!sym}(E^N)$ and $\pi \in \PPP_m(\PPP(E))$ such that $F^N \wto \pi$ weakly in  $\PPP_m(E^j)_{\forall j}$. Then 
\beqn\label{ineq:HFNtoHHpi}
\HH(\pi) \le \liminf_{N \to \infty} \, H (F^N).
\eeqn
\end{theo}

The proof of Theorem~\ref{th:Rob&Ruelle} uses the two following lemmas.

\begin{lem}\label{lem:RRH2}   For any  $\pi \in \PPP_m(\PPP(E))$ we define
$$
\HH'(\pi)  := \sup_{j \in \N^*} H(\pi_j) . 
$$
The functional $\HH'
: \PPP_m(\PPP(E)) \to \R \cup \{\infty \}$ is  affine, proper, and l.s.c. for the weak
convergence, and 
\beqn\label{eq:RRH'lim}
\HH'(\pi) =  \lim_{j\to\infty} H(\pi_j).
\eeqn
\end{lem}

The 
 proof of Lemma~\ref{lem:RRH2} is classical. For the sake of completeness we nevertheless present it.

\smallskip\noindent
{\sc Proof of Lemma~\ref{lem:RRH2}. } 
   Thanks to \eqref{lem:defENtropB}, for any $j \ge 1$, we have 
$$
H(\pi_j)  \ge \log c_m  - \int_E |v|^m \, d\pi_1
$$
so that $\HH'$ is proper on $\PPP_m(\PPP(E))$. It is also l.s.c. as the supremum of l.s.c. functions, since $H_j$ is l.s.c. on $\PPP_m(E^j)$ as it has been recalled in Lemma~\ref{lem:defENtropB} and since the inequality of the right of \eqref{ineq:MixtureW1WW1-2} shows that $\pi \mapsto \pi_j$ is also continuous for the weak convergence of measures. 

\smallskip
As a second step, we establish \eqref{eq:RRH'lim}. 
For any fixed $\ell \ge 1$ and any $j \ge \ell$ we introduce the  Euclidean decomposition $j = n \, \ell + r$, $0 \le r \le \ell-1$, 
and a direct iterative application of inequality \eqref{ineg:additiviteEntrop} together with \eqref{lem:defENtropB}  imply
\bean
H_j(\pi_j) &\ge& n \, H_\ell(\pi_\ell) + H_r(\pi_r)
\\
&\ge& n \, H_\ell(\pi_\ell)  + (j-1) \, [ (\log c_m)_-  - \MM_m(\pi)].
\eean
We deduce that for any $\ell \ge 1$ 
$$
\liminf_{j \to \infty} H(\pi_j) \ge \liminf_{j \to \infty} {n \over j} \, H_\ell(\pi_\ell)  = H(\pi_\ell),
$$
from which \eqref{eq:RRH'lim} follows. 

\smallskip
We conclude by establishing the affine property of $\HH'$. Let us consider $F, G \in \PPP_m(\PPP(E))$ and $\theta \in (0,1)$, and let us
assume that $H(F_j) < \infty$, $H(G_j) < \infty$ for any $j \ge 1$, the case when $H(F_j) = \infty$ or $H(G_j) = \infty$ being trivial. 
Using that  $s \mapsto \log s$ is an increasing function and that $s \mapsto s \, \log s$ is a convex function, we have 
\bean
H (\theta \, F_j + (1-\theta) \, G_j) 
&=& \frac1j \int_{E^j} (\theta \, F_j + (1-\theta) \, G_j)  \, \log (\theta \, F_j + (1-\theta) \, G_j) \\
&\ge& \frac1j \int_{E^j} \{ \theta \, F_j \, \log (\theta \, F_j ) + (1-\theta) \, G_j  \, \log ( (1-\theta) \, G_j) \} \\
&=& \theta \, H (F_j) + (1-\theta) \, H(G_j) + \frac1j [ \theta  \, \log \theta + (1-\theta)  \, \log (1-\theta)]  \\
&\ge& H( \theta \,  F_j +  (1-\theta) \, G_j) +  \frac1j [ \theta  \, \log \theta + (1-\theta)  \, \log (1-\theta)].
\eean
Passing to the limit $j\to\infty$ in the two preceding inequalities and using  \eqref{eq:RRH'lim}, we get 
$$
\HH' (\theta \, F + (1-\theta) \, G) \ge \theta \, \HH' (F) + (1-\theta) \, \HH' (G) \ge \HH' (\theta \, F + (1-\theta) \, G),
$$
which is nothing but the announced affine property. \qed

\medskip
We establish now in the following abstract lemma the last argument which allows us to prove
the first equality in \eqref{eq:Rob&Ruelle} and which will be useful  in the next section in order to 
get the same property for the similar functionals on $\PPP_m(\PPP(E))$ built starting from the Fisher information.  

\begin{lem}\label{lem:niv3abst} Consider a sequence  $(K_j)$ of functionals on $\PPP_m(E^j)$, $m \ge 0$,  such that 

\begin{itemize}

\item[(i)] $K_j : \PPP_m(E^j) \to \R \cup\{+\infty\}$ is convex, proper and l.s.c. for the weak convergence of measures on $\PPP_m(E^j)$ for any $j\ge1$.  Moreover, either $m=0$ and $K_j$ is positive for each $j$, or $m>0$ and there exists $k \in (0,m)$, a constant $C_k \in \R^+$ such that the functional  
$$
\PPP(E^j) \to \R \cup \{\infty\}, \quad G \mapsto  K_j(G) +  j [C_k + M_k(G)]
$$
is nonnegative and is l.s.c with respect to the weak convergence   in $\PPP(E)$.

\smallskip
\item[(ii)] $j^{-1} K_j (f^{\otimes j}) =  K_1 (f)$ for all $ f \in \PPP_m(E)$ and $j \ge 1$.

\smallskip
\item[(iii)]  $K_j (G) \ge K_\ell (G_\ell) + K_r(G_r)$ for any $G \in \PPP(E^j)$ and  any $\ell,r$ such that 
$j = \ell + r$. 
 
\smallskip
\item[(iv)]  The functional  $\KK' : \PPP_m(\PPP(E))  \to \R \cup\{+\infty\}$  defined for any $\pi \in \PPP_m(\PPP(E))$ by (this a part of the theorem that the $\sup$ equals the $\lim$)
$$
\KK'(\pi) := \sup_{j \ge 1} \frac1j \, K_j(\pi_j) = \lim_{j \to + \infty}  \frac1j \, K_j(\pi_j) ,
$$
where $\pi_j$ denotes j-th marginal defined thanks to Theorem~\ref{theo:H&S}, 
is affine in the following sense. For any probability measure  $\pi \in \PPP_m(\PPP(E))$
and any partition
partition of $\PPP_m(E)$ by some sets $\omega_i$, $1\le i \le M$, such that
$\omega_i$ is an open set in $E \backslash (\omega_1 \cup \ldots \cup
\omega_{i-1})$ for any $1\le i \le M-1$, $\omega_M = \PPP_m(E) \big\backslash (\omega_1 \cup \ldots \cup
\omega_{M-1})$  and $\pi(\omega_i) > 0$ for any 
$1\le i \le M$, defining 
$$
\alpha_i := \pi(\omega_i) \quad \hbox{and} \quad \gamma^i := {1 \over \alpha_i} \, {\bf 1}_{\omega_i} \, \pi \in \PPP_m(\PPP(E))
$$
so that 
$$
\pi = \alpha_1 \, \gamma^1 + ... + \alpha_M \, \gamma^M
\quad \hbox{and} \quad \alpha_1   + ... + \alpha_M =1,
$$
there  holds
$$
\KK'(\pi) = \alpha_1 \, \KK'(\gamma^1) + \ldots + \alpha_M \, \KK'(\gamma^M).
$$

\end{itemize}
Then under the above assumptions, for any $\pi \in \PPP_m(\PPP(E))$,  there holds 
$$
\KK'(\pi) = \KK(\pi) := \int_{\PPP(E)} K_1(\rho) \, \pi (d\rho).
$$
The functional $\KK
: \PPP_m(\PPP(E)) \to \R \cup \{\infty \}$ is affine, proper and  l.s.c. with respect to the weak convergence in $\PPP_m(\PPP(E))$.

Moreover, it satisfies the following $\Gamma$-l.s.c.~property. For any sequence $F^N$ of  $\PPP_{\!sym}(E^N)$ and $\pi \in \PPP(\PPP(E))$ such that $F^N \wto \pi$ weakly in  $\PPP_m(E^j)_{\forall j}$, then 
\beqn\label{ineq:KFNtoKKpi}
\KK(\pi) \le \liminf_{N \to \infty} \, K (F^N).
\eeqn
\end{lem}

\smallskip\noindent
{\sc Proof of Lemma~\ref{lem:niv3abst}. } We split the proof into five steps.  

 \medskip\noindent
{\sl Step 1. A fist inequality $\KK \ge \KK'$ } 
We skip the proof that the $\lim$ equals the $\sup$ in point $iv)$. This is a consequence of the hypothesis $iii)$ - and the bound by below in point $i)$ in the case $m>0$ - and has already been proved in the proof of Lemma~\ref{lem:RRH2} for the entropy.
We fix $\pi \in \PPP_m(\PPP(E))$. Thanks to assumptions (i) and (ii), we easily compute 
\bean
\KK(\pi) 
&=& \int_{\PPP(E)} \frac 1j K_j(   \rho^{\otimes j} )  \, \pi(d\rho)
\\
&\ge& \frac 1j K_j \Bigl( \int_{\PPP(E)}   \rho^{\otimes j}  \, \pi(d\rho) \Bigr) = \frac 1j K_j (\pi_j).
\eean
Taking the supremum over $j$ in this inequality,  we get a first inequality
$$
\KK (\pi) \ge  \sup_{j \ge 1} \frac 1j K_j(\pi_j) = \KK'(\pi).
$$

\smallskip\noindent
{\sl Step 2.  $\JJ$ is l.s.c.  on $\PPP_m(E)$ with respect to the $\WW_1$-metric.} 

 We consider the case when $m > 0$, and choose $k \in (0,m)$ such that (i) holds. We explain in the step $3'$ below the necessary adaptation to do in the case $m=0$.

For any $\delta > 0$, by compactness, we can find 
a family of finite cardinal $\NN$ of balls $\BB_i := \BB(\rho_i,\delta) =  \{ \rho \in \PPP_m(E); \, W_1(\rho,\rho_i) < \delta \}$, $\rho_i \in  \BB\PPP_{m,1/\delta}$, of radius $\delta$ so that 
$$
\BB\PPP_{m,1/\delta} \subset \bigcup_{i=1}^\NN \BB_i.
$$
We associate to that partition and "almost" partition of unity by
$$
\phi_i(\rho) := 2\, \Bigl[1 - \frac{W_1(\rho,\rho_i)}{2\delta} \Bigr]_+, \quad 
\theta_i (\rho) := \frac{\phi_i(\rho)}{\sum_{j=1}^\NN \phi_j(\rho) + \delta}. 
$$
Finally, we set for any $\rho \in \PPP_m(E)$ 
$$
J^\delta(\rho) := \sum_{i=1}^\NN \theta_i(\rho) \, J^\delta_i,
$$
where
$$
J^\delta_i := \inf_{\rho \in \BB(\rho_i,2\delta)} J(\rho)
\quad \text{and} \quad 
J(\rho) := 
K_1(\rho) +  C_k + M_k(\rho). 
$$
We claim that by construction the functional $J^\delta$ is Lipschitz with respect to the $W_1$ metric on $\PPP(E)$, and satisfies 
\beqn \label{encad:Jdelta}
\forall \, \rho \in \PPP_m(E), \quad 
\frac{\indiq_{\BB\PPP_{m,1/\delta}}(\rho)}{1+\delta} \inf_{\rho'  \in \BB(\rho,4\delta)} J(\rho')  \le J^\delta(\rho) \le J(\rho),
\eeqn
where $\indiq$ denote the indicator function. To obtain both inequalities, we introduce $I^\delta(\rho) := \{ 1\le i \le \NN, \; W_1(\rho, \rho_i^\delta) \le 2\delta \}$, and rewrite
$$
J^\delta(\rho) := \sum_{i \in I^\delta(\rho)} \theta_i(\rho) \, J^\delta_i.
$$
But for any $i$ such that  $W_1(\rho, \rho_i) \le 2\delta$ we have
$$
\inf_{\rho' \in \BB(\rho_i,4\delta)} J(\rho') \le J^\delta_i = \inf_{\rho' \in \BB(\rho_i,2\delta)} J(\rho') \le J(\rho).
$$
 The upper bound in \eqref{encad:Jdelta} follows form the second inequality (on the right). Since $J(\rho) \ge 0$ by hypothesis (i), the first above inequality implies that
 \begin{align*}
 J^\delta(\rho)  &\ge \Bigl(  \sum_{i \in I^\delta(\rho)} \theta_i(\rho)  \Bigr)  \inf_{\rho' \in \BB(\rho_i,4\delta)} J(\rho') 
 \ge  \frac{\sum_{j=1}^\NN \phi_j(\rho)}{\sum_{j=1}^\NN \phi_j(\rho) + \delta}
\;  \inf_{\rho' \in \BB(\rho_i,4\delta)} J_1(\rho').
 \end{align*}
The bound by below in \eqref{encad:Jdelta} then follows because any $\rho \in  \BB\PPP_{m,1/\delta}$ is at least in one of the $\BB_i$ for some $i$, and then 
 $\sum_{j=1}^\NN \phi_j(\rho) \ge \phi_i(\rho) \ge 1.$
 The inequalities \eqref{encad:Jdelta} and the hypothesis that $J$ is l.s.c. with respect to the weak convergence on $\PPP(E)$ implies that
 \beqn \label{eq:limJdelta}
 \forall \rho \in \PPP_m(E), \quad \lim_{\delta \to 0 } J^\delta(\rho) = J(\rho).
 \eeqn
We can now introduce the functionals $\JJ^\delta$ and $\JJ$ defined for all  $ \pi \in \PPP_m(\PPP(E))$ by
\begin{align*}
\JJ^\delta(\pi) &:=  \int_{\PPP_m(E)} J^\delta(\rho) \, \pi(d\rho) \\
\JJ(\pi) &:=  \int_{\PPP_m(E)} J(\rho) \, \pi(d\rho) = \KK(\pi) + C_k + \MM_k(\pi).
\end{align*}
Since $J^\delta$ is Lipschitz with respect to the $W_1$-metric, the Kantorovich-Rubinstein duality theorem \cite[Theorem 1.14]{VillaniTOT}
implies that the functionals $\JJ^\delta$ is continuous with respect to the $\WW_1$-metric. 
Moreover, the upper bound in \eqref{encad:Jdelta} implies that $\JJ^\delta(\pi) \le \JJ(\pi)$, for any $\pi \in \PPP_m(\PPP(E))$. Finally, an application of Fatou's Lemma together with  \eqref{eq:limJdelta} implies
\begin{align*}
\liminf_{\delta \to 0} \JJ^\delta(\pi) &= 
\liminf \int_{\PPP_m(E)} J^\delta(\rho) \, \pi(d\rho) \ge 
\int_{\PPP_m(E)}  \liminf  J^\delta(\rho) \, \pi(d\rho) \\
&\ge  
\int_{\PPP_m(E)}   J(\rho) \, \pi(d\rho) = 
\JJ(\pi).
\end{align*}
All in all, we get that
$$
\forall \pi \in \PPP_m(\PPP(E)), \quad \JJ(\pi) = \sup_{\delta >0 } \JJ^{\delta}(\pi),
$$
and that implies that $\JJ$ is l.s.c. with respect to the $\WW_1$-metric since the $\JJ^\delta$ are continuous with respect to that metric. 

\medskip\noindent
{\sl Step 2'. A necessary adaptation in the case $m=0$.}
In that case, things are in some sense simpler since the functional $K$ is already positive, so that we may try directly to apply Step $2$ with $J = K_1$. However, there is one difficulty : the compact sets $\BB\PPP_{m,1/\delta}$ does not covers $\PPP(E)$; even if we take their union for $\delta>0$ and $m>0$. 

However, we can still do a correct proof if we fix $\pi$  at the beginning. We then choose a increasing function 
$g : \R^+ \to \R^+$ such that
\beqn \label{def:GeneMoment}
 \lim_{v \to +\infty} g(v) = + \infty
\quad \text{and}\quad  M_g(\pi_1) := \int_E  g(\langle v \rangle) \, \pi_1(dv)  <  \infty.
\eeqn
Then we can restrict ourselves to the set  $\PPP_g := \{ \rho \in \PPP(E), \; M_g(\rho) < + \infty\}$, since 
the last hypothesis on $g$ implies that $\pi\bigl( \PPP_g(E) \bigr) = 1$. If we now replace in step $2$, the sets $\BB\PPP_{m,1/\delta}$ by the still compact sets
$$
\BB\PPP_{g,1/\delta} := \bigl\{ \rho, \; M_g(\rho) \le \delta^{-1} \bigr\},
$$
and follow the same strategy, we will conclude that $\KK(\pi) = \sup_{\delta >0} \KK^\delta(\pi)$ were the $\KK^\delta$ will be continuous with respect to the $\WW_1$-metric. It implies that $\KK$ is l.s.c. at $\pi$. Since $\pi$ is arbitrary,  $\KK$ is globally l.s.c.

\medskip\noindent
{\sl Step 3.  $\KK$ is l.s.c. with respect to the weak convergence of measures on $\PPP_m(\PPP(E))$. }

In the case $m=0$, that step is useless since in step $2'$ we proved that $\KK=\JJ$ is l.s.c.. So it remains only to treat the case $m>0$.
Since $\JJ = \KK + \MM_k + C_k $ is l.s.c. with respect to the $\WW_1$-metric on $\PPP_m(\PPP(E))$, the conclusion will follows if we show that $\MM_k$ is continuous with respect to the weak convergence on $\PPP_m(\PPP(E))$, defined in Definition \ref{def:CVG_PPP_k}. 

For this, we choose  $\rho,\mu \in \PPP_m(E)$. Since
$$
\forall \, v,v' \in E \qquad \bigl| \langle v\rangle^k - \langle v'\rangle^k  \bigr| \le k \, \min(1,|v-v'|)\, \bigl( \langle v\rangle^k + \langle v'\rangle^k \bigr),
$$
we obtain if we chose an optimal transference plan $\pi$ (for the distance $d_E$ on $E$) between $\rho$ and $\mu$
\begin{align*}
\bigl| M_k(\rho) - M_k(\mu)  \bigr| &\le \int \bigl| \langle v\rangle^k - \langle v'\rangle^k  \bigr| \, \pi(dv,dv') \\
& \le k \int  d_E(v,v') \bigl(  \langle v\rangle^k + \langle v'\rangle^k \bigr) \, \pi(dv,dv') \\
& \le  k  \left( \int  d_E(v,v')^{\frac m{m-k}} \, \pi(dv,dv') \right)^{1 - \frac km}
\left( M_m(\rho) + M_m(\mu) \right)^{ \frac km} ,
\end{align*}
so that 
\begin{align*}
\bigl| M_k(\rho) - M_k(\mu)  \bigr| & \le k \left( M_m(\rho) + M_m(\mu) \right)^{ \frac km} W_1(\rho, \mu)^{1 - \frac km},
\end{align*}
where we have used H\"older inequality and the fact that $d_E \le 1$. Choosing now two $\alpha,\beta \in \PPP_m(\PPP(E))$ and an 
 optimal transference plan
 $\pi$ (for the distance $W_1$ on $\PPP(E)$) between them, we get 
\begin{align*}
\bigl|  \MM_k(\alpha) - \MM_k(\beta) \bigr|  & =  
\biggl| \int  M_k(\rho) - M_k(\rho')  \, \pi(d\rho,d\rho') \biggr|
\\
& \le k \int \left( M_m(\rho) + M_m(\rho') \right)^{ \frac km} W_1(\rho, \rho')^{1 - \frac km} \pi(d\rho,d\rho'),
\end{align*}
and then 
\begin{align*}
\bigl|  \MM_k(\alpha) - \MM_k(\beta) \bigr| & \le k \left( \MM_m(\alpha) + M_m(\beta) \right)^{ \frac km} \WW_1(\alpha, \beta)^{1 - \frac km},
\end{align*}
where we have used H\"older inequality. This concludes the step since weak convergence on $\PPP_m(\PPP(E))$ exactly means that $\WW_1$ goes to zero and the moment of order $m$ are bounded. 

\medskip\noindent
{\sl Step 4. Proof of the remaining inequality $\KK' \ge  \KK$. } 
 Because $\PPP_m(E)$ endowed with the MKW distance $W_1$ is a
Polish space, for any fixed $\eps > 0$, we can cover it by a countable union of balls $\BB_n := \BB(f_n,\eps)$ of radius $\eps$. 
For a given $\pi \in \PPP_m(\PPP(E))$, we can choose $M$ such that 
$$\omega_M := \PPP_m(E) \backslash (\BB_1 \cup \ldots \cup \BB_{M-1}) \quad  \text{satisfies} \quad  \pi(\omega_M ) \le \eps$$ 
and denote $\omega_i := \BB_i \big\backslash (\BB_1 \cup \ldots \cup \BB_{i-1})$ for all $1\le i \le \NN-1$. 
We define then
$$
\alpha_i := \pi(\omega_i), \quad \gamma_i := \frac1{\alpha_i} \pi_{| {\omega_i} } , \quad 
\pi^M := \sum_{i=1}^M \alpha_i \,  \delta_{\gamma^i_1}, 
\qquad 
\gamma^i_1 = \int_{\PPP(E)} \rho \, \gamma^i(d\rho).
 $$
For any $1 \le i \le M$, we have
$$
\KK'(\gamma^i) := \sup_{j \ge 1} {1 \over j} \, K_j  (\gamma^i_j) \ge K_1 (\gamma^i_1).
$$
Using the affine property $(iv)$ of $\KK'$, the above inequality  
and the definitions of  $\pi^N$ and $\KK$, we get
\bear\nonumber
\KK' (\pi) 
&=& \alpha_1 \, \KK'(\gamma^1) + ... + \alpha_M \, \KK'(\gamma^M)
\\ \label{ineq:HH>HH'}
\KK' (\pi) 
 &\ge& \alpha_1 \, K_1(\gamma^1_1) + ... + \alpha_M \,K_1(\gamma^M_1)
 = \KK(\pi^M).
\eear

We observe that because $\pi^M_1 = \pi_1$, we have
\bean
\langle \pi^M_1,|v|^m \rangle
= \langle \pi_1,|v|^m \rangle
=  M_m(\pi) < \infty,
\eean
and in particular $\pi^M \in \PPP_m(\PPP(E))$. Moreover, defining $T^M : \PPP(E)
\to \{\gamma^1, ... , \gamma^M\}$ by $T^M(\rho) = \gamma^i$ for any $\rho \in
\omega_i$, we have $\pi^M = (T^M) _\sharp \pi$ and then 
$$
\WW_1(\pi,\pi^M) \le \bigl\langle (id \otimes T^M)_\sharp \pi, W_1(.,.) \bigr\rangle \le 2 \eps.
$$ 

\smallskip
We consider now a sequence $\eps\to0$ and the corresponding sequence $(\pi^M)$
for which we then have by construction $\pi^M \wto \pi$ weakly in
$\PPP_m(\PPP(E))$. Inequality \eqref{ineq:HH>HH'}, the above convergence and the
l.s.c. property of  $\KK$ proved in step $2$ and $3$ imply the second (and reverse) inequality 
$$
\KK(\pi)  \le  \liminf_{M \to \infty} \KK(\pi^M)  \le \KK' (\pi).
$$

\medskip\noindent
{\sl Step 5. The $\Gamma$-l.s.c. property of $\KK$. } We give the proof only in the case $m>0$, the case $m=0$ being simpler.  
We consider $(F^N)$ a sequence of $\PPP_{\! sym}(E^N)$ and $\pi \in \PPP(\PPP(E))$ such that $F^N \wto \pi$ weakly in $\PPP_m(E^j)_{\forall j}$, in particular  $M_m(F^N_1) \le a$ for some $a \in (0,\infty)$.  For any fixed $j \ge 1$, using the l.s.c. property of $K_j$, introducing the  Euclidean decomposition $N = n \, j + r$, $0 \le r \le j-1$ and using iteratively the inequality $(iii)$ of the hypothesis  as in the proof of Lemma~\ref{lem:RRH2} as well as the lower bound on $K_r$ provided by hypothesis (i), we get 
\bean
\frac1j K_j(\pi_j) 
&\le& \liminf_{N \to \infty}\frac1j K_j(F^N_j) 
\\
&\le&  \liminf_{N \to \infty} {1 \over n \, j}  \{ K_N(F^N) - r \, K_r(F^N_r) \} \\
&\le&  \liminf_{N \to \infty}  \{ {1 \over n \, j}  \, K_N(F^N) + {1 \over n} \, (C_m + a) \} \\
&=&  \liminf_{N \to \infty} N^{-1} K_N(F^N).
\eean
We deduce \eqref{ineq:HFNtoHHpi} thanks to \eqref{eq:Rob&Ruelle}. 
That concludes the proof. \qed

\medskip\noindent
{\sc Proof of Theorem~\ref{th:Rob&Ruelle}. }  The proof is just an application of the two previously proved lemmas.
First, let us observe that $H_j$, $\HH$ and $\HH'$ fulfill  the assumptions of Lemma~\ref{lem:niv3abst} since that 
$(i)$ is nothing but \eqref{eq:EntropBoltzHi&ii}, $(ii)$ is a consequence of Lemma~\ref{lem:defENtropB}, $(iii)$
is nothing but \eqref{ineg:additiviteEntrop}, 
and  a stronger version of $(iv)$  has been established in Lemma~\ref{lem:RRH2}. Then \eqref{eq:Rob&Ruelle} and \eqref{ineq:HFNtoHHpi} are exactly the conclusion of  Lemma~\ref{lem:niv3abst} applied to the entropy.
 \qed


\subsection{Level-3 Fisher information for mixtures} 

We state now a similar result for the Fisher information for mixtures of probability measures. 

 \smallskip
Let us assume that  $E = \R^d$ or  $E$ is an  open connected and bounded set of $\R^d$ with smooth boundary. 
Then, for any $\pi \in \PPP(\PPP(E))$ we define 
\beqn\label{eq:defII}
\II(\pi) := \int_{P(E)} I(\rho) \, \pi(d\rho), 
\eeqn
where $I$ is the Fisher information  defined on $\PPP(E)$.

\begin{theo} \label{th:Fisher3}  
${\bf (1)}$ The functional $\II: \PPP(\PPP(E)) \to \R \cup \{\infty \}$ is
affine, nonnegative and l.s.c. for the weak convergence. 
Moreover, for any $\pi \in \PPP(\PPP(E))$,    there holds
\beqn\label{eq:Rob&RuelleFisher}
\II(\pi)  = \sup_{j \in \N^*} I(\pi_j) 
= \lim_{j\to\infty} I(\pi_j),
\eeqn
where $I$ stands for the normalized Fisher information defined in $\PPP(E^j)$ for any $j\ge1$. 

\smallskip
${\bf (2)}$
Consider $(F^N)$ a sequence of  $\PPP_{\!sym}(E^N)$ 
and $\pi \in \PPP(\PPP(E))$ such that $F^N \wto \pi$ weakly in  $\PPP(E^j)_{\forall j}$. Then 
\beqn\label{ineq:IFNtoIIpi}
\II(\pi) \le \liminf I(F^N).
\eeqn
\end{theo}

As for Theorem~\ref{th:Rob&Ruelle}, the proof of Theorem~\ref{th:Fisher3} relies on the abstract lemma \ref{lem:niv3abst}.  The hypothesis of that lemma are proved to be true in the lemma~\ref{lem:RRI2} below.  Two useful intermediate results are stated in the next two lemmas.

\begin{lem} \label{lem:rhot} There exist :

- a family of regularizing operators $S_t :\PPP(E) \to \PPP(E)$ defined for any $t>0$,  

- a family $(C_t)$ of positive constants 

- a family $\eps_t$ of positive constants such that $\eps_t \to 0$ when $t \to 0$ 

- for any $k>0$, a family $(\eps'_{kt})$ of positive constants so that $\eps'_{kt} \to 0$ when $t \to 0$ 

such that for any $\rho \in \PPP(E)$ and any $t > 0$, denoting $\rho_t := S_t(\rho)$ we have
 \begin{align}\label{eq:rhot}
 &I(\rho_t) \le I(\rho),
\quad M_k(\rho_t) \le 2^k \,\bigl( M_k(\rho) + \eps_{kt}\bigr), \quad 
\|\nabla \ln \rho_t \|_\infty \leq {C_t} \\
& \text{and}\quad W_1(\rho,\rho_t) \le \eps_t . \nonumber
\end{align}

\end{lem}

\smallskip\noindent
{\sc Proof of Lemma~\ref{lem:rhot}. } We only consider the case $E=\R^d$. The case when $E$ is a  smooth bounded open set can be handled 
similarly by using for $(\rho_t)$ the  solution of the heat equation (with Neumann boundary conditions) and  the strong maximum principle. 
We define   
$$
\eta_t(z) := \frac {C_d} {t^d} e^{- <z/t>} = \frac {C_d} {t^d} e^{- \sqrt{1+ (|z|/t)^2}} \quad\hbox{and}\quad \rho_t := \eta_t * \rho.
$$
Observing that  
$$
 \frac{|\nabla \eta_t(z)|}{\eta_t(z)} = \frac1t \, \frac{|z|}{\langle z \rangle} \le \frac1t,
$$
we deduce that  for any $x \in \R^d$, we have
\bear\nonumber
|\nabla \rho_t(x) | & \le & \frac1t  \int_{ \R^d}    \eta_t(x-y) \rho(y) \,dy  = \frac1t \, \rho_t. 
\eear

The inequality on the moment of order $k$ is a consequence of the inequality 
$$
\langle x+y \rangle^k \le 2^k \bigl( \langle x \rangle^k + \langle y \rangle^k  \bigr),
$$ 
which leads to the claimed inequality with $\eps_{kt} =  M_k(\eta_t)= t^k M_k(\eta_1)$.

As $\rho_t$ is also an average of translations of $\rho$ (which has the same Fisher information as $\rho$), the convexity of the Fisher information implies that 
$$I(\rho_t) = I \left( \int \rho(\cdot - z) \, \eta_t(dz)\right) \le \int  I\bigl( \rho(\cdot - z) \bigr) \, \eta_t(dz) = I(\rho).
$$
We finally observe that  for any $\rho \in \PPP(E)$ there holds
$$
W_1(\rho, \rho_t) = W_1(\rho, \rho \ast \eta_t) \le \int_{\R^d} |z| \, \eta_t(z) \, dz = C_d \, t,
$$
and that proves the last estimate.
\qed

\begin{lem}\label{lem:pitj} Consider $\pi \in \PPP(\PPP(E))$ and  define the regularized family $\pi_t \in \PPP(\PPP(E))$, for $t > 0$,   by  push-forward by $S_t$, $\pi_t := {S_t}_\# \pi$ or equivalently
$$
\langle \pi_t,\Phi \rangle = \langle \pi,\Phi_t \rangle \quad \forall \, \Phi \in C_b(\PPP(E))
$$
where $ \Phi_t \in C_b(\PPP(E))$ is defined by $\Phi_t(\rho) := \Phi(\rho_t)$ and $\rho_t$  is the defined in Lemma~\ref{lem:rhot}.
Also denote by $\pi_{tj} \in \PPP(E^j)$ the  j-th marginal of $\pi_t$ defined thanks to Theorem~\ref{theo:H&S}. 
For any $t >0$ and any $X^j := (x_1, ..., x_j) \in E^j$ there holds
\beqn\label{eq:pitj}
\bigl | \nabla_1 \ln \pi_{tj}(X^j) \bigr| \le C_t.
\eeqn
\end{lem}

\smallskip\noindent
{\sc Proof of Lemma~\ref{lem:pitj}. } Thanks to Lemma~\ref{lem:rhot}, we write  
\bean
\frac{\bigl | \nabla_1 \pi_{tj}(X^j) \bigr|}{\pi_{tj}(X^j)} & =& \frac{ \biggl | \int \nabla_1 \rho_t(x_1) \rho_t^{\otimes  j-1}(x_2,\ldots x_j)\pi(d\rho) \biggr |}{\pi_{tj}(X^j)} \\
& \leq & C_t  \frac{ \int  \rho_t^{\otimes  j}(X_j)\pi(d\rho)}{\pi_{tj}(X^j) }  = C_t,
\eean
which is nothing but \eqref{eq:pitj}.  
\qed

\begin{lem}\label{lem:RRI2}   For any  $\pi \in \PPP(\PPP(E))$ we define
$$
\II'(\pi)  := \sup_{j \in \N^*} I(\pi_j). 
$$
The functional $\II'
: \PPP(\PPP(E)) \to \R \cup \{\infty \}$ is nonnegative, l.s.c. for the weak convergence, satisfies 
\beqn\label{eq:RRIlim}
\II'(\pi) =  \lim_{j\to\infty} I(\pi_j)
\eeqn
and is affine in the same sense as formulated in point (iv) of
Lemma~\ref{lem:niv3abst}. \end{lem}

\smallskip\noindent
{\sc Proof of Lemma~\ref{lem:RRI2}. } The fact that $\II'$ is nonnegative and l.s.c. is clear and 
\eqref{eq:RRIlim}  comes from the monotony property 
$I(\pi_{j-1}) \le I(\pi_{j})$, $\forall \, j \ge 2$ established  in Lemma~\ref{lem:FisherCondiv} $(i)$. 
It remains only  to prove the linearity property of $\II'$. 
For the sake of simplicity we only consider the case when $M=2$ and $\omega_1$
is a ball. The case when 
$\omega_1$ is a general open set can be handled in a similar way and the case
when $M \ge 3$ can be deduced 
by an iterative argument. For some given $\pi \in \PPP(\PPP(E))$ which is not a Dirac mass, $f_1 \in \PPP(E)$ and $r \in (0,\infty)$
so that 
$$
\theta := \pi(\BB_r) \in (0,1), \quad \BB_r := \BB(f_1,r) = \{ \rho, W_1(\rho,f_1) < r\},
$$
we define 
$$
F := \frac1 \theta \indiq_{\BB_r} \pi, \quad G := \frac1{1-\theta} \indiq_{\BB_r^c} \pi
$$
so that 
$$
F,G \in \PPP(\PPP(E)) 
\quad\hbox{and} \quad \pi = \theta F + (1- \theta )G,
$$
and we have to prove that 
\beqn\label{eq:linearIIprim}
\II'(\pi) = \theta \, \II'(F) + (1-\theta) \, \II'(G).
\eeqn
We split the proof of that claim in four steps. 

\smallskip\noindent
{\sl Step 1. Approximation and estimation of the affinity defect. } As explained  for $\pi$ in the statement of  Lemma~\ref{lem:pitj}, we define $F_t$ and $G_t$ to be the push-forward of the 
measures $F$ and $G$ by the regularisation operator $S_t$, and then $F_{tj}$ and
$G_{tj}$ are their projections on $\PPP(E^j)$
$$
F_{tj} := \int_{\PPP(E)} \rho^{\otimes j} F_t(d\rho) = \int_{\PPP(E)} \rho_t^{\otimes j} F(d\rho),
\quad \text{or } \quad 
\langle F_{tj},\varphi \rangle = \int_{\PPP(E)} R_\varphi (\rho) F_{t}(d\rho), 
$$
via duality, for any $\varphi \in C_b(E^j)$ where $R_\varphi$ is the polynomial on $\PPP(E)$ associated to $\varphi$ thanks to \eqref{def:polyproba}.  The same holds for $G$.
We also remark that these two above operations (regularisation and projection on $E^j$) commute if we define the regularisation operators $S_t$ on $E^j$ by the convolution with $\eta_t^{\otimes j}$.  
It is worth emphasizing that we do not need here, in order  to define these objects, that
$F$ and $G$ are probability measures, but only that they are Radon measures on $\PPP(E)$.  

\smallskip
For any given $j \in \N$, we define
\begin{align*}
A_{tj} &:=  \theta \, I (F_{tj}) + (1-\theta) \, I(G_{tj}) 
- I(\theta \, F_{tj} + (1-\theta) \, G_{tj}), \\
& = \theta \int \frac{|\nabla F_{tj}|^2}{F_{tj}} + (1-\theta) \int \frac{|\nabla G_{tj}|^2}{G_{tj}} - 
 \int \frac{|(1-\theta) \nabla G_{tj} + \theta \nabla F_{tj}|^2}{(1-\theta)G_{tj} + \theta F_{tj}}.
\end{align*}
After reduction to the same denominator, and some simplification, we end up with
\bean
A_{tj} &=&  \, \theta (1-\theta) \int \frac{G_{tj} F_{tj}}{(1-\theta)G_{tj} + \theta F_{tj}}
\biggl| \nabla_1 \ln \frac{G_{tj}}{F_{tj}}\biggr|^2. \\
&\le &
 2 \, \theta (1-\theta) \int \frac{G_{tj} F_{tj}}{(1-\theta)G_{tj} + \theta F_{tj}} \left(\bigl| \nabla_1 \ln F_{tj} \bigr|^2+ \bigl| \nabla_1 \ln G_{tj}\bigr|^2 \right).
\eean
We can estimate the r.h.s.  term thanks to Lemma~\ref{lem:pitj} by
$$
A_{tj} 
\le
4 \theta (1-\theta) C_t \int \frac{G_{tj} F_{tj}}{(1-\theta)G_{tj} + \theta F_{tj}}   .
$$
\smallskip\noindent
{\sl Step 2. Disjunction of the supports.} Let us introduce for any $s \in (0,r)$ the two measures on $\PPP(E)$ (which 
are not necessarily probability measures)
$$
F' :=  \indiq_{\BB_s} F  =\frac1 {\theta} \indiq_{\BB_s} \pi,\qquad 
F'' := \indiq_{\BB_r\backslash \BB_s} F, \quad \text{so that} \quad F' + F'' = F
$$
and let us observe that 
$$
\lim_{s \to r } \int F''(d\rho) = \lim_{s \to r } \int \indiq_{\BB_r\backslash \BB_s}(\rho)  F(d\rho) = 0  ,
$$
by Lebesgue's dominated convergence theorem.  
For any $t> 0$ and $j \ge 1$  there holds $F_{tj} '+ F_{tj}''=F_{tj}$
with $F''_{tj} \ge 0$, so that we may write for any $\eps > 0$ 
\bean
A_{tj} &\le & 
4\theta (1-\theta) C_t \int \frac{G_{tj} F'_{tj}}{(1-\theta)G_{tj} + \theta F'_{tj}}   +  
4\theta  C_t \int  F''_{tj},\\
&\le & 
4  \theta (1-\theta) C_t \int \frac{G_{tj} F'_{tj}}{(1-\theta)G_{tj} + \theta F'_{tj}}   +  \eps,
\eean
taking $s$ close enough to $r$, and this independently of $j$ and $t$ because 
$$\int_{E^j} F''_{tj} = \int_{\PPP(E)} F''_t = \int_{\PPP(E)} F.$$ 

\smallskip\noindent
{\sl Step 3. Concentration.} We introduce the real numbers $u= \frac{r+s}2$ and  $\delta = \frac{r-s}2$, 
depending on $\eps$, as well as the set 
$$
\tilde B_{u} := \{ X^j=(x_1,\ldots,x_{j}) \;, \;  W_1(\mu^j_{X^j}, f_1) < u \}
\subset E^j
$$
which is nothing but the reciprocal image of the ball $\BB_u \subset \PPP(E)$ by the empirical measure map. 
Using that  
$$
\frac{G_{tj} F'_{tj} }{(1-\theta)G_{tj}
+ \theta F'_{tj} } \le \frac1\theta G_{tj} \, {\bf 1}_{\tilde B_{u}} + 
\frac{1}{ 1-\theta} F'_{tj}  \, {\bf1}_{\tilde B_{u}^c}, 
$$
we get 
\beqn\label{eq:Atj}
A_{tj} \leq 
4C_t \left( 
 (1-\theta) \int_{\tilde B_{u}}  G_{tj} +  \theta \int_{\tilde B_{u}^c}  F'_{tj}  \right)
 +  \eps.
\eeqn
 If $\rho$ belongs to the support of $F'$ and  $X^j \in \tilde
B_{u}^c$, we have thanks to the last estimate in Lemma~\ref{lem:pitj}
\bean
W_1(\mu^j_{X^j},\rho_t ) 
&\ge&
 W_1(\mu^j_{X^j},f_1) - W_1(f_1,\rho) - W_1(\rho,\rho_t) 
 \\
 &\ge&  u - s  - C_d \, t \ge \delta/2,
\eean
for any $t \in [0,T(\eps)]$, $T(\eps) > 0$.  We first assume that $\pi \in \PPP_m(\PPP(E))$ for some $m>0$, which implies also that $F,G \in \PPP_m(\PPP(E))$.
Gathering this information with  the Chebychev inequality, estimate~\eqref{rem:Omega0FNf}
and estimate~\eqref{eq:rhot}, 
 we conclude that  
 \bear\nonumber
 \int_{\tilde B_{u}^c}  F'_{tj}   
&= & 
 \int_{\PPP(E)} \langle \rho_t^{\otimes j}, {\bf 1}_{\tilde B_{u}^c} \rangle
F'(d\rho)
 \\ \nonumber
& \le & \frac2 \delta 
\int_{\PPP(E)}  \left( \int_{E^j} W_1(\mu_{X^j}^j,\rho_t)  \rho_t^{\otimes j}(dX^j)
\right) F'(d\rho) 
\\ \label{eq:BucF'}
& \le & \frac {C}  {\delta j^\gamma}
\int_{\PPP(E)} M_m(\rho_t)^{1/m} F'(d\rho) 
\le 
\frac C{\delta j^\gamma} \,  \Bigl( \MM_m(F) + \eps_{tm} \Bigr)^{1/m},
\eear
with $\gamma := 1/(d+2+d/m)$.   With exactly the same arguments, we prove that for any $\eps > 0$ and any $t \in [0,T(\eps)]$
\beqn\label{eq:BuG}
 \int_{\tilde B_u } G_{tj}
\le \frac{2C}{\delta j^\gamma} \,   \MM_m(G)^{1/m} + \eps.
\eeqn

 Gathering \eqref{eq:Atj} with \eqref{eq:BucF'} and \eqref{eq:BuG}, we get 
that for any $\eps > 0$, $t \in (0,T(\eps)]$ and $j \ge 1$, 
$$
A_{tj}   \le 
\frac{ 4C_t M_m(\pi)^{1/m} }{\delta j^\gamma} +  3\eps,
$$
and then for any $\eps > 0$, $t \in (0,T(\eps)]$ 
\beqn\label{eq:limAtj}
\limsup_{j \to \infty} A_{tj}   \le  3\eps.
\eeqn

 \smallskip\noindent
{\sl Step 3'. Adaptation for $\pi \notin \PPP_m(\PPP(E))$. }
In the case when $\pi \notin \PPP_m(\PPP(E))$ whatever is $m>0$, we can still prove that
\eqref{eq:limAtj} holds true for any $\eps$ and $t \in (0,T(\eps))$ where $T(\eps)$ is small enough.   Remark that it cannot be the case if $E$ is a smooth bounded open set, so we have only to deal with the case $E=\R^d$ here.

The idea is the same as in the proof of Lemma~\ref{lem:niv3abst}. We choose a function $g:\R_+ \to \R_+$  satisfying \eqref{def:GeneMoment} together with $g(2x) \le 2 g(x)$ for all $x \ge 0$. We argue by using moment with respect to $g(\langle \cdot \rangle )$ rather than to $\langle \cdot \rangle^m$. The property $g(2x) \le 2g(x)$ ensures that the estimate on the moments in Lemma~\ref{lem:rhot} is still true with the moment  $M_g$. 

Next, for any $R >0$, we introduce the mapping from $P_R :\R^d \to \R^d$ defined by 
$$
P_R(x) := \begin{cases}
x & \text{if } |x| \le R \\
R \frac x{|x|} & \text{else}
\end{cases}.
$$
using the concentration estimatge~\eqref{rem:Omega0FNf} for the probability $(\rho \circ P_R^{-1})^{\otimes N} = ({P_R}_\# \rho)^{\otimes N}$,  it is still possible to deduce that
$$
\int_{E^j} W_1(\mu_{X^j}^j,\rho_t)  \rho_t^{\otimes j}(dX^j) \leq \frac C{\delta j^\gamma}R + \frac{2M_g(\rho_t)}{g(R)}.
$$
Summing up with respect to $\pi$, choosing $R$ large enough and letting $j \to + \infty$, we get the claimed inequality for the limsup (with maybe a $4\eps$ in place of the $3\eps$).

\smallskip\noindent
{\sl Step 4. Conclusion. }
The regularization by convolution (or with the heat flow) implies that for any  $\alpha \in \PPP(\PPP(E))$
$$
\II(\alpha_t) = \sup_{j \le 1} I(\alpha_{tj}) = \sup_{j \le 1} I(\alpha_{jt}) \le
\sup_{j \le 1} I(\alpha_j ) = \II(\alpha).
$$
Moreover,  the last point in Lemma~\ref{lem:rhot}, implies that 
 $\alpha_t \wto \alpha$ with respect to the $\WW_1$-metric. Thanks to the previous inequality and the l.s.c. property of  $\II'$, we obtain
\beqn \label{eq:CVGIIt}
\II'(\alpha) = \lim_{t \to 0} \II'(\alpha_t). 
\eeqn
Turning back to the definition of $A_{tj}$, the estimate \eqref{eq:limAtj} and the above properties, we obtain
for any $\eps > 0$, $t \in (0,T(\eps)]$ 
\bean
\II'(\pi) &\ge& \II'(\pi_t) 
\\
&\ge& \theta \, \II'(F_t) + (1-\theta) \, \II'(G_t) -  3\eps.
\eean
First passing to the limit $t \to 0$ and using \eqref{eq:CVGIIt} we get 
\bean
\II'(\pi) \ge  \theta \, \II'(F) + (1-\theta) \, \II'(G) -  3\eps,
\eean
for any $\eps > 0$, which concludes the proof of 
\eqref{eq:linearIIprim} since the reverse inequality is just
a consequence of the convexity of the functional $\II'$. 
 \qed

\medskip
\noindent{\sc Proof of Theorem~\ref{th:Fisher3}.  }
We only have to observe that $I_j$, $\II$ and $\II'$ fulfil  the assumptions of Lemma~\ref{lem:niv3abst}. But the assumption $(i)$ is a  consequence of Lemma~\ref{lem:FisherCondii},  the assumption $(ii)$ is a consequence of
 Lemma~\ref{lem:Fisher1},   the assumption $(iii)$  is proved in Lemma~\ref{lem:FisherCondii} and 
assumption $(iv)$  in Lemma~\ref{lem:RRI2}. 
Then  \eqref{eq:Rob&RuelleFisher}  and \eqref{ineq:IFNtoIIpi} are exactly the conclusion of  Lemma~\ref{lem:niv3abst}
adapted to the Fisher information.
\qed

\begin{prop} \label{prop:suppLpPPPE}
Consider $\pi \in \PPP(\PPP(E))$ and $(\pi_j)$ the associated family of compatible and symmetric probability measures in $\PPP(E^j)$ defined
as in the De Finetti, Hewitt \& Savage theorem. For any $p \in[1,+\infty]$, the following equality holds
\beqn \label{eq:suppLpPPPE}
\pi-\text{Suppess } \{  \| \rho \|_p, \; \rho \in \PPP(E) \} = \sup_{j \in \N} \| \pi_j \|_p^{\frac1j} = \lim_{j \to +\infty} \| \pi_j \|_p^{\frac1j}. 
\eeqn
It is part of the result that the limit exists. In particular, it implies the equivalence
$$ 
\forall j \in \N, \; \| \pi_j \|_{L^p (E^j)} \le C^j \Longleftrightarrow  \pi-\text{Suppess } \{  \| \rho \|_p, \; \rho \in \PPP(E) \} \leq C. 
$$
\end{prop}

\noindent
{\sc Proof of Proposition~\ref{prop:suppLpPPPE}. }  First remark that there is nothing to prove for $p=1$ since we are dealing with probability measures.  Now, one inequality is a simple consequence of the De Finetti, Hewitt \& Savage theorem. In fact, using the definition of $\pi_j$, we get
$$
\| \pi_j \|_p = \left\|  \int_{\PPP(E)} \rho^{\otimes j}\,\pi(d\rho) \right\|_p 
\leq   \int_{\PPP(E)} \| \rho^{\otimes j} \|_p\,\pi(d\rho)   
=  \int_{\PPP(E)} \| \rho \|_p^j  \,\pi(d\rho) ,
$$
and the last quantity is clearly bounded by $ M^j$, $M: = \pi-\text{Suppess } \{  \| \rho \|_p, \; \rho \in \PPP(E) \}$. 

For the reverse inequality, we denote by $q \in (1,+\infty]$ the real conjugate
to $p$.
Because $L^q(E) = (L^p(E))'$, the Hahn-Banach separation theorem infers that for any $\lambda <  M$ there exists $f$ in the unit ball of $L^q(E)$ so that the set
$$
\mathcal B :=  \{  \rho \in \PPP(E) \text{ s.t. } \; \int f(x) \rho(dx) \geq \lambda  \}
$$
is of $\pi$-measure positive : $ \delta := \int_{\mathcal B} \pi(d\rho) >0$. Now for any $j \in \N$
$$
\| \pi_j \|_p \geq \int_{E^j} f^{\otimes j} \;d\pi_j = \int_{\PPP(E)} \left( \int_{E^j} f^{\otimes j} \rho^{\otimes j} \right) d\pi(\rho) \geq \delta \lambda^j,
$$
which implies the reserve inequality $M \leq \lim_{j \to +\infty} \| \pi_j \|_p^{\frac1j}$.
\qed


 \subsection{Strong version of  De Finetti, Hewitt and Savage theorem  and strong convergence in $\PPP(E^N)$}
\label{subsec:EntropieCombFiniEtatsPurs}

We begin that section by an HWI inequality valid on $\PPP(\PPP(E))$, which is just a "summation" of the usual one and will be very useful in the sequel. 

\begin{prop} \label{prop:HWI_PPE} Assume $E=\R^d$  or more generally that \eqref{eq:GalHWI}  holds for $N=1$.
For any $\alpha,\beta \in \PPP(\PPP_2(E))$, we have
\beqn\label{ineq:HHWWII}
\HH(\alpha) \leq \HH(\beta) + C_E \sqrt{ \II(\alpha)} \, \WW_2(\alpha,\beta). 
\eeqn
As a consequence, the entropy  $\HH$ is continuous on bounded sets
relatively to $\II$. In more precise words, 
 if  $(\pi_n)$ is a bounded sequence of $\PPP_m(\PPP(E))$, $m > 0$,  such that 
$$
\pi_n \wto \pi \,\,\hbox{ weakly in } \,\, \PPP(\PPP(E)) \quad \hbox{and} \quad
\II(\pi_n) \le C,
$$
then $\HH(\pi_n) \to \HH(\pi)$.
\end{prop}

\noindent
{\sc Proof of proposition~\ref{prop:HWI_PPE}. }
A first way in order to prove \eqref{ineq:HHWWII} is just to pass in the limit in the HWI inequality \eqref{eq:GalHWI} for $\alpha_N$ and $\beta_N$ and use the inequality stated 
in lemma~\ref{lem:WassNleqW1} for the quadratic cost, and the result of the previous section about level 3 entropy and Fisher information \ref{eq:Rob&Ruelle} et \ref{eq:Rob&RuelleFisher}.

\smallskip
Another possibility is to sum up the HWI inequality \eqref{ineq-HWI} for $\rho \in \PPP(E)$. Choosing an optimal transference plan $\Pi$ for $\WW_2$ between $\alpha$ and $\beta$, we have 
\bean
\int_{\PPP(E)} H(\rho) \, \Pi(d\rho,d\eta) & \leq & \int_{\PPP(E)} H(\eta)  \, \Pi(d\rho,d\eta)+ \int_{\PPP(E)} \sqrt{I(\rho)} W_2(\rho,\eta)  \, \Pi(d\rho,d\eta) ,
\eean
so that 
\bean
H(\alpha) & \leq & H(\beta) + \left(  \int_{\PPP(E)} I(\rho)  \, \Pi(d\rho,d\eta) \right)^{\frac12}
\left(  \int_{\PPP(E)}  W_2(\rho,\eta)^2  \, \Pi(d\rho,d\eta) \right)^{\frac12},
\eean
thanks to Cauchy-Schwarz inequality. It leads to the desired inequality. 

The second point is obtained by two applications of the previous
inequality, leading to
$$
\left| \HH(\pi_n) - \HH(\pi) \right| \leq \left( \sqrt{\II(\pi)} +
\sqrt{\II(\pi_n)}\right) \WW_2(\pi_n,\pi),
$$
and then using the l.s.c. property of the level 3 Fisher information in order to prove that $\II(\pi) < \infty$. We conclude
by remarking that the RHS converges to $0$ as $n$ tends to $\infty$.
\qed

\medskip
The  results of the preceding section and  the HWI inequality  make possible to compare different senses of convergence for sequences of $\PPP(E^N)$, $N \to \infty$, without any assumption of chaos.

\begin{theo}\label{theo:FisherVsEntropGal} Assume $E=\R^d$ or  $E \subset \R^d$ is a bounded connected open subset with 
smooth boundary and  that \eqref{eq:GalHWI} 
holds.
Consider $(F^N)$ a sequence of $\PPP_{\! sym}(E^N)$ and $\pi \in \PPP (\PPP(E))$ such that 
$F^N \wto \pi$ weakly in  $\PPP_k(E^j)_{\forall j}$, $k  >2$. 

\smallskip\noindent
{\bf (1) } In the list of assertions below, each assertion implies the one which follows:

(i) $I(F^N) \to \II(\pi)$, $\II(\pi) < \infty$;

(ii) $I(F^N)$ is bounded; 

(iii) $H(F^N) \to \HH(\pi)$, $\HH(\pi) < \infty$.

\smallskip

\noindent
{\bf (2) } More precisely, the following version of the implication (ii)
$\Rightarrow$ (iii) holds. There exists a 
numerical constant $C$ such that for any $k>2$ and $K>0$, and  for any any sequence $(F^N)$ of
$\PPP_{\! sym}(E^N)$ satisfying
$$
\forall \, N \qquad M_k(F^N_1) \le K^k, \quad I(F^N) \le K^2,
$$
there holds 
\beqn\label{estim:FisherVsEntropGal}
\forall \, N \ge 4^{2d} \qquad 
|H(F^N) - \HH(\pi) | \le  K \,   W_2(F^N,\pi_N) +  C K^{d'} \frac{\ln(K N) }{N^\gamma}, 
\eeqn
with $\gamma :=\frac{k-2}{k(1+2d) + 4d - 2}$ and as usual $d' =\max(2,d)$.

\smallskip\noindent
{\bf (3) } In particular, for any
sequence $(\pi_j)$ of symmetric and compatible probability measures of
$\PPP(E^j)$ satisfying 
$$
M_k(\pi_1) \le K^k, \qquad \forall \, j \ge 1 \quad   I(\pi_j) \le K^2,
$$
there holds 
\beqn\label{estim:FisherVsEntropGal2}
\forall j \ge 4^{2d}, \qquad  |H(\pi_j) - \HH(\pi) | \le C  K^{d'} \frac{\ln(K j) }{j^\gamma}
\eeqn
for the same value of $\gamma$. 
In other words, \eqref{estim:FisherVsEntropGal2} gives a rate of convergence for the limit \eqref{eq:Rob&Ruelle}. 
\end{theo}

The fact that the constant $C$ does not depend on $k$ is interesting when the space $E$ is compact or the measures $F^N$ have strong integrability properties, for instance an exponential moment. It allows to choose large $k$ and get almost the largest exponent $\gamma$ possible. Precise versions of the point (iii) are stated (without proofs) in the corollary below.

\begin{cor} 

(i) In the case where $E$ is compact, we denote $K := \max( \mathrm{diam}(E), \sqrt{\II(\pi)})$. Then there holds  for all $j \ge 4^{2d}$ 
\beqn
 |H(\pi_j) - \HH(\pi) | \le C  K^{d'} \frac{\ln(K j) }{j^\gamma} \quad \text{with} \; \gamma = \frac{1}{2d+1}
\eeqn

(ii) If $M_{\beta,\lambda} (\pi_1) := \int_{E} e^{\lambda |x|^\beta} \pi_1(dx) <+\infty$ for some $\lambda >0$ and $\II(\pi) < +\infty$, there exists a constant $C(d,\beta,\lambda,\II(\pi))$ such that for $j$ large enough ( $\geq C' \ln M_{\beta,\lambda} (\pi_1)$)
\beqn
|H(\pi_j) - \HH(\pi) | \le C \frac{[\ln j]^{1+d'/\beta} }{j^\gamma} \quad \text{with} \; \gamma = \frac{1}{2d+1}. 
\eeqn 
\end{cor}

\noindent
{\sc Proof of Theorem~\ref{theo:FisherVsEntropGal}. }  We split the proof into
four steps. 

\smallskip\noindent{\sl Step 1. } i)
implies ii) is clear. For ii) implies iii), we use the HWI
inequality~\eqref{eq:GalHWI} and we write
\bean
|H(F^N) - \HH(\pi) | 
&=& |H(F^N) - H(\pi_N) + H(\pi_N) - \HH(\pi) |
\\
&\le& C_E \Bigl( \sqrt{I(F^N)}+ \sqrt{I(\pi_N)} \Bigr) \, W_2(F^N,\pi_N) + |
H(\pi_N) - \HH(\pi) |.
\eean
We know from \eqref{eq:Rob&Ruelle} that $\HH(\pi) = \lim H(\pi_N)$  and from  \eqref{eq:Rob&RuelleFisher}
and \eqref{ineq:IFNtoIIpi} that \\ \mbox{$I(\pi_N) \le \II(\pi) \le \liminf
I(F^N) \le K$}, from which we conclude that 
there exist a sequence  $\eps_\pi(N) \to 0$ such that 
$$
|H(F^N) - \HH(\pi) |  \le  2\,C_E \,  K  \, W_2(F^N,\pi_N) + \eps(N).
$$
We now aim to estimate $\eps(N)$ more explicitly as claimed in point {\bf
(3)}. Then  {\bf (2)} will be a direct consequence of {\bf (3)} and the above estimate. 
 
 \medskip
 From now on, we only consider the case $E= \R^d$ since the general case is similar (and the case when $E$ is compact is
 even simpler). 
 
\medskip\noindent
{\sl Step 2. } 
   From \cite[Theorem A.1]{BGV} we know that for any $R,\delta > 0$ we may cover $\PPP(B_R)$ by $\NN(R,\delta/2)$ balls of
radius $\delta/2$ in $W_1$ distance (which is less accurate than the one considered in the above quoted result) with 
$$
\NN(R,\delta) \le  \Bigl( { C'_1 \, R \over \delta} \Bigr)^{C'_2 (R/\delta)^d},
$$
where the constant $C_1'$ and $C_2'$ are numerical. Let us fix $a \ge 1$ and recall that we define  $\BB\PPP_{k,a}(E) := \{ \rho \in \PPP(E) \text{ s.t. }
M_k(\rho) \leq a\}$. Next, for any $\rho \in \BB\PPP_{k,a}(E)$, we define $\rho_R
\in \PPP(B_R)$ by $\rho_R = \rho(B_R)^{-1} \, \rho \, {\bf 1}_{B_R}$ for $R$
large enough (so that it defines a probability measure),
and we observe that for any $f \in \PPP(E)$ we have 
$$
W_1(\rho,f) 
\le W_1(\rho_R,f)+ W_1(\rho_R,\rho),
$$
and that for any $R$ such that $R^k > 2a$
\bean
W_1(\rho_R,\rho) &\le&  \|\rho_R - \rho\|_{TV} \le  \Bigl|1- {1 \over \rho(B_R)}
\Bigr| + \rho(B_R^c) 
\\
&\le& \left( 1 + \frac1{\rho(B_R)}\right) \,  \rho(B_R^c) \le 3 \, {a \over
R^k}, 
\eean
since then   $ \rho(B_R) \geq 1 - \frac a{R^k} \ge \frac12$. 

As a consequence, for any $\delta \leq 1$ and $a \ge 1$,  choosing $R$
such that $3 \, {a/R^k} = \delta/2$ in the two preceding estimates, we may 
cover $\BB\PPP_{k,a}(E)$ by $\NN_a(\delta) = \NN(R,\delta/2)$ balls of radius $\delta$ in
$W_1$ distance, with 
$$
{1 \over \delta} \le \NN_a(\delta) \le  \Bigl(  C_1  a^{1 \over k} \delta^{-1-
{1 \over k}}  \Bigr)^{C_2\,  a ^{d \over k} \delta^{- d - {d \over k}} }.
$$
The above lower bound on $\NN_a(\delta)$ is straightforwardly obtained by
considering balls centered on Dirac masses
distributed on a line.  In the sequel, we shall often use the shortcut $\NN = \NN_a(\delta)$.
Let us then introduce a covering family $\omega^\delta_i \subset \BB\PPP_{k,a}(E)$, $1
\le i \le \NN_a(\delta)$, such that 
$$
\sup_{\rho,\eta \in \omega^\delta_i} W_1(\rho,\eta) \le 2  \delta, \quad
\omega^\delta_i \cap \omega_j^\delta = \emptyset \,\, \hbox{if} \,\, i \not=j,
\quad \BB\PPP_{k,a}(E) = \bigcup_{i=1}^{\NN_a(\delta)} \omega^\delta_i,
$$
as well as the masses and centers of mass
$$
\alpha^\delta_i := \int_{\omega^\delta_i} \pi, \quad f^{\delta}_i := 
\frac1{\alpha^{\delta}_i }  \int_{\omega^\delta_i} \rho \, \pi (d\rho).
$$
We also denote $\omega_0^\delta := \left[ \BB\PPP_{k,a}(E)
\right]^c$ and $\alpha^\delta_0 := \int_{\omega^\delta_0} \pi$, so that
$\sum_{i=0}^\NN \alpha^\delta_i = 1$. 
Denoting $\ZZ := \{i =1, ..., \NN_a(\delta); \; \alpha^{\delta}_i \ge
\NN_a(\delta)^{-2} \} $,  we finally define
$$
\pi^\delta :=  \sum_{i=1}^{\NN_a(\delta)}  \beta^{\delta}_i \delta_{f^{\delta}_i}, \quad
\hbox{with}\quad \beta^{\delta}_i := {\alpha^{\delta}_i  \over \sum_{j \in \ZZ} \alpha^{\delta}_j}   \,\, \hbox{if} \,\, i \in \ZZ
\; \text{and } \beta^{\delta}_i := 0 \,\, \hbox{if} \,\, i \notin \ZZ.
$$
Remark that by our moment assumption
$$
 \alpha_0^\delta \leq \int_{ (\omega_0^\delta)^c} \pi(d\rho) \le \int_{ \PPP(E) } {M_k(\rho) \over a} \,  \pi(d\rho)  = {M_k(\pi_1) \over a} .
  $$
Since $\sum_{i \notin \ZZ, i \ge 1} \alpha_i^\delta  \leq \NN^{-1} \leq \delta$, we
necessarily have $\ZZ \not= \emptyset$ if $\delta + \frac{M_k(\pi_1)}a \le {1 \over 2} <1$, an assumption that we will make in the sequel.
We fix now the value of $a$ to be so that
$$
\delta = \frac{M_k(\pi_1)}a .
$$
As we shall see that will lead to the optimal inequality. With that particular choice, the condition above simply writes $\delta \le \frac14$, and the upper bound on $\NN$ may be rewritten
\beqn \label{eq:BoundNN2}
\NN(\delta) := \NN_a(\delta) \le  \Bigl(  C_1 \,K \,  \delta^{-1-
{2 \over k}}  \Bigr)^{C_2\,  K^d \delta^{- d\left(1 + \frac2k \right)} }.
\eeqn

In that case, we have
\beqn \label{estim:ZZ}
 \sum_{j \in \ZZ, j \ge 0} \alpha^\delta_j \leq 2 \, \delta , \qquad 1 \ge \sum_{j \in \ZZ} \alpha^\delta_j 
 \ge 1 - 2 \, \delta \ge \frac12.
\eeqn
Now, by convexity of the Fisher information 
$$ I(f^{\delta}_i) \leq  \frac1{\alpha^{\delta}_i }  \int_{\omega^\delta_i}
I(\rho) \, \pi (d\rho),
$$
which in turns implies that 
$$\II( \pi^\delta) =  \sum_{i=1}^{\NN_a(\delta)}  \beta^{\delta}_i I(f^{\delta}_i) \le 
{1  \over \sum_{j \in \ZZ} \alpha^{\delta}_j}  \sum_{i \in \ZZ} 
\int_{\omega^\delta_i} I(\rho) \, \pi (d\rho)  \le  2\, \II(\pi)  . 
$$
Similarly, for the moment of order $k$ :
$$
M_k(\pi^\delta_1) =  \sum_{i=1}^{\NN_a(\delta)}  \beta^{\delta}_i
M_k(f^{\delta}_i) \leq {1  \over \sum_{j \in \ZZ} \alpha^{\delta}_j}  \sum_{i
\in \ZZ} 
\int_{\omega^\delta_i} M_k(\rho) \, \pi (d\rho)  \le  2\, M_k(\pi_1).
$$
In order to prove \eqref{estim:FisherVsEntropGal2}, we introduce the splitting 
\bear\label{estim:HpijHHpi2}
 | H(\pi_{j}) - \HH(\pi) | &\le&  | H(\pi_{j}) - H(\pi^\delta_{j}) | 
 \\&& \nonumber
 + | H(\pi^\delta_{j}) -  \HH(\pi^\delta) |  +   | \HH(\pi^\delta) - \HH(\pi) | ,
\eear
where we have written  $\pi^\delta_{j} := (\pi^\delta)_j$. We now estimate each term separately.

\medskip\noindent
{\sl Step 3. } 
On the one hand, defining $T^{\delta} : \PPP(E) \to \{f^{\delta}_0, ...,
f^{\delta}_\NN \}$, $T^{\delta}(\rho) = f^{\delta}_i$ if $\rho \in
\omega^{\delta}_i$, 
$T^{\delta}(\rho) = f^{\delta}_0 = \delta_0$ if $\rho \in \omega^\delta_0$ and
$\beta^{\delta}_0 := 0$, 
we compute
\bean
\WW_1(\pi,\pi^\delta) 
&\le& \WW_1\Bigl(\pi, \sum_{i=0}^\NN \alpha_i^{\delta} \, \delta_{f^{\delta}_i} \Bigr)
 +  \WW_1 \Bigl( \sum_{i=0}^\NN \alpha_i^{\delta} \, \delta_{f^{\delta}_i},  \sum_{i=1}^\NN \beta_i^{\delta} \, \delta_{f^{\delta}_i} \Bigr) 
\\
&\le& \int_{\PPP(E) \times \PPP(E)} W_1(\rho,\eta) (Id \otimes T^{\delta}) \sharp \pi
+ \Bigl\|  \sum_{i=0}^\NN (\alpha_i^{\delta} - \beta_i^{\delta}) \, \delta_{f^{\delta}_i} \Bigr\|_{TV} 
\\
&\le& \int_{\PPP(E)  } W_1(\rho,T^{\delta}(\rho)) \, \pi (d\rho)  +   
\sum_{i=1}^\NN |\alpha_i^{\delta} - \beta_i^{\delta}| +  |\alpha_0^{\delta} | 
\\
&\le&  \hspace{1.5cm} \delta + \frac{M_1(\pi)}a \qquad + \qquad 6 \delta
\hspace{1cm} \le 8  \delta,
\eean 
where we have used several times estimation~\eqref{estim:ZZ}, in particular in order to get the inequality
$$
\sum_{i=1}^\NN |\alpha_i^{\delta} - \beta_i^{\delta}|  + \alpha_0^\delta   = 
\left( 1 - \frac1{\sum_{i \in \ZZ} \alpha_i^{\delta}}\right) \sum_{i \in \ZZ}
\alpha_i^{\delta} + \sum_{i \notin \ZZ} \alpha_i^{\delta} 
\le  3 \sum_{i \notin \ZZ} \alpha_i^{\delta}.
$$
Using   lemma~\ref{lem:ComparDistPPE} and the bound on $M_k(\pi^\delta_1)$, we
obtain a bound on $\WW_2(\pi^\delta,\pi)$ as follows (we recall that the constant $C$ that appears is numerical : $C = 2^{3/2}$)
$$
 \WW_2 (\pi^\delta,\pi) \le C\,  M_k(\pi_1)^{1/k} \,  \WW_1
(\pi^\delta,\pi)^{1/2 - 1/k} \leq 4 \, C \,  K  \delta^{1/2 - 1/k}. 
$$
Now, we use the HWI inequality on $\PPP(E)$ stated in Proposition~\ref{prop:HWI_PPE} and we bound  the
first term in \eqref{estim:HpijHHpi2} by 
\bean
 | \HH(\pi^\delta) - \HH(\pi) | 
 &\le& \left[ \sqrt{\II(\pi^\delta)} + \sqrt{ \II(\pi)} \right] \, \WW_2
(\pi^\delta , \pi) \leq 2 \, K \, \WW_2 (\pi^\delta , \pi) , 
 \eean
and the third term in \eqref{estim:HpijHHpi2} very similarly 
\bean
 | H(\pi^\delta_{j}) - \HH(\pi_{j}) | 
 &\le& \left[ \sqrt{I(\pi^\delta_{j})}+ \sqrt{ I(\pi_{j})} \right] \, W_2 (\pi^\delta_{ j},
\pi_{j}) 
 \\
 & \le &  \left[ \sqrt{\II(\pi^\delta)}+ \sqrt{ \II(\pi)} \right]  \, \WW_2 (\pi^\delta ,
\pi) \leq 2\, K\,  \WW_2 (\pi^\delta , \pi),
\eean
where we have used the properties~\eqref{eq:Rob&RuelleFisher} of the level 3
Fisher information and Lemma~\ref{lem:WassNleqW1} in order to bound $W_2$ by
$\WW_2$. All together, we have proved 
\beqn\label{ineq:estimThHcvgceDFHS}
 | \HH(\pi^\delta) - \HH(\pi) | + 
 | H(\pi^\delta_{j}) - \HH(\pi_{j}) | 
 \le C \,K^2\,  \delta^{1/2 - 1/k},
 \eeqn
 for some numerical constant $C \le 2^6$. 

\medskip\noindent
{\sl Step 4. } We estimate the second term in \eqref{estim:HpijHHpi2}. 
Using that $\pi_j^\delta = \beta^\delta_1 \, (f_1^\delta)^{\otimes j} + ... +
\beta_\NN^\delta (f_\NN^\delta)^{\otimes j}$, we write
\bean
H(\pi_j^\delta) 
&=& {1 \over j}\int_{E^j} \pi^\delta_j \log \pi^\delta_j 
\\
&=& \sum_{i=1}^\NN \beta^\delta_i \, H(f^\delta_i) + {1 \over j} \int_{E^j} \pi^\delta_j \, \Lambda \left( {\beta^\delta_1 \, (f_1^\delta)^{\otimes j} \over  \pi^\delta_j}, ... , {\beta^\delta_\NN \, (f_\NN^\delta)^{\otimes j} \over  \pi^\delta_j}\right),
\eean
with  $\Lambda : \{ U = (u_i) \in \R^\NN_+, \; \sum_i u_i =1\} \to \R$ defined
by 
$$
\Lambda(U) :=  u_1 \, \log  \left( {\beta^\delta_1 \over u_1}\right)  + ... +
u_\NN \, \log  \left( {\beta^\delta_\NN \over  u_\NN}\right) . 
$$
Observing that $\Lambda$ is in fact (the opposite of ) a discrete relative
entropy, we have for any $U \in \R_+^N$ with $\sum_i u_i = 1$  
$$
- \log (\NN^2) \le  \log (\min \beta^\delta_i)  \le \Lambda(U) \le 0,
$$
we deduce 
$$
| H(\pi^\delta_j)   -  \HH(\pi^\delta) | \le {2 \over j}\, \log
\NN_a(\delta). 
$$

\medskip\noindent 
{\sl Step 5. } All in all, observing that thanks to~\eqref{eq:BoundNN2}
$$
\log \NN(\delta) \le C \, K^d
\, \delta^{- d\left(1+\frac 2k \right)  } \left[ 1 + \ln K - \ln \delta \right] ,
$$
 we have
$$
C^{-1} \, | H(\pi_j) - \HH(\pi) | 
  \le    K^2 \delta ^{1/2 - 1/k} 
+ {1 \over j} \, {K^d \over \delta^{d \left(1+\frac 2k
\right)}} \left[1 + \ln ( K \delta^{-1}) \right].
$$
We can now (almost) optimize by choosing $\delta = j^{-r}$, with 
$r^{-1} := \frac12  - \frac1k +d \, (1 + \frac2k) $ we obtain 
$$
C^{-1} \, | H(\pi_j) - \HH(\pi) |   \le   K^{\max(2,d)} \frac{\ln(K j) }{j^\gamma}
$$
for the integers $j \ge 4^{1/r}$ so that the  the condition $\delta  \leq \frac14$ is fulfilled (in order to ensures that $\ZZ \neq \emptyset$). But it can be checked that for $k \in [2, +\infty)$, $d \le \frac1r \le 2d$. So that the previous condition  on $j$ is fulfilled for $j \ge 4^{2d}$. 

\qed



\bibliographystyle{acm}


 \signmh\signsm
\end{document}